\documentclass[11pt, leqno]{amsart}
\usepackage{amsthm, amsmath}
\usepackage{amssymb} 
\usepackage{amscd}
\usepackage[curve,matrix,arrow,frame,tips]{xy}
\usepackage{verbatim}
\newtheorem{thm}{Theorem}[section]
\newtheorem{prop}[thm]{Proposition}
\newtheorem{lemma}[thm]{Lemma}

\newtheorem{Definition}[thm]{Definition}
\newtheorem{Remarknumb}[thm]{Remark}
\newtheorem{Remarks}[thm]{Remarks}
\newtheorem{Remark}[thm]{Remark}

\newtheorem{conjecture}[thm]{Conjecture}
\newtheorem{cor}[thm]{Corollary}

\newcounter{ex}[section]

\newcommand{\cal}{\mathcal}

\newcommand{\E}{{\mathcal E}}

\newcommand{\C}{{\bf C}}

\newcommand{\Q}{{\bf Q}}

\newcommand{\ep}{\epsilon}

\newcommand{\Hom}{{\rm Hom}}
\newcommand{\X}{{\mathcal X}}

\newcommand{\Gm}{{{\bf G}_m}}
\newcommand{\Z}{{\bf Z}}
\newcommand{\F}{{\mathcal F}}
\newcommand{\ti}{\tilde}
\newcommand{\Spec}{{\rm Spec }\, }
 \renewcommand{\O}{{\mathcal O}}

\newcommand{\und}{\underline}

\newcommand{\gfr}{{\mathfrak g}}

\renewcommand{\L}{{\mathcal L}}

\newcommand{\Res}{{\rm Res}}

\newcommand{\Tr}{{\rm Tr}}

\def\thfill{\null\nobreak\hfill}

\def\endproof{\thfill\vbox{\hrule
  \hbox{\vrule\hbox to 5pt{\vbox to 5pt{\vfil}\hfil}\vrule}\hrule}}

\newcommand{\tiuS}{\und{\ti S}}

\newcommand{\aff}{{\rm aff}}

\newcommand{\Gal}{{\rm Gal}}

\begin{document}

\title[twisted loop groups]{Twisted loop groups and their affine flag varieties}
\author[G. Pappas]{G. Pappas*}
\thanks{*Partially supported by NSF
grants \# DMS05-01049 and  \# DMS01-11298 (via the Institute for Advanced Study)}
\address{Dept. of
Mathematics\\
Michigan State
University\\
E. Lansing\\
MI 48824-1027\\
USA}
\email{pappas@math.msu.edu}
\author[M. Rapoport]{M. Rapoport}
\address{Mathematisches Institut der Universit\"at Bonn,  
Beringstrasse 1\\ 53115 Bonn\\ Germany.}
\email{rapoport@math.uni-bonn.de}

\date{\today}
 
\maketitle

\section*{Introduction}

Loop groups are familiar objects in several branches of mathematics. Let us mention here three variants.

The first variant is differential-geometric in nature. One starts with a Lie group $G$ (e.g., a compact
Lie group or its complexification). The associated loop group is then the group of ($C^0$-, or $C^1$-, or $C^\infty$-)
maps of $S^1$ into $G$, cf.~[P-S] and the literature cited there. A twisted version arises from an automorphism $\alpha$ of $G$. The associated
twisted loop group is the group of maps $\gamma : {\bf R} \to G$ such that
\begin{equation*}
\gamma (\theta + 2\pi)\ =\ \alpha(\gamma(\theta))\ .
\end{equation*}

The second variant is algebraic and arises in the context of Kac-Moody algebras. Here one constructs an 
infinite-dimensional algebraic group variety with Lie algebra equal or closely related to a given Kac-Moody algebra.
(This statement is an oversimplification and the situation is in fact more complicated: there exist
various constructions  at a formal, a minimal, and a maximal level which produce infinite-dimensional
groups with Lie algebras closely related to the given Kac-Moody Lie algebra,  see [Ma2], also [T2], [T3] 
and the literature cited there). 

The third variant is algebraic-geometric in nature and is our main concern in this paper. Let us recall the
basic definitions in the untwisted case. Let $k$ be a field and let $G_0$ be an algebraic group over $\Spec (k)$.
We consider the functor $LG_0$ on the category of $k$-algebras,
\begin{equation*}
R\mapsto LG_0(R)=G_0(R((t))).
\end{equation*}
Here $R((t))=R[[t]][1/t]$ is the ring of Laurent series with coefficients in $R$. One proves that this functor
is representable by an ind-scheme ($=$ inductive limit of $k$-schemes), called {\it the algebraic loop group} associated
to $G_0$. At the same time, one considers the associated flag varieties, mostly the affine Grassmannian $\mathcal G$
and the affine flag variety $\mathcal F$, the $fpqc$-sheaves associated to the functors,
\begin{equation*}
R\mapsto G_0(R((t)))/G_0(R[[t]]),\  \text{resp.}\ \ 
R\mapsto G_0(R((t)))/\mathcal B(R).
\end{equation*}
Here $\mathcal B (R)$ is the inverse image of a Borel subgroup under the reduction map $G_0(R[[t]])\to G_0(R)$.
These sheaves are also representable by ind-schemes. Important results on the structure of algebraic loop groups
and their associated flag varieties are due to Beauville, Laszlo, Sorger, Faltings, Beilinson, Drinfeld and
Gaitsgory, cf.~[B-L], [BLS],  [L-S], [Fa1], [B-D], [G]. These results have applications in the theory of vector bundles on algebraic curves
[B-L], [L-S], [Fa1], in geometric Langlands theory [B-D], [G], and to local models of Shimura varieties
[Go1], [Go2], [P-R].

The aim of the present paper is to develop a similar theory of twisted algebraic loop groups
and of their associated flag varieties.
As it turns out, this theory gives a geometric interpretation of many aspects of Bruhat-Tits 
theory in the equal characteristic case. 

The basic definition is very simple. Let $G$ be a linear algebraic group over $K=k((t))$ (we always assume
that $G$ is connected reductive). Then $G$ defines its associated algebraic loop group $LG$, which is the
ind-scheme representing the functor
\begin{equation*}
R\mapsto LG(R) = G(R((t))).
\end{equation*}
Note that, since $R((t))$ is a $k((t))$-algebra, the right hand side makes sense. When $G=G_0\otimes_k K$, one recovers the definition in the untwisted case. To any facet $\underline{a}$ in
the Bruhat-Tits building of $G$ over $K$ there is associated a unique smooth group scheme with connected fibers 
$P_{\underline{a}}$ over $k[[t]]$ such that $P_{\underline{a}}(k[[t]])$ is equal to the parahoric subgroup of
$G(K)$ attached to $\underline{a}$. To $P_{\underline{a}}$ corresponds an infinite-dimensional affine group scheme
$L^+P_{\underline{a}}$ over $\Spec (k)$ with
\begin{equation*}
L^+P_{\underline{a}}(R)=P_{\underline{a}}(R[[t]]).
\end{equation*}
The quotient $\mathcal F_{\underline{a}}=LG/L^+P_{\underline{a}}$ (in the sense of $fpqc$-sheaves) is
representable by an ind-scheme over $k$ which we call 
the (partial) affine flag variety associated to $G$ and $\und a$. In the untwisted case, the affine Grassmannian is associated to the special vertex
arising from an isomorphism $G\simeq G_0\otimes_k k((t))$, and the affine flag variety is associated to an alcove in the Bruhat-Tits building containing this vertex. Note that even in the untwisted case, the definition of  $\mathcal F_{\underline{a}}$ is useful, as it re-establishes the natural symmetry attached to the $K$-group scheme $G_0\otimes_kK$ which has no preferred special vertex. 

Our main results may now be formulated as follows.

Let $\pi_1(G)$ denote the algebraic fundamental group of $G$ 
(see [Bo]). The group $\pi_1(G)$  can be
  defined as the quotient $P^\vee(G)/Q^\vee(G)$ of the coweight lattice
by the coroot lattice of $G(K^{\rm sep})$. If
$k$ is algebraically closed, let  $I=\Gal(\overline{K}/K)$ denote the inertia group. Let
$\pi_1(G)_I$ be
the co-invariants under $I$.

\begin{thm}\label{A}  Let $k$ be algebraically closed. The Kottwitz homomorphism ([Ko], see \ref{parahoric}) induces bijections on the set of connected
components,
\begin{equation*}
\pi_0(LG)=\pi_1(G)_I,\ \  \pi_0(\mathcal F_{\underline{a}})=\pi_1(G)_I.
\end{equation*}
\end{thm}
Note that in the untwisted case $G=G_0\otimes_k K$, the action of $I$ is trivial so that $\pi_1(G)_I=\pi_1(G)$
and we recover the result of Beauville, Laszlo and Sorger [BLS], Lemma~1.2., resp.~of Beilinson and Drinfeld [B-D], Prop.~4.5.4, which is analogous to the corresponding result in topology.

\begin{thm}\label{B}  Suppose that $G$ is semi-simple and splits over a tamely ramified extension of $K$ and that the order of
the fundamental group of the derived group $\pi_1(G_{\rm der})$ is prime to the characteristic of $k$. Then the ind-schemes $LG$ and
$\mathcal F_{\underline{a}}$ are reduced.
\end{thm}

In the case of $SL_n$ this is due to Beauville and Laszlo [B-L], and in the  case of a general split group and when $k$ has characteristic $0$, to Beilinson and
Drinfeld [B-D], Thm.~4.5.1.  In characteristic $p>0$, this is due to Laszlo and Sorger [LS], and to Faltings [Fa1], when $G$ is split and simply connected.

We next consider the Schubert varieties contained in a partial affine flag variety $\mathcal F_{\underline{a}}$. By definition, these are the reduced closures of orbits
of $L^+P_{\underline {a}}$. They are finite-dimensional projective varieties. 

\begin{thm}\label{C}  Suppose that $G$ splits over a tamely ramified extension of $K$
and that the order of
the fundamental group of the derived group $\pi_1(G_{\rm der})$ is prime to the characteristic of $k$. Then all Schubert varieties in
$\mathcal F_{\underline{a}}$ are normal and have only rational singularities. In positive characteristic, they
are Frobenius-split.
\end{thm}

In characteristic zero, normality of Schubert varieties in the context of Kac-Moody partial 
affine flag varieties is due
to  Kumar [Ku1], Littelmann [Li] and Mathieu [Ma1]. By an argument of Faltings these 
coincide with the Schubert varieties in this paper.
In positive characteristic the above theorem is due to
Faltings [Fa1] in the split simply connected case. (The normality of Schubert varieties in the case of $SL_n$ is proved in [P-R1]). Normality of Schubert varieties 
 in characteristic $p$ has also been shown by 
Mathieu [Ma1] and Littelmann [Li] (in  the context of Kac-Moody theory). Their definition of Schubert varieties is a priori different.
They are given using an embedding into the infinite dimensional projective space associated to a 
highest-weight representation of the Kac-Moody algebra. In contrast, we are essentially defining Schubert varieties as 
varieties of certain lattices and so [Ma1] and  [Li] does not imply our result.
Nevertheless, it is a corollary of  Theorem \ref{C} combined 
with the  results of Mathieu and Littelmann 
that the two notions of Schubert varieties coincide (see Section \ref{CompAtp}).
From this, it also follows a posteriori that (for an absolutely simple, simply connected 
$G$ that splits over a tame extension) our partial affine flag varieties are 
isomorphic to the ones defined in the Kac-Moody theory.

Let us now comment on our proofs of the above theorems.  The proofs of theorems \ref{B} and \ref{C} are closely intertwined.  The proofs  go  along the lines of the proof of the theorem corresponding   to Theorem \ref{C} of Faltings [Fa1]
(for more details on  this proof, cf. [Go2]). The main method is to lift the situation in positive
characteristic to characteristic zero and compare the lifted situation with the Kac-Moody setting. It is therefore somewhat reminiscent of the method introduced by Kazhdan to relate the representation theory of a group over a local field of characteristic $p$ to that of a group over a local field of unequal characteristic. This also
explains why we have to exclude groups which split over a wild extension (although we strongly believe that 
appropriate versions of Theorems \ref{B} and \ref{C} also hold in this case). We note that in the previous proofs of theorem \ref{B} (in [BL] for the case of $SL_n$, in [LS] for a split semi-simple, simply connected group in characteristic $0$, in [Fa1, Go2] for   a semi-simple split simply connected group in characteristic $p$)  all use  the big cell.  Our   argument,  which we learned in Faltings' course  [Fa2] is different.  Finally, the proof of Theorem \ref{A}  generalizes the proof of
the corresponding theorem in [B-D], using the construction of the Kottwitz map [Ko].

We would like to stress that  interesting problems arise when one tries to 
generalize to the setting of twisted loop groups classical results in the untwisted case. The theory 
of perverse sheaves on their partial flag varieties has still to be worked out. Also, we expect a  theory of 
bundles on algebraic curves for corresponding non-constant algebraic groups, of their Picard groups and their theta functions.   

Our main motivation for developing the theory of twisted loop groups is its application to   the theory of local models of
Shimura varieties. These are  schemes defined in linear algebra terms that describe the \'etale local structure of integral models for 
Shimura varieties and other moduli spaces [R]. Consider the group $G$ of unitary similitudes corresponding to a ramified quadratic
extension of $K$. Then $G$ is not of the form $G=G_0\otimes_k K$. Generalizing the method of G\"ortz [Go1], we
wish to embed the special fiber of the naive  local model associated to a ramified group of unitary
similitudes over a $p$-adic local field into the affine flag variety $\mathcal F$ of $G$. At this point we require knowledge of the
structure of $\mathcal F$ which we then can in turn use to deduce structure results on the local model. This is the way in which we will use the results of this paper. To give the taste of this application we  explain in the last section   this method in a representative simple case. 

This method hinges on a coherence conjecture on the dimensions of the spaces of global sections of the natural ample invertible sheaves on partial flag varieties attached to a fixed group $G$ over $K$.  This conjecture, which may be stated in purely combinatorial terms, seems to us of independent interest. We show that in the cases of $SL_n$ and $Sp_{2n}$ it follows from previous results [Go1, Go2, PR1, PR2] on local models. The main point of the application we have in mind here
is to turn the logic around and deduce facts on local models from this coherence conjecture. 

We now explain the lay-out of the paper. In section \ref{twiloop} we give the construction 
of twisted loop groups and their quotients. In section \ref{affineflags} we recall the definition of parahoric subgroups and introduce the (generalized) affine flag varieties. In section \ref{tori} we make these notions explicit in the case of tori. In section \ref{unitary} we discuss in some detail the case of a unitary group and relate their partial affine flag varieties to the lattice model. As explained above, this example plays a central role in the applications of our results. Section \ref{kottwitzmap} is devoted to the proof of Theorem \ref{A}. In section \ref{reduced} we discuss Theorem \ref{B} and reduce its proof to the case of a simple simply connected   group. In section \ref{tamepar} we explain how to lift the group $G$ and its partial flag varieties to mixed characteristic, i.e., from $k((t))$ to $W((t))$, where $W$ denotes the ring of Witt vectors of $k$. In section \ref{weyl} we introduce 
the Schubert varieties in twisted affine flag varieties and their Demazure resolutions. The proofs of Theorems \ref{B} and \ref{C} are given in section \ref{proofmain}. In section \ref{coherence} we state and discuss the coherence conjecture mentioned above. In the final section, we discuss how we   apply our results to the local models of Shimura varieties  associated to ramified unitary groups. As mentioned above, these applications depend on a positive solution to the coherence conjecture. 

In conclusion, we would like to stress how much we profited from the fact that one of us was able to attend Faltings' course at the University of Bonn. We also thank C. Kaiser and J.-L.~ Waldspurger for helpful discussions and G. Prasad for some comments. Finally, we would like to thank 
the University of Bonn,    Michigan State University, and the Institute for Advanced Study for making our collaboration possible.

\section{General constructions}\label{twiloop}

\setcounter{equation}{0}

Let $k$ be a field. We set $K=k((t))$ for the field of Laurent power
series with indeterminate $t$ and coefficients in $k$. Let $V=k[[t]]\subset K$ be the discretely valued ring 
of power series with coefficients in $k$. 

\subsection{}

Let $X$ be a scheme over $K=k((t))$. We denote by $LX$
the functor from the category of $k$-algebras to that of sets
given by
\begin{equation*}
R\mapsto LX(R):=X(\Spec(R((t))))\ .
\end{equation*}
If ${\cal X}$ is a scheme over $V$, we denote by $L^+{\cal X}$ the
functor from the category of $k$-algebras to that of sets
given by
\begin{equation*}
R\mapsto L^+{\cal X}(R)={\cal X}(\Spec(R[[t]]))\ .
\end{equation*}
The functors $LX$, $L^+{\cal X}$ give sheaves of sets for the fpqc topology
on $k$-algebras. In what follows, for simplicity, we will call
such functors ``$k$-spaces".
Recall that a $k$-space is called an ind-scheme (resp. a strict ind-scheme) if it is the inductive limit
of the functors associated to a directed family
of $k$-schemes (resp., via transition morphisms which are closed embeddings).
A group ind-scheme is an ind-scheme which is a group object in the category
of ind-schemes.

If $\X = \Spec (A)$ is affine of finite type, $L^+\X$ is represented by an affine scheme. If $\X = {\bf A}^r_V$ is the affine space of dimension $r$ over $\Spec (V)$, then $L^+\X$ is the infinite-dimensional affine space 
$L^+\X = \underset{i = 0}{\overset{\infty}{\prod}} ({\bf A}^r)$, via
\begin{equation*}
L^+\X (R) = \Hom_{k[[t]]} (k[[t]][T_1, \dots , T_r], R[[t]])
= R[[t]]^r = \prod^{\infty}_{i = 0} R^r =\prod^{\infty}_{i = 0} {\bf A}^r (R)\ .
\end{equation*}
Let $\X$ be the closed subscheme of ${\bf A}^r_V$ defined by the vanishing of a polynomial $f$  in $k[[t]][T_1, \dots , T_r]$. Then $L^+\X (R)$ is the subset of 
$L^+ {\bf A}^r(R)$ of $k[[t]]$-algebra homomorphisms $k[[t]][T_1,\ldots,T_r]\to R[[t]]$ which factor through $k[[t]][T_1,\ldots,T_r]/(f)$.     Write  
$f = \sum a_{\underline{n}} (t)\cdot T^{\underline{n}}$.  Then $L^+\X$ is defined by
the vanishing of the infinite set of polynomials in the variables $X_{1, \bullet}, X_{2, \bullet}, \dots, X_{r, \bullet}$
with $\bullet = 0, 1, 2, \dots ,$
\begin{equation*}
\sum_{\underline{n}} C_i^{( \underline{n})} (a_{\underline{n}, \bullet}; X_{1, \bullet}, ..., X_{r, \bullet}) = 0,\ i = 0,1, \dots \ .
\end{equation*}
Here the polynomials $C_i^{( n_1, \dots , n_r)}$ in the indeterminates $a_0, \dots , a_i, X_{1, 0}, \dots , X_{1, i}, \dots , X_{r, 0}, \dots , X_{r, i}$ are defined by the identity
\begin{equation*} 
(\sum^\infty_{i = 0} a_i\cdot t^i) \cdot (\sum^\infty_{i = 0} X_{1, i}\cdot t^i)^{n_1}\cdot\cdot\cdot (\sum^\infty_{i = 0} X_{r, i}\cdot t^i)^{n_r} = \sum^\infty_{i = 0} C_i^{( n_1, ..., n_r)}\cdot t^i\ .
\end{equation*}
Note that if $\X$ is an affine $k$-scheme, which we regard as an affine scheme over $V$
via $V\to k$ which sends $t$ to $0$,  then $L^+ \X =\emptyset$.
 If $X$ is an affine $K$-scheme, then $L^+ X =\emptyset$. On the other hand, in this case $LX$ is 
represented by a strict ind-scheme. Finally, we mention that $L(X\times_{{\rm Spec}\, K}Y)
=LX\times_{{\rm Spec}\, k}LY$, and similarly for $L^+$.
\begin{Definition}\label{1.1} {\rm Let $G$ be a linear algebraic group over $K$. The {\sl loop group}
associated to $G$ is the ind-group scheme $LG$ over $\Spec(k)$.}
\end{Definition}

In the case $G=GL_n$, the fact that $LG$ is  an ind-group scheme is well-known, comp. [B-L],
Prop. 1.2. 
In fact, if $G$ is of the form $G=G_0\times_{\Spec(k)}\Spec(K)$ for a linear algebraic group
$G_0$ over $k$, then the notion of a loop group is in [B-D],  [B-L], [Fa1].

\subsubsection{} We now list some easy functoriality properties of this construction. 

a) If $k'$ is a $k$-field extension then we have an isomorphism of ind-schemes over $k'$
\begin{equation}
LG\times_{{\rm Spec}(k)}{{\rm Spec}(k')}\simeq L(G\times_{{\rm Spec}(k((t)))}{{\rm Spec}(k'((t)))})\ .
\end{equation}

b) Assume that $K'/K$ is a finite extension of $K$; then a choice of uniformizer $u$ for $K'$
allows us to write $K'=k'((u))$.
If $G={\rm Res}_{K'/K}H$ for some linear algebraic group $H$ over $K'$, then the choice of $u$
gives an isomorphism of ind-schemes over $k'$,
\begin{equation}\label{weilre1}
LG \simeq {\rm Res}_{k'/k}(LH).
\end{equation}

In particular, if $k'=k$ then we have
\begin{equation}\label{weilre2}
L({\rm Res}_{K'/K}H)\simeq LH\ .
\end{equation}
Indeed, for a $k$-algebra $R$, we have
\begin{eqnarray*}
(LG)(R)=G(R((t)))=H(R((t))\otimes_{k((t))}k'((u)))
=H((R\otimes_kk')((u)))\\=LH(R\otimes_k k')={\rm Res}_{k'/k}(LH)(R).\ \ \ \ \ \
\end{eqnarray*}
(Note that $R((t))\otimes_{k((t))}k'((u))
=(R\otimes_kk')((u))$ uses that $k'((u))/k((t))$ is finite.)

\subsubsection{}\label{lie}  We will denote by $Lie(LG)$ the Lie algebra functor of the group functor of $LG$ over $k$
(see for example [D-G], II, \S 4, 1). By definition
\begin{equation}\label{lie1}
Lie(LG)(R)={\rm ker}(G(R[\epsilon]((t))\to G(R((t)))).
\end{equation}
The points $Lie(LG)(k)$ form a $k$-vector space.  If $R$ is a $k$-algebra which is finite dimensional over $k$  and $J\subset R$
is an ideal of square $0$ then there is an exact sequence
\begin{equation}\label{lie2}
0\to Lie(LG)(k)\otimes_k J\to LG(R)\to LG(R/J)\ .
\end{equation}

\subsubsection{} \label{steinberg} Suppose that $k$ is algebraically closed and $G$ is a connected reductive group over $K=k((t))$.
By a theorem of Steinberg ([St2], extended by Borel and Springer to the case of a non-perfect field, cf. [Se2],  III, 2.3, Remark 1 after Theorem 1$'$),
every $G$-torsor over $K$ is trivial. Let $R$ be an Artinian $k$-algebra with algebraically closed residue field
$k'$;
then $R((t))$ is a complete local ring with residue field $k'((t))$. Hence, $R((t))$ is   henselian and by the above, every
$G$-torsor over $R((t))$ is trivial. It follows that if
\begin{equation*}
1\to G'\to G\to G''\to 1
\end{equation*}
is an exact sequence of connected reductive groups over $K$ and $R$ is as above, the sequence
\begin{equation*}
1\to LG'(R)\to LG(R)\to LG''(R)\to 1
\end{equation*}
is exact.

\subsection{}

Now let $P$ be a flat affine group scheme of finite type over $\Spec(V)$, with associated affine group scheme $L^+P$ over $k$. 
We wish to form the quotient of the generic fiber $LP_\eta$ by $L^+P$. 

\begin{prop}\label{line}
Suppose that $Q\subset G$ are both flat affine group schemes of finite type over $\Spec(V)$
and that $Q$ is a proper closed subgroup scheme of $G$. Then there 
exists a finitely generated free $V$-module $\Lambda\simeq V^n$ and a direct summand $\Lambda'\subset \Lambda$
of $\Lambda$ of $V$-rank $1$ such that:

a) There is a representation $\rho: G\to {\rm GL}(\Lambda)$ that identifies $G$ with
a closed subgroup scheme of ${\rm GL}(\Lambda)$,

b) The representation $\rho$ identifies $Q$ with the closed subgroup scheme of $G$
that normalizes the $V$-line  $\Lambda'$.
\end{prop}

\begin{Proof}
This follows the argument of the proof of [D-G], II \S 2, 3.3. Let $A(G)$, $A(Q)$ be the $V$-algebras
of regular functions on $G$ and $Q$ respectively. Since $Q$ and $G$ are $V$-flat, these are $V$-torsion free.
Denote by $I$ the ideal of definition of $Q$ in $G$ so that $A(Q)=A(G)/I$. The group scheme $G$ acts on $A(G)$ via the 
``regular representation". Using [Se1], \S 1, Prop. 2, and the fact that $A(G)$ is $V$-flat, we see that there is a finitely generated free $V$-submodule $L$ of
$A(G)$ with the following properties:

i) $G$ acts on $L$, 

ii) $A(G)/L$ is $V$-torsion free,

iii) $L$ generates the $V$-algebra $A(G)$ and is such that $I\cap L$ generates the ideal $I$.

 Note that since $Q$ is flat, $I$ is also $V$-cotorsion free in $A(G)$. 
Set $L'=I\cap L$. Then we can see by following the proof of [D-G], II \S 2, 3.3 that   $\Lambda=\wedge^{{\rm rank}(L')}L$ and 
$\Lambda'= \wedge^{{\rm rank}(L')}L'=\det(L')$ satisfy the requirements of the proposition.\endproof
\end{Proof}

\begin{prop}\label{quasiaffine}
Suppose that $P$ is a flat affine group scheme of finite type over $\Spec(V)$. Then

a) $P$ is linear, i.e., it is a closed subgroup scheme of the group scheme ${\rm GL}_n$ over $\Spec(V)$,
for some $n$.

b) There exists a group scheme closed immersion $P\subset {\rm GL}_n\times \Gm$ such that  the quotient fppf sheaf $({\rm GL}_n\times\Gm)/P$ is representable 
by a quasi-affine scheme over $\Spec(V)$.
\end{prop}

\begin{Proof}
To obtain (a) we  apply Proposition \ref{line} to $Q=\{1\}$ and $G=P$. Now let us discuss the proof of (b).
We apply Proposition \ref{line} to $Q=P$ and $G={\rm GL}_n$.
Let $\chi: P\to \Gm={\rm GL}(\Lambda')$ be the character of $P$ giving the action on the $V$-line
$\Lambda'$ of Proposition \ref{line}. Now consider the embedding $i:P \to G':={\rm GL}_n\times\Gm$  given
by $p\mapsto (p, \chi^{-1}(p))$
 and the representation $\rho': G'\to {\rm GL}(\Lambda\otimes_V \Lambda') $
given by $\rho'((g,\lambda))=\rho(g)\otimes \lambda$. Let  $v\in \Lambda'$ be a generator of the $V$-line $\Lambda'$, i.e, 
such that $\Lambda'=V\cdot v$. Then the subgroup scheme $P$ is identified with
the stabilizer of $v\in \Lambda$. Hence, the quotient $({\rm GL}_n\times\Gm)/P$
(which is representable by a $V$-scheme by a general result
of [An],  Th. 4.C) is identified with the ${\rm GL}_n\times\Gm$-orbit
of $v\in \Lambda$. A standard argument (cf. [D-G],  proof of II \S 5, Prop. 3.1)
shows that this  orbit  is open in its Zariski closure and hence
it is quasi-affine.\endproof
\end{Proof}

\subsection{}
Let $P$ be a flat affine group scheme of finite type over $V$. Let $P_\eta$ denote the generic fiber of $P$
(a group scheme over $K$). We consider the quotient fpqc sheaf over $\Spec(k)$
\begin{equation}
\F_P:=LP_\eta/L^+P\ .
\end{equation}
By definition, this is the fpqc sheaf associated to the presheaf which to a $k$-algebra $R$ associates
the quotient $P(R((t)))/P(R[[t]])$. 
\medskip

\subsubsection{} \label{afgraexample} When $P={\rm GL}_n$, the fpqc sheaf $\F_P$ is (represented by) the affine Grassmanian
of lattices in the vector space $K^n$.  More precisely, for any $k$-algebra $R$, a lattice 
in $R((t))^n$ is a sub-$R[[t]]$-module $\L$ of $R((t))^n$ which is 
projective of rank $n$ 
and is such that $\L\otimes_{R[[t]]}R((t))=R((t))^n$. Consider the functor $\F_n$ which to 
a $k$-algebra $R$ associates the set of lattices in $R((t))^n$. Then there is a natural isomorphism
$\F_{{\rm GL}_n}\simeq \F_n$. It is well-known [BL], Prop. 2.3, that $\F_n$ is represented by an ind-$k$-scheme which is ind-proper over $k$
(i.e.,  given as an increasing union of proper schemes over $k$).  In this case ($P={\rm GL}_n$), the quotient morphism
\begin{equation*}
L{\rm GL}_n\to LGL_n/L^+GL_n=\F_{{\rm GL}_n}
\end{equation*}
admits   sections locally for the Zariski topology, cf. [B-L], Thm. 2.5. In the general case we have the following statement.

\begin{thm}\label{represent}
Let $P$ be an affine group scheme which is smooth  over $V$. Then the fpqc sheaf
$\F_P$ is represented by an ind-$k$-scheme of ind-finite type over $k$. The quotient
morphism $LP_\eta\to \F_P$ admits sections locally for the \'etale topology (i.e., $L^+P$-equivariant
isomorphisms
  $\Spec(R)\times_{\F_P}LP_\eta\simeq \Spec(R)\times_{\Spec(k)} L^+P$ for each point of
  $\F_P$ with values in a strictly henselian ring $R$.)
\end{thm}

\begin{Proof}
The result is shown in [B-D], Theorem 4.5.1, or [G], Appendix, in the case when $P$ is obtained by base change
$P=P_0\otimes_k V$, with $P_0$ an affine smooth group scheme over $k$. In fact, in these
references it is shown that if $P_0$ is reductive, then $\F_P$ (which is then the affine Grassmannian
for $P_0$) is ind-proper over $k$. Here, we will also follow the argument in the proof of [B-D], Theorem 4.5.1. By Proposition \ref{quasiaffine},
there is an embedding $P\to G={\rm GL}_n\times \Gm$ such that the quotient $U:=G/P$ (over $\Spec(V)$)
is quasi-affine, i.e.,  $U$ is an open subscheme of an affine $V$-scheme $Z$. We will consider the
functors $LU$, $LZ$, $L^+U$, $L^+Z$, defined as above. Then $L^+U$, $L^+Z$ are schemes over $k$
and $L^+U$ is an open subscheme of $L^+Z$.  Note that $L^+U$ represents the fpqc quotient $L^+G/L^+P$: 
Actually,  if $R$ is strictly henselian so is $R[[t]]$ and since $P$ is smooth, we have $U(R[[t]])=G(R[[t]])/P(R[[t]])$. 
We can see that, since $Z$ is affine over $V$, $L^+Z$ is a closed subscheme of the ind-scheme $LZ$. Therefore, $L^+U$ is a locally closed sub-ind-scheme of $LZ$. Now consider the natural morphism $f: LG\to LZ$ given by the quotient
$G\to U=G/P$ followed by the open immersion $U\to Z$. Let $Y:=f^{-1}(L^+U)$ be the inverse image of $L^+U$; it is a locally closed sub-ind-scheme of $LG$. Suppose that we have $g\in Y(R)\subset LG(R)=G(R((t)))$. We have $LG(R)/LP_\eta (R)\subset LU(R)$ and after an \'etale cover $R\to R'$,  $f(g)\in L^+U(R)$ comes from an element of  $L^+G(R')$ modulo $L^+P(R')$. Therefore, locally for the \'etale topology, we can write $g$ as a product of an element from $LP_\eta $ with one from $L^+G$. The ind-scheme $Y$ is invariant under right translation by $L^+G$. Hence, and since $LG\to \F_G$ admits a section locally for the Zariski topology, $Y$ is the preimage of a locally closed sub-ind-scheme $Y'\subset LG/L^+G=\F_G$; in fact, if the quotient $G/P$ is affine, then $U=Z$ and $Y'$ is closed in $\F_G$. In any case, $Y'$ is of ind-finite type. Now observe that $LP_\eta $ is contained in $Y$, since it maps to the trivial coset under $f$.  This gives a natural morphism
\begin{equation}\label{LPY}
\pi: LP_\eta\to Y'\ .
\end{equation}
We claim that $\pi$ identifies $Y'$ with the fpqc quotient $\F_P=LP_\eta/L^+P$.
It is clear that $\pi$ induces a map of fpqc sheaves $LP_\eta/L^+P\to Y'$.
To see that this  gives an isomorphism note that $Y\to Y'$ is a $L^+G$-torsor, that the morphism
$f: Y\to L^+U$ is $L^+P$-equivariant and that $LP_\eta=f^{-1}(e\, {\rm mod}\,  L^+P)$.
This combined with the above gives the result. In fact, the proof shows that $\pi$ admits 
a section locally for the \'etale topology.
\endproof
\end{Proof}
\bigskip
\bigskip

\section{Affine flag varieties}\label{affineflags}
\setcounter{equation}{0}

\subsection{}Let $G$ be a (connected) reductive group over the local field $K$.
(Recall, that then, by [SGA3] XXIII 5.6, $G$ splits over the separable closure $K^{\rm sep}$ of $K$.)
Our main examples of group schemes $P$ to which we want to apply Theorem \ref{represent}  are given by the group schemes associated to the parahoric subgroups of $G(K)$
by Bruhat-Tits [B-T], [B-TII]. These group schemes are affine smooth
over $\Spec(V)$ and their generic fibers coincide with $G$. Since these are basic for what
follows we take some space to properly set the definition and notation
by tracing down the corresponding sections in [B-TII]. The corresponding quotients  $\F_P$ are the partial  affine flag
varieties and will be our main concern.

\subsubsection{} Set $K^{\flat}=\bar k((t))$ with $\bar k$ the algebraic closure of $k$.
Let ${\cal B}={\cal B}(G(K))$ be the Bruhat-Tits building of $G(K)$
([B-T], 7.4). Suppose that $\und a$ is a facet
in ${\cal B}$.  This corresponds to a 
facet $\und a^\flat$ of the building ${\cal B}(G(K^{\flat}))$ which is fixed under $\Gamma:=\Gal(K^\flat/K)=\Gal(\bar k/k)$ ([B-TII], Th. 5.1.25).
Now Bruhat-Tits associate to $\und a^\flat$ the ``connected stabilizer" group scheme
$P_{\und a}^\flat$ (see [B-TII], 4.6), comp.~also [Yu], section 7.  This comes with $\Gamma$-descent data
(see [B-TII], 5.1.9) which give the group scheme $P_{\und a}$. The $V$-valued points
$P_{\und a}(V)$ give the corresponding parahoric subgroup of $G(K)$
(see [B-TII], p. 165), and  $P_{\und a}$ is the unique extension to $\Spec (V)$ of $G$ as a smooth affine group scheme with this property, cf. [B-TII], 1.7.

Note that by construction $P_{\und a}(V)=P_{\und a^\flat}(V^\flat)\cap G(K)$.

\subsubsection{}\label{parahoric}
Assume now that $k$ is algebraically closed in which case $K^\flat=K$.
Then we can also obtain the parahoric subgroup $P_{\und a}(V)$
as follows [H-R], Prop. 3 (see also [R]): Consider the corresponding facet $\und a'$ in the Bruhat-Tits building of the adjoint group
${\cal B}(G_{\rm ad}(K))$. Then $P_{\und a}(V)$ is the intersection of the stabilizer of ${\und a'}$
in $G(K)$ with the kernel of the Kottwitz  homomorphism
\begin{equation}
\kappa_G: LG(k)=G(K)\to \pi_1(G)_{I}.
\end{equation}
Here again we denote by $\pi_1(G)$ the $\Gal(K^{\rm sep}/K)$-module
given by the algebraic fundamental group ([Bo]).
Recall $\pi_1(G)$ is the quotient $P^\vee(G)/Q^\vee(G)$ of the coweight lattice
by the coroot lattice of $G(K^{\rm sep})$.
The inertia group $I=\Gal(K^{\rm sep}/K)$ acts on $\pi_1(G)$ and we denote by $\pi_1(G)_I$
the coinvariants under this action. 

Let us quickly recall the construction of the
Kottwitz homomorphism in this case (cf. [Ko] \S 7, where the case is considered, in which $k((t))$ is replaced by the
completion of the maximal unramified extension of a local $p$-adic field.) It proceeds in four steps:
\smallskip

{\sl Step 1.} The case that $G=\Gm$. Then $\pi_1(G)=\Z$ with trivial $I$-action.
Write $f\in LG(k)=k((t))^\times$ in the form $f=t^k\cdot u$, $u\in k[[t]]^\times$. Then set $\kappa(f)=k$.
This definition extends in the obvious way to the case when $G={\rm Res}_{K'/K}(\Gm)$
and more generally to the case of an induced torus $G=\prod_i{\rm Res}_{K_i/K}(\Gm^{n_i})$.
\smallskip

{\sl Step 2.} The case that $G=T$ is a torus. Note that in this case we can identify
$\pi_1(T)$ with the cocharacter lattice $X_*(T)$. As in loc. cit. one uses an exact sequence
\begin{equation*}
R\xrightarrow{f} S\xrightarrow{g} T\to 1,
\end{equation*}
where $R$ and $S$ are induced tori and ${\rm ker}(f)$, ${\rm ker}(g)$ are tori. This induces an exact sequence
\begin{equation*}
\pi_1(R)_I\to \pi_1(S)_I\to \pi_1(T)_I \to 1\
\end{equation*}
and  this together with \ref{steinberg} and Step 1 suffices to construct $\kappa_T$.
\smallskip

{\sl Step 3.} The case where the derived group $G_{\rm der}$ of $G$ is simply connected.
In this case, $\pi_1(G)=\pi_1(D)=X_*(D)$ where $D=G/G_{\rm der}$ and we
define $\kappa_G$
as the composition
\begin{equation*}
LG(k)\to LD(k)\to \pi_1(D)_I\ .
\end{equation*}
\smallskip

{\sl Step 4.} The general case. One uses a $\sl z$-extension
\begin{equation*}
1\to S\to G'\to G\to 1
\end{equation*}
i.e a central extension with $G'_{\rm der}$ simply connected and $S$ an induced torus ([MS], Prop. 3.1), to reduce to the previous case, cf.~loc.~cit.
\smallskip

By  loc.~cit.~the Kottwitz homomorphism $\kappa_G$ is surjective.
 We will give a geometric interpretation of the Kottwitz homomorphism  and additional properties of it 
in the following sections.

The parahoric subgroup scheme $P_{\und a}$ leads by the procedure of the previous section to the fpqc sheaf
$\F_{P_{\und a}}=LG/L^+P_{\und a}$ which in the sequel we will often denote simply by $\F_{\und a}$;
we will call the corresponding ind-schemes (cf. Theorem \ref{represent})
{\it generalized affine flag varieties}.

\begin{Remarknumb}
{\rm   More generally, one can also consider the group schemes $P$ corresponding to the Moy-Prasad type
subgroups of $G(K)$ which are constructed by Yu in [Yu]. These are also smooth and affine over $\Spec(V)$
with generic fiber $G$. As an example,  consider the case $G={\rm GL}_n$. Take the group scheme $P_m$ given in [Yu] with $V$-valued points
\begin{equation}
P_m(V)=\{g\in {\rm GL}_n(k[[t]])\ |\ g\equiv{\rm Id}_n\mod t^m\}.
\end{equation}
Then  the fpqc-quotient $\F_{P_m}$ is a
${\rm Res}_{k[[t]]/(t^m)/k}{\rm GL}_n$-bundle over the affine Grassmannian
$\F_{P_0}$. The  morphism $\F_{P_m}\to \F_{P_0}$ represents the
forgetful functor
\begin{equation}
(\L, \ \sigma: \L/t^m\L \xrightarrow{\sim} (R[[t]]/(t^m))^n)\mapsto \L 
\end{equation}
where $\L$ denotes a lattice as in \ref{afgraexample}.}
\end{Remarknumb}
\bigskip
\bigskip

\section{Tori and groups of multiplicative type	}\label{tori}

 In this section we present various  examples.

\subsection{} a) For  $P=\Gm$, we have $L^+P=R[[t]]^\times$, $LP_\eta=R((t))^\times$.
For
\begin{equation*}
a=\sum_{i>>-\infty}^{+\infty} a_i t^i \in R((t))^\times
\end{equation*}
denote by ${\rm ord}(a)$ the smallest integer $n\in \Z$ for which $a_n$ is invertible in $R$. (We are assuming that $\Spec(R)$ 
is connected.) Then for all $m<n$, $a_m$ is nilpotent. When $R$ is reduced we can write $a=t^{\rm ord(a)}\cdot u$
with $u\in R[[t]]^\times$. The map ${\rm ord}$ gives a Zariski locally constant
map
\begin{equation*}
{\rm ord}: \F_{\Gm}\to \Z\ .
\end{equation*}
Consider the ind-scheme $L^{--}\Gm$
over $k$ which represents the functor which to a $k$-algebra $R$ associates
\begin{equation}
L^{--}\Gm(R)= (1+t^{-1}R[t^{-1}])^\times
\end{equation}
i.e.,  the elements of $L^{--}\Gm(R)$ are polynomials $1+a_1t^{-1}+\cdots +a_it^{-i}$ with $a_i$ nilpotent in $R$.
Then it is easily seen that $L^{--}\Gm$ can be identified with the open and closed sub-ind-scheme of  $\F_{\Gm}$ for which
${\rm ord}=0$. From this it follows that
that the fibers of ${\rm ord}$ are the connected components of $\F_{\Gm}$.
In fact, the topological space underlying the ind-scheme $\F_{\Gm}$
is $\Z$ and ${\rm ord}$ gives a homeomorphism. The connected
components of $\F_{\Gm}$ are infinitesimal ind-schemes; they are all
isomorphic to each other via morphisms given
by multiplication by the appropriate power of $t$. 
\smallskip

b) Consider $K'=k((u))$ with $u^2=t$ as a degree $2$ extension of $K$, and suppose ${\rm char}(k)\neq 2$. Let $P$ be the
group scheme over $V=k[[t]]$ which is given as the kernel of the norm
\begin{equation*}
P={\rm ker}({\rm Res}_{k[[u]]/k[[t]]}(\Gm)\xrightarrow{\rm N} \Gm)\  .
\end{equation*}
Note that if $R$ is reduced, then any $a\in R((u))^\times$ is of the form $u^m\cdot v$ with
$v\in R[[u]]^\times$; the equation $N(a)=1$ then gives $m=0$ and $N(v)=1$. This shows
$\F_P(R)=\{1\}$ for all reduced $R$ and therefore the corresponding
topological space is a point. On the other hand, if $R=\Spec(k[x]/(x^2))$ then
$a=xu^{-m}+1$, for any odd $m>0$, is a non-trivial element of $\F_P(R)$. We conclude that $\F_P$
is not reduced and is an infinitesimal ind-scheme.
\smallskip

c) Let $P=\mu_n$. Then $LP_\eta(R)=\mu_n(R((t)))$,
$L^+P(R)=\mu_n(R[[t]])$. If the characteristic of $k$ does not divide $n$, then $\mu_n(R((t)))=\mu_n(R[[t]])=\mu_n(R)$ and $\F_P=\Spec(k)$.
However, this is not true in general. In fact, if $n=p={\rm char}(k)$, then we can see that $\F_{\mu_p}$ is represented
by the ind-scheme
\begin{equation}
\prod_{i=-1}^{-\infty} \Spec(k[x_{i}]/(x^p_i))\ .
\end{equation}

\subsection{} Let $G=T$ be a torus over $K$.
In this case, the connected Neron model
of $T$ over $V$ (also called the Neron-Raynaud model in [Yu]) gives the
(unique) parahoric subgroup scheme of $T$ ([B-TII], 4.4). 
Let us briefly recall this construction below.
Let $K'$ denote a splitting field for $T$, so that $T_{K'}$ is a split torus.
We can assume that $K'/K$ is finite and separable. The
natural homomorphism
\begin{equation*}
T\to {\rm Res}_{K'/K}T_{K'}
\end{equation*}
realizes $T$ as a subgroup scheme of an induced torus. Denote by ${\cal T}_{V'}$ the
natural extension of the split torus $T_{K'}$ to a torus over $V'$ (if $T_{K'}\simeq ({\bf G}^h_m)_{K'}$, then
${\cal T}_{V'}\simeq ({\bf G}^h_m)_{V'}$.) Consider the Zariski closure ${\cal T}$
of $T$ in ${\rm Res}_{V'/V}({\cal T}_{V'})$; this is an affine flat group scheme
of finite type over $V'$ which is independent of all choices (see [B-TII], 4.4.6). One can construct its smoothening (``lissification") ${\cal T}^R$ following Raynaud (see [B-TII], 4.4.12). This is a (commutative) group scheme
affine and smooth over $V$ with generic fiber $T$; its special fiber is not necessarily connected.
(Let us remark here that if $K'/K$ is tamely ramified, then it is known by [Ed]
that ${\cal T}$ is  smooth over $V$ and so ${\cal T}={\cal T}^R$.)
We have ([B-TII], 4.4.12)
\begin{equation}\label{trv}
{\cal T}^R(V^\flat)=T_b(K^\flat):=\{t\in T(K^\flat)\ |\ {\rm val}(\chi(t))=0, \hbox{\rm for all }\chi\in X^*_{K^\flat}(T)\}.
\end{equation}
The connected Neron  model ${\cal T}^0$ of $T$ is  the connected component
of the smooth group scheme ${\cal T}^R$. This agrees with the connected component
of the locally finite type (lft) Neron model of the torus $T$ 
(see [BLR], \S 10), so there is no confusion in our terminology. Furthermore, ${\cal T}^R$ is the maximal subgroup scheme
of the lft Neron model of finite type over $V$, comp. [R], Prop.~on p.~314.

If $T$ is an induced torus, i.e.,  of the form
$T\simeq \prod_i{\rm Res}_{K_i/K}({\bf G}_m^{n_i})_{K_i}$, then by [B-TII],  4.4. we have
${\cal T}^0={\cal T}^R={\cal T}$.
Then by (\ref{weilre1}) and the above construction, there is a finite extension $k'$ of $k$  such that
\begin{equation}\label{induced}
LT\times_kk'\simeq (L{\bf G}_m)^n, \quad L^+{\cal T}\times_kk'\simeq (L^+\Gm)^n, \quad \F_{{\cal T}}\times_kk'\simeq (\F_{\Gm})^n\, ,
\end{equation}
with $n=[k':k]\cdot \sum_i n_i$.

Let us consider $\F_{{\cal T}^0}$ and $\F_{{\cal T}^R}$ for a general torus $T$.
By Theorem \ref{represent}, the $\bar k$-valued points of these ind-schemes
are dense and we have
\begin{equation}
\F_{{\cal T}^R}(\bar k)=T({\bar k}((t)))/{\cal T}^R({\bar k}[[t]]),\quad
\F_{{\cal T}^0}(\bar k)=T({\bar k}((t)))/{\cal T}^0({\bar k}[[t]])\ .
\end{equation}
By (\ref{trv}), $\F_{{\cal T}^R}(\bar k)=T(K^\flat)/T_b(K^\flat)\subset {\rm Hom}(X^*_{K^\flat}(T),\Z)$.
We have an exact sequence
\begin{equation}
0\to  {\cal T}^R({\bar k}[[t]])/{\cal T}^0({\bar k}[[t]])\to  \F_{{\cal T}^0}(\bar k)    \to \F_{{\cal T}^R}(\bar k)\to 0\ .
\end{equation}
It will follow from our discussion of the Kottwitz homomorphism
(\S \ref{kottwitzmap}, Cor. \ref{torusDiscrete}) that
this exact sequence can be identified with the exact sequence
\begin{equation}
0\to (X_*(T)_I)_{\rm tor}\to X_*(T)_I\to (X_*(T)_I)_{\rm cotor}\to 0\ .
\end{equation}
(Here  $I$ is the inertia subgroup of $\Gal(K^{\rm sep}/K)$.)

\subsubsection{}\label{normtorus}  
Consider $K'=k((u))$ with $u^2=t$ as
a degree $2$ extension of $K$ and suppose that
${\rm char}(k)\neq 2$.
Let $T$ be the torus over $K$ given as the kernel of the norm, 
$
T:={\rm ker}({\rm Res}_{K'/K}(\Gm)\xrightarrow{N} \Gm) .
$
In this case, the torus ${\cal T}$ above is the group scheme over $k[[t]]$
which is also given as the kernel of the norm
\begin{equation}\label{intnormtorus}
{\cal T}={\rm ker}({\rm Res}_{k[[u]]/k[[t]]}(\Gm)\xrightarrow{\rm N} \Gm)\  .
\end{equation}
Notice that ${\cal T}$ is the fixed point scheme of an involution
on ${\rm Res}_{k[[u]]/k[[t]]}(\Gm)$. Hence, it follows by [Ed], Prop 3.5, that
${\cal T}$ is smooth over $V$ and so ${\cal T}={\cal T}^R$. However, ${\cal T}$ is
not connected; the special fiber has two connected components corresponding to
units $a\in k[[u]]^\times$ with $a\cdot \bar a=1$ having residue class $1$, resp. $-1$ modulo $u$.
Hilbert's theorem 90 implies that there is an exact sequence
\begin{equation}\label{H90}
1\to {\Gm} \to {\rm Res}_{k[[u]]/k[[t]]}\Gm\xrightarrow{b\mapsto \bar b\cdot b^{-1}}{\cal T}^0\to 1
\end{equation}
over $k[[t]]$.
In this case, the Kottwitz   homomorphism
\begin{equation}
\kappa_T: T(\bar k((t)))\to \{\pm 1\}
\end{equation}
is given as follows: If $a\in \bar k((u))^\times$ with $a\cdot \bar a=1$, then write $a=\bar b\cdot b^{-1}$, $b\in \bar k((u))^\times$ and we have $\kappa_T(a)=(-1)^{{\rm ord}_u(b)}$.
Note that for $a\in  \bar k[[u]]^\times$ we have $a\equiv  \kappa_T(a)\mod (u)$
and so we have  ${\cal T}^0( \bar k[[t]])={\cal T}( \bar k[[t]])\cap {\rm ker}(\kappa_T)$, see also \S \ref{kottwitzmap} below.
\bigskip
\bigskip

\section{Affine flag varieties for quasi-split unitary groups}\label{unitary}
\setcounter{equation}{0}

We continue with an example which is essential for our applications. 

\subsection{} Let ${\rm char}(k)\neq 2$. Let $K'/K$ be a ramified quadratic extension
and fix a uniformizer $u$ of $K'$ with image under the Galois involution $\bar u=-u$, $u^2=t$. Let $W$ be a
$K'$-vector space of dimension $n\geq 2$ and let
\begin{equation*}
\phi: W\times W\to K'
\end{equation*}
be a $K'/K$-hermitian form. Let
\begin{equation*}
 U(W,\phi)=\{g\in GL_{K'}(W)\ |\ \phi(gv, gw)=\phi(v,w), \forall v, w\in W\}
\end{equation*}
and consider also $SU(W,\phi)=\{g\in U(W, \phi)\ |\  \det(g)=1\}$
(a semi-simple simply connected group over $K$).
In these cases, the group is not of the form $G_0\times_{\Spec(k)}\Spec(K)$.
Let us consider the case when the form $\phi$ is split. This means that there exists
a basis $e_1,\ldots ,e_n$ of $W$ such that
\begin{equation}\label{hermbasis}
\phi(e_i, e_j)=\delta_{i, n+1-j},\quad \forall\ i, j=1,\ldots, n\ .
\end{equation}
In what follows, we fix a choice of such a form and denote the groups by $U_n$, $SU_n$.
Since $SU_n$ is semi-simple simply connected,
 the Kottwitz homomorphism for $U_n$ is the composition
\begin{equation}
\kappa_{U_n}: U_n({\bar k}((t)))\xrightarrow{\rm det} T({\bar k}((t)))\xrightarrow {\kappa_T} \{\pm 1\}\ .
\end{equation}
Here $T$ is the torus considered in \ref{normtorus}.
Note that there are two associated $\O_K$-bilinear forms,
\begin{equation}
(x,y)={\rm Tr}_{K'/K}(\phi(x, y)),\quad <x,y>={\rm Tr}_{K'/K}(u^{-1}\phi(x,y))\ .
\end{equation}
The form $(, )$ is symmetric while $<,>$ is alternating. They satisfy the identities,
\begin{equation}
(x,uy)=-(ux, y), \quad <x,uy>=-<ux, y>\, .
\end{equation}
For any $\O_{K'}$-lattice $\Lambda$
in $W$, we set
\begin{equation*}
\hat\Lambda=\{w\in W\, |\, \phi(w,\Lambda)\subset \O_{K'}\}= \{v\in W\, |\, \ <v,\Lambda>\,\subset \O_{K}\}\, .
\end{equation*}
Similarly, we set
\begin{equation*}
\hat\Lambda^s=\{w\in W\, |\,  (w,\Lambda)\subset \O_{K}\},
\end{equation*}
so that $\hat\Lambda^s=u^{-1}\hat\Lambda.$
 For $i=0,\ldots, n-1$, let us set
\begin{equation*}
\lambda_i:={\rm span}_{\O_{K'}}\{u^{-1}e_1,\ldots, u^{-1}e_i, e_{i+1}, \ldots , e_n\}\ .
\end{equation*}
The lattice $\lambda_0$ is self-dual for the alternating form $<,>$.
We now distinguish two cases:

a) $n=2m+1\geq 3$ is odd. [Then $SU_n$ is quasi-split of type C-BC$_m$; see [T1], p. 60.]
We have  $\hat\lambda^s_{m}=\lambda_{m+1}$. Let $I$ be a non-empty
subset of $\{0,\ldots, m\}$ and consider the lattice chain consisting of $\lambda_i$, $i\in I$.
In fact, we can extend $\lambda_i$, $i\in I$, to a periodic self-dual lattice chain by first including the duals
$\hat\lambda^s_i=\lambda_{n-i}$ for $i\neq 0$, and then all the $u$-multiples of our lattices:
For $j\in \Z$ of the form $j=k\cdot n+i$ with $0\leq i<n$
we put
\begin{equation*}
\lambda_j=u^{-k}\cdot \lambda_i\ .
\end{equation*}
Then $\{\lambda_j\}_j$ forms a periodic lattice chain $\lambda_I$
($u\cdot \lambda_j=\lambda_{j-n}$) which satisfies
$\hat\lambda_j=\lambda_{-j}$. Now consider the subgroups
\begin{equation*}
P_I=\{g\in U_n\ |\ g\lambda_i=\lambda_i,\ \forall\ i\in I\}, \quad P'_I=\{g\in P_I\ |\ \det(g)=1\}
\end{equation*}
of $U_n(K)$ and $SU_n(K)$ respectively. We then have the following statements:

{\it The subgroup $P'_I$
is a parahoric subgroup of $SU_n(K)$. Any parahoric subgroup of $SU_n(K)$ is conjugate to a subgroup $P'_I$ for a unique subset $I$.  The sets $I=\{0\}$ and $I=\{m\}$ correspond to the special maximal parahoric subgroups.}

We indicate briefly how these assertions can be proved by following the method of Kaiser [K], \S 5.2 (in loc.~cit.~the symplectic group is considered). (Another possible reference is [T1], 1.15 and 3.11, or the work of Gan-Yu [G-Y].) We consider the embedding of the Bruhat-Tits building $\mathcal B$ of $SU_n(K)$ into
the Bruhat-Tits building $\mathcal B'$ of $SL_n(K')$. 
Denote by $z\mapsto \hat z$
the involution of $\mathcal B'$ induced by the map which to a lattice $\Lambda$ associates $\hat \Lambda$. The following points can be checked by imitating Kaiser's reasoning in the case at hand.

\begin{enumerate}
\item  Every maximal parahoric subgroup $P$ of $SU_n(K)$ fixes a vertex $z$ in $\mathcal B'$ such that $z$ and $\hat z$
are neighbouring vertices (or coincide), and, conversely, the unordered pair $(z, \hat z)$ determines $P$ uniquely. 

\item Conversely, the stabilizer of a vertex $z$ of $\mathcal B'$ of the kind described above is a maximal parahoric subgroup.

\item To a vertex $z$ as above we associate as follows a flat group scheme $\mathcal G$ over $k[[t]]$. 
Let $z$ be the homothety class of  a lattice $Z$ with $uZ\subset \hat Z \subset Z$. Let
$\mathcal G$ be the Zariski closure of $SU_n$ in $GL(Z)$. Then $\mathcal G$ is a smooth group scheme whose special 
fiber  has as its maximal reductive quotient the group $\bar {\mathcal G}=Sp(Z/\hat Z)\times SO(\hat Z/uZ)$. The map
$P\to \bar{\mathcal G}(k)$ is surjective. Furthermore, the link of the vertex in $\mathcal B$ corresponding to $P$ can 
be identified with the Tits building of $\bar{\mathcal G}$. 

\item It follows that the maximal parahoric subgroup $P$ is special if and only if either $Sp(Z/\hat Z)$ or $SO(\hat Z/uZ)$
is trivial, i.e., if either $Z=\hat Z$ or ${\rm dim}\  \hat Z/uZ=1$.

\item Fix a lattice $Z$ with $\hat Z=Z$ and a standard basis $e_1, \dots,e_n$ of it satisfying (\ref{hermbasis}). Then any
simplex in $\mathcal B$ is conjugate under $SU_n(K)$ to a simplex in the Tits building of the corresponding
group  $\bar{\mathcal G}$ which is unique up to conjugation; and hence is conjugate to a 
unique simplex in the fundamental chamber of the corresponding apartment of the Tits building of $SO(Z/uZ)$.
It follows that any parahoric subgroup of $SU_n(K)$ is conjugate to $P'_I$ for a unique subset $I$.

\end{enumerate}

Note that the subgroups $P_I$ are {\sl never} parahoric subgroups of $U_n(K)$.  One can see this by observing that $P_I$ always contains the unitary transformation $e_m\mapsto -e_m$, $e_i\mapsto e_i$ for $i\neq m$; the Kottwitz invariant of this element of $U_n(K)$ is non-trivial
and so the special fiber of $P_I$ is not connected. To obtain the corresponding parahoric subgroups of $U_n(K)$ we need to
intersect with the kernel of the Kottwitz homomorphism, i.e.,  consider
\begin{equation*}
P^0_I=\{g\in P_I\ |\ \kappa_{U_n}(g)=1\}\ .
\end{equation*}

b) $n=2m$. [We will always assume $m\geq 2$, in which case $SU_n$ is  quasi-split with local Dynkin diagram of type B-C$_m$ if $m\geq 3$, and  C-B$_2$ if $m=2$.
If $m=1$ then $SU_n\simeq SL_2$ and  is split.]
 In this case, we also introduce  
\begin{equation*}
\lambda_{m'}={\rm span}_{\O_{K'}}\{ u^{-1}e_1,\ldots, u^{-1}e_{m-1}, e_{m},  u^{-1}e_{m+1}, e_{m+2},\ldots , e_n\}\ .
\end{equation*}
Then both $\lambda_m$ and $\lambda_{m'}$ are self-dual for the symmetric form
$(, )$. Now consider non-empty subsets  $I\subset \{0,\ldots, m-2, m, m'\}$
(we view  $m'$   as an additional symbol)  and as above
consider the subgroups
\begin{equation*}
P_I=\{g\in U_n(K)\ |\ g\lambda_i=\lambda_i,\ \forall\ i\in I\},\quad P'_I=\{g\in P_I\ |\ \det(g)=1\}\ .
\end{equation*}
{\it The subgroup $P'_I$ is a parahoric subgroup of $SU_n(K)$.
Any parahoric subgroup of $SU_n(K)$ is conjugate to a subgroup $P'_I$
for a unique subset $I$. The sets $I=\{m\}$ and $I=\{m'\}$ correspond to the special maximal parahoric subgroups.}

The way to see these facts is similar to the case of odd $n$ considered above. But there are differences. In (b) above,
there is the exceptional case when $z$ is represented by a lattice $Z$ with $uZ\subset \hat Z\subset Z$, where ${\rm dim}\ \hat Z/uZ=2$. The action of the stabilizer of $z$  on the link of $z$ factors  through $\bar{\mathcal G}$ and fixes in this case  the two isotropic lines
in the 2-dimensional quadratic space $\hat Z/uZ$. Hence this stabilizer is not a maximal parahoric subgroup,
but is the intersection of the two maximal parahoric subgroups which fix the lattices corresponding to the isotropic lines. Therefore we have to exclude this kind of pair $(z, \hat z)$ in the correspondence
described in (a). Similarly, in (d), to identify the special maximal parahoric subgroups, the condition in this case becomes 
that $SO(\hat Z/uZ)$ is trivial, which now means that  $\hat Z=uZ$. But other than these minor differences the argument is the same.

Out of the lattices $\lambda_i$, $i\in I$, we can create a periodic self-dual lattice chain whose stabilizer in $SU_n(K)$ is equal to $P'_I$ 
in the following manner.  We first introduce a set $I^\sharp\subset \{0,\ldots , m-1, m\}$ that depends on $I$ and is defined as follows: If both $m$ and $m'$ belong to $I$ we set $I^\sharp$ to be the set $(I-\{m'\})\cup \{m-1\}$
which is obtained from $I$ by replacing $m'$ by $m-1$. We can see that we have $P'_I=P'_{I^\sharp}$, $P_I=P_{I^\sharp}$.
%(This is true only because we are considering the {\sl special} unitary group.) 
In all other cases, i.e., when not both $m$ and $m'$ are in $I$, set $I^\sharp=I$. Now  we can saturate $\lambda_i$, $i\in I^\sharp$, as above to create
a periodic self-dual lattice chain $\lambda_I$.

Note again that the subgroups $P_I$ are not parahoric subgroups of $U_n(K)$ unless
the set $I$ contains either $m$ or $m'$. Indeed, suppose that $I$  contains neither
$m$ nor $m'$. Then the unitary transformation $e_m\mapsto e_{m+1}$, $e_{m+1}\mapsto e_m$,
$e_i\mapsto e_i$ for $i\neq m, m+1$, belongs to $P_I$ and has
non-trivial Kottwitz invariant. Conversely, assume that $I$ contains for example $m$
so that $P_I\subset P_{\{m\}}$. Let $g\in P_{\{m\}}$ and consider the action of $g$ on $\lambda_m/u\lambda_m$. We can see that 
$g$ respects the perfect alternating form on $\lambda_m/u\lambda_m$ given by $s(\bar x, \bar y)=(x, u^{-1}\cdot y)$.
Therefore the determinant of $g$ is $1$ modulo $u$; by \ref{normtorus} the Kottwitz invariant of any such $g$ is trivial.
This is enough to imply that $P_I$ is parahoric. When  $I$  contains neither
$m$ nor $m'$,
to obtain the corresponding parahoric subgroup of $U_n(K)$ we need again to
consider the kernel of the Kottwitz homomorphism, i.e.,
\begin{equation*}
P^0_I=\{g\in P_I\ |\ \kappa_{U_n}(g)=1\}\ .
\end{equation*}

In both cases, we will also denote by $P_I$, $P'_I$ the affine group schemes
over $k[[t]]$ whose $S$-valued points for a $k[[t]]$-algebra $S$ are the
$k[[u]]\otimes_{k[[t]]}S$-automorphisms of the (polarized) chain $\lambda_{I}\otimes_{k[[t]]}S$,
that preserve the form $\phi\otimes_{k[[t]]}S$ (and in the case of $P'_I$ have 
$k[[u]]\otimes_{k[[t]]}S$-determinant $1$).  By [R-Z],  App. to ch. 3,  the group scheme $P_I$ is
smooth over $k[[t]]$ (This is a variant of Case III of loc.~cit., with $\O_F$
replaced by the d.v.r. $k[[u]]$ and $\O_{F_0}$ by $k[[t]]$; actually [R-Z] would
only allow $S$ for which $t$ is locally nilpotent, but this is enough to imply
our result). Therefore, since, as we will see below in Theorem \ref{components},  the Kottwitz homomorphism
is locally constant, the group scheme $P^0_I$ is also smooth and affine.
A similar argument shows that $P'_I$ is also smooth
over $k[[t]]$; therefore, by appealing to [B-TII], 1.7.,  $P'_I$ coincides with  the corresponding parahoric group scheme
 constructed in  [B-TII]. 
We have an  exact sequence
of group schemes over $k[[t]]$
\begin{equation}
1\to P'_I\to P_I \xrightarrow{\rm det} {\cal T},  \
\end{equation}
where ${\cal T}$ is the torus of (\ref{intnormtorus}). 
The neutral component $P_I^0$ is the preimage of the
neutral  component ${\cal T}^0$ of ${\cal T}$. This fact, together with
the description of the Kottwitz invariant in this case, implies that
the $\bar k[[t]]$-valued points of the neutral component $P^0_I$
of $P_I$ are given by the $\bar k[[t]]$-valued points of $P_I$ that
have trivial Kottwitz invariant (so our notation is consistent). This in turn
implies that $P^0_I$ coincides with  the corresponding parahoric group scheme
that is constructed in  [B-TII]. We have a short exact sequence
of (connected, smooth) group schemes over $k[[t]]$,
\begin{equation}
1\to P'_I\to P_I^0\xrightarrow{\rm det} {\cal T}^0\to 1\ .
\end{equation}

\subsection{}\label{4b}
Recall the definition above of $I^\sharp$ when $n$ is even. We set $I^\sharp=I$ if $n$ is odd. In both cases $I^\sharp$ is totally ordered, if in the case that $n=2m$ is even,
we consider $\{0,1,\ldots, m, m'\}=\{0,1,\ldots, m\}\cup \{m'\}$ with the partial order which extends the
standard order and for which $m'>j$ for $j=0,\ldots, m-1$ .
Let us now write $I^\sharp=\{i_0<i_1<\cdots <i_k\}$.
Then we can consider the ``standard" lattice chain
\begin{equation}\label{unichainSt}
\lambda_{i_0}\subset \lambda_{i_1}\subset \cdots \subset \lambda_{i_k}\subset u^{-1}\lambda_{i_0}\ .
\end{equation}
We can consider the functor $\F_{I^\sharp}$ which to a $k$-algebra $R$ associates
the set of $R[[u]]$-lattice chains
\begin{equation}\label{chains1}
L_{i_0}\subset L_{i_1}\subset \cdots \subset L_{i_k}\subset u^{-1}L_{i_0}
\end{equation}
in $W\hat\otimes_{K'}R=R((u))^n$ which satisfy the following conditions:

a) For any $q\in \{0,\ldots, k\}$, we have
\begin{equation}
L_{i_q}\subset u^{-1}\hat L_{i_q}\subset u^{-1}L_{i_q}\ ,
\end{equation}

b) The quotients $L_{i_{q+1}}/L_{{i_q}}$, $u^{-1}\hat L_{i_q}/L_{i_q}$, $u^{-1}L_{i_0}/L_{i_{q+1}}$ are projective $R$-modules
of rank equal to the rank of the corresponding quotients for the standard chain (\ref{unichainSt}) (when $q=k$, these conditions have to be interpreted in the obvious way).
\medskip

The ind-group scheme $LU_n$ over $k$ given by $LU_n(R)=U_n(R((t)))$  acts naturally on $\F_{I^\sharp}$.

\begin{thm}\label{latmod}
There is an $LU_n$-equivariant isomorphism
$ LU_n/L^+P_I\simeq \F_{I^\sharp}$ of sheaves for the fpqc topology.
\end{thm}

\begin{proof}
%Note  here that the quotient $LSU_n/L^+P_I$ was explained in the previous %section.

To prove the theorem it suffices to check the following two statements:

i) For any $R$,  the subgroup of $LU_n(R)$ consisting of elements that stabilize the lattice chain
$\lambda_{I^\sharp}\hat\otimes R$ agrees with $L^+P_I(R)=P_I(R[[t]])$.

ii) Let $L_{I^\sharp}\in \F_I(R)$ be a lattice chain as above. Then
locally for the \'etale topology on $R$, there exists an isomorphism
$\lambda_{I^\sharp}\hat\otimes R\simeq L_{I^\sharp}$
realized by multiplication by an element $g\in LU_n(R)$.

Both of these statements essentially follow from the work in [R-Z], app. \!to ch. \!3. 
Indeed, to see (i) observe that we have $\lambda_{I^\sharp}\hat\otimes_kR=\lambda_{I^\sharp}\otimes_{k[[t]]}{R[[t]]}$
and we can consider $\lambda_{I^\sharp}\otimes_{k[[t]]}{R[[t]]}$ as a polarized periodic
chain of $k[[u]]\otimes_{k[[t]]}R[[t]]$-modules (terminology of loc. cit.)
One can see that, even in the case when $n=2m$ is even and when $\{m, m' \}\subset I$,
 the  group of automorphisms
of this chain that preserve the form $\phi\otimes_{k[[t]]}R[[t]]$ is equal to  $L^+P_I(R)=P_I(R[[t]])$. Using
$\lambda_i\otimes_{k[[t]]}R[[t]]\subset R((u))^n$ we see that these automorphisms
are exactly the  elements of $U_n(R((u)))=LU_n(R)$ that stabilize the lattice chain $\lambda_{I^\sharp}\hat\otimes_kR$.

Let us now discuss part (ii). We can see {\sl a priori} that $\F_{I^\sharp}$ is represented
by an ind-scheme of ind-finite type over $k$. By part (i), we obtain $LU_n/L^+P_I\to \F_{I^\sharp}$;
it is enough to show (ii) for $R$ strictly henselian. Now observe that [R-Z], app. \!to ch. \!3, modified for this situation, shows that there is  an isomorphism
\begin{equation}
\sigma_i: \lambda_{i}\hat\otimes_k R \xrightarrow{\sim} \L_i
\end{equation}
of polarized chains of $R[[u]]$-modules ( here we need to observe that the arguments of [R-Z]
extend to this case in which $t$ is topologically nilpotent.) Since both
$\lambda_i\hat\otimes_kR$ and $\L_i $ are lattices in $R((u)^n$,
the isomorphism $\sigma=\sigma_i[t^{-1}]$ is independent of $i$ and is given by an element $\sigma$ of
$U_n(R((t)))=LU(R)$.
\end{proof}

\begin{Remarknumb}\label{uniRem}
{\rm A. By our discussion above, when $P_I$ is not a parahoric subgroup of $U_n(K)$, then its special fiber has two connected components and $L^+P_I/L^+P^0_I\simeq \{\pm 1\}$, where the isomorphism is given by the Kottwitz invariant. We will see in the following paragraphs that $LU_n$ has
two connected components also distinguished by the value $\pm 1$ of the Kottwitz invariant. This will
allow us to conclude that if $P_I$ is not parahoric,
then for the parahoric subgroup $P^0_I$ we have $LU_n/L^+P^0_I\simeq \F_{I^\sharp}\sqcup \F_{I^\sharp}$.

 B. Let us discuss here the variant of the special unitary group $SU_n=SU(W,\phi)$.
In this case, let $\F'_{I^\sharp}$ be the functor which to a $k$-algebra $R$ associates
the set of $R[[u]]$-lattice chains as in (\ref{chains1}) which satisfy, in addition
to (a) and (b) above, also

(c) For all $q=0,\ldots, k$, we have
\begin{equation}\label{chains2}
\det(\L_{i_q})=\det(\lambda_{i_q}) \subset R((u))=\wedge^nR((u))^n.
\end{equation}

The ind-group $LSU_n$ acts on $\F'_{I^\sharp}$ and we obtain $\pi_I: LSU_n/L^+P'_I\to \F'_{I^\sharp}$.
We will see that $\pi_I$ is an isomorphism when $n$ is odd, or when $n=2m$
is even and $I$  contains neither $m$ nor $m'$. Indeed, we can always  identify as before
 $L^+P'_I$ with the stabilizer in $LSU_n$ of the standard lattice chain $\lambda_{I^\sharp}$;
the issue is again if each point $\L_{I^\sharp}$ in $\F'_{I^\sharp}(R)$ for $R$ strictly henselian
is in the orbit of the standard lattice chain $\lambda_{I^\sharp}$ by an element
$g\in LSU_n(R)=SU_n(R((t))$. By the proof of the theorem above, we can find
 $\sigma\in U_n(R((t)))$ such that $\sigma\cdot (\lambda_i\hat\otimes_kR)=\L_i$, for all $i$.
 Our additional condition (c) now guarantees that the determinant $a=\det(\sigma)$ is a unit in
 $R[[u]]^\times$. Since $\sigma$ is unitary, we also have ${\rm Norm}(a)=1$ and the exact sequence (\ref{H90})
 now implies  that then there is $b\in R[[u]]^\times$ such that $a=\pm \bar b\cdot b^{-1}$. Now suppose first that $n=2m+1$ is odd. Then we can compose $\sigma$ with the unitary automorphism $e_m\mapsto a^{-1}e_m$,  $e_i\mapsto e_i$, for $i\neq m$, which stabilizes $\lambda_{I^\sharp}$. The resulting unitary automorphism has determinant $1$ and takes $\lambda_i\hat\otimes_kR$ to $\L_i$. Suppose
 now that $n=2m$ and that $I$  contains neither $m$ nor $m'$. Assume that $a=-\bar b\cdot b^{-1}$.
 Then the unitary transformation $e_m\mapsto b e_{m+1}$, $e_{m+1}\mapsto \bar b^{-1}e_m$,
 $e_i\mapsto e_i$ if $i\neq m$, $m+1$, has determinant $a^{-1}$ and stabilizes $\lambda_{I^\sharp}$.
 Composing it with $\sigma$ gives the desired result. If $a=\bar b\cdot b^{-1}$ we can use
$e_m\mapsto b e_{m }$, $e_{m+1}\mapsto \bar b^{-1}e_{m+1}$,
 $e_i\mapsto e_i$ if $i\neq m$, $m+1$ instead.
 
 C. In the case when $n=2m$ is even and $\{m, m' \}\subset I$, our definition of $I^{\sharp}$ is asymmetric (we favor $\lambda_m$ over $\lambda_{m'}$). Note, however,
 that the element $g\in U_n(K)$ with $e_m\mapsto e_{m+1}, e_{m+1} \mapsto e_m$, and with $ e_i \mapsto e_i$ for $i\neq m, m+1$ takes $\lambda_m$ to $\lambda_{m'}$. 
 The element $g$ can be used to define an isomorphism between $\F_{I^\sharp}$ and 
 $\F_{I^{\sharp '}}$, where $\F_{I^{\sharp '}}$ denotes the version where we favor $\lambda_{m'}$ over $\lambda_m$.
 A slightly different point of view is explained in [P-R3]; there we consider the 
 groups of unitary similitudes.
 
  }
\end{Remarknumb}
\bigskip
\bigskip

\section{The Kottwitz map and connected components}\label{kottwitzmap}

\setcounter{equation}{0}

\subsection{} We continue with the following assumptions:
 $G$ is a connected reductive group over $K=k((t))$ with
 $k$ algebraically closed.
 Let $P$ be a parahoric group scheme of $G$
 in the sense of \S \ref{parahoric} and consider $\F_P:=LG/L^+P$.
By Theorem \ref{represent} the morphism $LG\to \F_P=LG/L^+P$ splits locally for the
 \'etale topology. Therefore, we have
 \begin{equation}
 \F_P(k)=LG(k)/L^+P(k)=G(k((t)))/P(k[[t]]).
 \end{equation}
Since by [H-R], Prop. 3, the Kottwitz homomorphism
\begin{equation*}
 \kappa=\kappa_G: G(k((t)))\to \pi_1(G)_I
 \end{equation*}
  is trivial on all parahoric subgroups
of $G(k((t)))$ we also obtain a map
 \begin{equation*}
\bar \kappa: \F_P(k)\to \pi_1(G)_I\ .
 \end{equation*}
 Note that since $\F_P$ is of ind-finite type, the $k$-points of $\F_P$ are dense for the Zariski topology.
 In what follows, we will denote by $\pi_0$ the set of connected components.
 
\begin{thm}\label{components}
The maps $\kappa : LG(k)\to \pi_1(G)_I$, $\bar \kappa: \F_P(k)\to \pi_1(G)_I$ are locally constant for the Zariski topology of the ind-schemes $LG$ and $\F_P$. In fact, they
induce bijections
 \begin{equation*}
 \pi_0(LG)\xrightarrow{\sim} \pi_1(G)_I, \quad \pi_0(\F_P)\xrightarrow{\sim} \pi_1(G)_I\ .
 \end{equation*}
\end{thm}

\begin{Proof}
Note that since $P$ is a smooth affine group scheme over $k[[t]]$ with connected fibers,
it follows from [GrII], p. 263, that $L^+P$ is also a connected group scheme
over $k$. Therefore, since by Theorem \ref{represent} the morphism $LG\to \F_P=LG/L^+P$
splits over points $\Spec(R)\to \F_P$ with values in strictly henselian rings $R$,
we see that it induces a bijection $\pi_0(LG)\simeq \pi_0(\F_P)$ on the set of connected components.

We will divide the proof into three cases, that essentially correspond to the steps in the definition
of the Kottwitz homomorphism:
\smallskip

A. Let us first discuss the case that $G=T$ is a torus. Then the connected component
of the Neron model $P={\cal T}^0$ is the unique parahoric subgroup
scheme of $T$. By  [R], appendix, p.~380,   ${\cal T}^0(k[[t]])$ 
coincides with the kernel of the Kottwitz homomorphism
\begin{equation*}
\kappa_T: T(K)\to \pi_1(T)_I=X_*(T)_I \,
\end{equation*}
which is surjective.
Therefore, our claims for $G=T$ will follow if we show that
$L^+P$ is the neutral component of the ind-scheme $LT$.
Alternatively,
it is enough to show that the connected component of $\F_{P}$ that contains
the identity coset $e\, {\rm mod}\, L^+P$ is equal to $\{e\, {\rm mod}\, L^+P\}$.
Note that $\F_P$ is of ind-finite type. Therefore, to show this
statement it is enough to show the following:
\smallskip

{\bf Claim:} {\it Let 
$\phi: Y=\Spec(R)\to \F_P$ be a morphism with $Y$  a connected affine smooth curve over $k$
and let $y$ be a point in $Y(k)$  such that $\phi(y)=e\, {\rm mod}\, L^+P$.
Then the morphism $\phi$ factors through the identity coset $\Spec(R)\to \Spec(k)\to \F_P$.}
\smallskip

To prove the claim observe that there is a (connected) \'etale cover $Z\to Y$ and a morphism $Z\to LT$ which lifts $Y\to \F_P$. By translating  we can assume that there is a $k$-valued point $z\in Z(k)$ which maps to the unit element in $LT$. It will be enough to show that when $Z=\Spec(S)$ is a connected affine smooth curve and $z\in Z$ is a closed point, any morphism $a: Z\to LT$ with $a(z)=e$ factors through $Z\to L^+P$. Now  $a$ corresponds to an $S((t))$-valued point $a$ of $T$. Let $\hat S_z\simeq k[[x]]$ be the
completion of the local ring of $S$ at $z$.   Consider the
scheme $W=\Spec(\hat S_z[[t]])$ with generic fiber $W_\eta=\Spec(\hat S_z((t)))$.
Consider the composed morphism
\begin{equation*}
f: W_\eta=\Spec(\hat S_z((t)))\to \Spec(S((t)))\to T
\end{equation*}
obtained from $Z\to LT$. By [BLR] 7.2/3  extended to locally finite type
Neron models (this is done by replacing in the proof of 7.2/3
the use of loc.~cit.~7.2/1 (ii) by 10.1/3)  we can now see that the morphism $f$ extends
to
\begin{equation*}
\tilde f: W=\Spec(\hat S_z[[t]])\to {\cal T}^{lft N} \ ,
\end{equation*}
where ${\cal T}^{lft N}$ is the locally finite type Neron model of $T$
(this is a strengthening of the Neron
extension property to this situation).
By our assumption, the specialization $\Spec(k[[t]])\to {\cal T}^{lft N}$ of this morphism factors
through the neutral component, $\Spec(k[[t]])\to P={\cal T}^0$; it follows that $f$ actually extends to
$\tilde f: W\to P={\cal T}^0$. Consider the corresponding ring homomorphism
$\tilde f^*: \O(P)\to \hat S_z[[t]]$; by our construction,
the composition of this homomorphism with $\hat S_z[[t]]\subset \hat S_z((t))$
is equal to
\begin{equation*}
\O(P)\subset \O(P)\otimes_{k[[t]]}k((t))=\O(T)\xrightarrow {a^*}S((t))\ .
\end{equation*}
We  conclude that this last homomorphism takes values in $\hat S_z[[t]]\cap S((t))=S[[t]]$;
this shows that $Z\to LT$ factors through $Z\to L^+P$. This completes the proof in the
case that $G=T$ is a torus. In particular, we see that in this case the Kottwitz homomorphism
induces a bijection,
\begin{equation}
\bar\kappa_T: \F_{{\cal T}^0}(k)\xrightarrow{\sim} X_*(T)_I\ .
\end{equation}

Before we continue with the other cases let us remark that our basic task is to
show that the Kottwitz homomorphism $\kappa_G: LG(k)\to \pi_1(G)_I$, or equivalently
$\bar\kappa_G: \F_P(k)\to \pi_1(G)_I$ is locally constant. (Our argument in case A
establishes this when $G$ is a torus.)
Assuming this we see that $\kappa_G$ induces a group homomorphism
\begin{equation*}
 \pi_0(LG)\simeq \pi_0(\F_P)\to \pi_1(G)_I\ .
\end{equation*}
This is surjective, since $\kappa_G$ is surjective. To show that $\pi_0(LG)\to \pi_1(G)_I$ is injective
it is enough to show that the kernel  $LG(k)_1\subset LG(k)$ of
$\kappa_G$ lies in the neutral component of $LG$.
By [H-R],   Lemma 17, $LG(k)_1$ is generated by  the union of all parahoric subgroups
$L^+P(k)=P(k[[t]])$ of $G$. As mentioned in the beginning of our proof,  $L^+P$ is connected for each parahoric $P$, hence the
 images of all multiplication morphisms $L^+P_{i_1}\times \cdots \times L^+P_{i_n}\to LG$ (with source an arbitrary finite product
 of parahorics) are contained in
the neutral component of $LG$. We conclude that $LG(k)_1$ is contained in the neutral component
of $LG$ and so $\pi_0(LG)\to \pi_1(G)_I$ is injective.

In what follows, we will only deal with the  claim that $\kappa_G$ is locally constant
since by the above this is enough to complete the proof.
\smallskip

B. The case that the derived group $G_{\rm der}$ is simply-connected.
Then the Kottwitz homomorphism is the composition
\begin{equation*}
LG(k)\to LD(k)\xrightarrow{\kappa_D} \pi_1(D)_I=\pi_1(G)_I
\end{equation*}
where $D$ is the torus $G/G_{\rm der}$. Our claim that $\kappa_G$ is locally constant then
follows from case A applied to $T=D$, since $LG(k)\to LD(k)$ is induced by a morphism of ind-schemes.
\smallskip

C. The general case. We can find a ${\sl z}$-extension
(see  \ref{parahoric}, Step 4)
\begin{equation*}
1\to S\to G'\to G\to 1\ .
\end{equation*}
 As in
case A it is enough to show that if $Z$ is a connected affine smooth curve over $k$
and a morphism $a: Z\to LG$ is given, the Kottwitz homomorphism is constant on the image
$a(Z(k))\subset LG(k)$. Let $z, z'\in Z(k)$ and consider $W=\Spec(\O_{Z,{z;z'}}((t)))$
where $\O_{Z,{z;z'}}$ is the semi-local ring of $Z$ at $z$, $z'$.
Since $\kappa_G$ is surjective, we can assume without 
loss of generality that $\kappa_G(a(z))=0$. 
Since $\O_{Z,z;z'}$ is a UFD, the Picard group of $W$ is trivial. Observe
that for every finite separable extension $L/K$, we have $L\simeq k((t))$
and so $W\times_KL\simeq W$. 

A standard argument using the
Hochschild-Serre spectral sequence now gives ${\rm H}^1_{et}(W, S_W)=(0)$
for the induced torus $S$:
Indeed, it is enough to explain this for $S={\rm Res}_{K'/K}\Gm$ with $K'/K$
a finite separable extension. Let $\tilde K=k((\tilde u))$ be the Galois hull of $K'$ over $K$;
set $\Gamma=\Gal(\tilde K/K)$ and denote by $H$ the subgroup that corresponds to 
$K'$. The Hochschild-Serre spectral sequence for $W\times_K\tilde K\to W$
now is
\begin{equation*}
{\rm H}^p(\Gamma, {\rm H}^q_{et}(W\times_K\tilde K, S_{W\times_K\tilde K}))\Rightarrow {\rm H}^{p+q}_{et}(W, S_W)\ .
\end{equation*}
Now $S_{W\times_K\tilde K}$ is a split torus. Hence, by the above
\begin{equation*}
{\rm H}^1(W\times_K\tilde K, S_{W\times_K\tilde K})=(0),
\end{equation*} 
while also
\begin{equation*}
{\rm H}^0(W\times_K\tilde K, S_{W\times_K\tilde K})={\rm Ind}_H^\Gamma(\Z)\otimes (\O_{Z,{z;z'}}((\tilde u)))^\times\ .
\end{equation*}
Shapiro's lemma now shows that the cohomology group ${\rm H}^1(\Gamma, -)$ of this last $\Gamma$-module is
isomorphic to ${\rm H}^1(H, (\O_{Z,{z;z'}}((\tilde u)))^\times)$. A similar argument as above applied
to the $H$-cover $W\times_K\tilde K\to W\times_K K'$ now shows that since ${\rm H}^1(W\times_KK', S_{W\times_KK'})=(0)$,
we have ${\rm H}^1(H, (\O_{Z,{z;z'}}((\tilde u)))^\times)=(0)$. 
We  conclude that indeed ${\rm H}^1_{et}(W, S_W)=(0)$, as desired.

Since $S\to \Spec(k[[t]])$ is smooth, ${\rm H}^1_{et}(W, S_W)=(0)$ 
implies that the $S$-torsor
given by $G'\to G$ splits after base-changing along the morphism
$W\to G$ given by $a$. Therefore, 
the morphism $W\to G$ lifts to $W\to G'$; this gives $\Spec(\O_{Z,z;z'})\to LG'$
that factors $a_{z;z'}: \Spec(\O_{Z,z;z'})\to LG$. Now the image of $\Spec(\O_{Z,z;z'})$ lies in 
a single connected component of $LG'$. Since we know the conclusion of the Theorem for $G'$ 
(case B) we deduce that $\kappa_{G'}(z)=\kappa_{G'}(z')$. By the definition of
$\kappa_G$ this implies $\kappa_G(z)=\kappa_G(z')$.
\endproof
\end{Proof}
\smallskip

We obtain:

\begin{cor}\label{torusDiscrete}
Assume that $T$ is a torus over $K=k((t))$, where $k$ is algebraically closed.
Let ${\mathcal T}^0$ be the connected component of the N\'eron model of $T$ over $k[[t]]$.
Then the reduced locus $(LT/L^+{\mathcal T}^0)_{\rm red}=(\F_{{\mathcal T}^0})_{\rm red}$ is a discrete set which can
be identified with $X_*(T)_I$.
\end{cor}

\begin{cor}
Assume that $G$ is a reductive group over $K=k((t))$ with $k$ algebraically closed.

a) If $G$ is split over $K$, then $\pi_0(LG)=\pi_1(G)$.

b) If $G$ is semi-simple and simply connected then $LG$ is connected.
\end{cor}

We point out the analogy of the statement in (a) to the corresponding statement in topology. This
statement (for $k=\C$) appears in [BLS], 1.2, and in [B-D], 4.5.4. The heuristic that semi-simple simply connected groups over
a non-archimedean field behave like connected groups in the classical theory appears quite often in the
papers of Bruhat and Tits, comp., e.g., [B-TII], 5.2.11.

\bigskip
\bigskip

\section{Reducedness;  reductions and examples} \label{reduced}

\setcounter{equation}{0}

We continue with the following assumptions: $k$ is a perfect field of characteristic $p$ and
$G$ a connected reductive group over $K=k((t))$.
Consider the generalized affine flag variety $\F_{\und a}$ associated to
$G$ and to the facet $\und a$ of the Bruhat-Tits building of $G(K)$. The following
 theorem is one of the main results of this paper.

\begin{thm}\label{redu}
Suppose that $G$ is semi-simple, splits over a tamely ramified extension of $K$
and that the order of the fundamental group
$\pi_1(G)$ is prime to the characteristic of $k$.  
Then the ind-schemes $LG$ and $\F_{\und a}$, for any facet
$\und a$, are reduced.
\end{thm}

In this section, we will show that it will be enough to prove Theorem \ref{redu} under the additional
assumptions that $k$ is algebraically closed and $G$ is simply connected. The proof of this result
and therefore also of Theorem \ref{redu} will then
be completed in \S \ref{proofmain}. First of all, let us observe that the reduction to the case that
$k$ is algebraically closed
is immediate since $k$ is perfect. Observe also that, since $P_{\und a}$ is a smooth affine
scheme over $k[[t]]$, it follows  from  [GrII], p. 264, that
$L^+P_{\und a}$ is a reduced scheme. Hence, since $LG\to \F_{\und a}=LG/L^+P_{\und a}$
splits locally for the \'etale topology, $\F_{\und a}$ is reduced if and
only if  $LG$ is reduced.

\begin{Remarknumb}
{\rm a)  When $k$ has characteristic $0$ and $G$ is split, then the result follows
from Proposition 4.6 in [L-S]. This proof uses a theorem of Shafarevich
which is particular to characteristic $0$. For $G=SL_n$, the result is in [B-L]. 

b) We do not know if the assumption that $G$ splits over a 
tamely ramified extension of $K$ is necessary. As we shall see
in Remark \ref{pgl2}, $LPGL_2$ is non-reduced in characteristic $2$,
so the assumption $p\nmid\#\pi_1(G)$ appears necessary. We shall also see
in Proposition \ref{nonsemi}
that $LG$ is non-reduced for any reductive group $G$ that is not semi-simple. 
}
\end{Remarknumb}

We will now show how to reduce the proof of Theorem \ref{redu} to the case that $G$
is simply connected.

\subsection{} \label{Parcentral}
 Let $k$ be algebraically closed
and assume that $G$ is semi-simple over $K$. Then by a result of Steinberg
 $G$ is quasi-split ([St2]). Denote by $\tilde G$
its simply connected cover,
\begin{equation}\label{zg}
1\to Z\to \tilde G\to G\to 1\ .
\end{equation}
Here the center $Z$ is a finite commutative group scheme over $K$.
Let $T$ be a maximal torus of $G$ and consider  its inverse image $\tilde T$ in
$\tilde G$; this is a maximal torus of $\tilde G$.
The group scheme $Z$ is a subgroup  of  $\tilde T$ 
and hence it is a finite multiplicative group scheme over $K$
(i.e.,  the
Cartier dual $Z^D$ of $Z$ is a finite \'etale
commutative group scheme over $K$). Actually, since $\tilde G$
is semi-simple and simply connected, $\tilde T=\prod_{i}{\rm Res}_{K_i/K}(\Gm)$ is an induced torus, e.g. [B-TII], Prop. 4.4.16.
Hence, if $Z$ is annihilated by $n$, we have 
\begin{equation}\label{zindu}
Z\subset \prod_{i}{\rm Res}_{K_i/K}(\mu_n).
\end{equation}
Observe that
\begin{equation}
\pi_1(G)=Z(K^{\rm sep})(-1)
\end{equation}
as $I$-modules, where $(-1)$ denotes a (negative) Tate twist (and where $I={\rm Gal}(K^{\rm sep}/K)$).

 Now let $R$ be a $k$-algebra and consider the exact sequence
derived from (\ref{zg}) above,
\begin{equation}\label{longzg}
1\to Z(R((t)))\to \tilde G(R((t)))\to G(R((t)))\to {\rm H}^1_{\rm fppf}(R((t)), Z)\to {\rm H}^1_{\rm fppf}(R((t)), \tilde G) .
\end{equation}
By a result of Grothendieck [Gr], Thm. 11.7,   ${\rm H}^1_{\rm fppf}(R((t)), \tilde G)={\rm H}^1_{\rm et}(R((t)), \tilde G)$.
When $R$ is Artinian with algebraically closed residue field this latter group is trivial by
\S \ref{steinberg}. In this case, we have an exact sequence
\begin{equation}\label{longzg2}
1\to Z(R((t)))\to \tilde G(R((t)))\to G(R((t)))\to {\rm H}^1_{\rm fppf}(R((t)), Z)\to 1  .
\end{equation}

\subsubsection{}\label{Parcentral2} In this subsection, we suppose in addition that 
the characteristic of $k$ does not divide the order of $\pi_1(G)$.
Then $Z$ is finite \'etale over $k((t))$; for simplicity, we will identify $Z$ with 
the $I$-module given by $Z(K^{\rm sep})$. Then we can write
 $Z=\pi_1(G)(1)$.
Let us   assume that $R$ is  a local $k$-algebra with residue field $k'$, where $k'$ is some algebraically closed field extension of $k$. Then,
we see from (\ref{zindu}) that we have
\begin{equation}\label{zconstant}
Z(R((t)))=Z(k((t)))=Z(k[[t]])=(\pi_1(G)(1))^I.
\end{equation}
These equalities imply that $Z(R((t)))\subset \tilde {\cal T}^{0}(R[[t]])=(\prod_i{\rm Res}_{\O_{K_i}/\O_K}\Gm)(R[[t]])$ for each maximal torus $\tilde T$ of $\tilde G$.
Therefore, since each parahoric subgroup scheme $\tilde P_{\und a}$ of $\tilde G$
contains the connected component of the Neron  
model $\tilde {\cal T}^{0}$ for a maximal torus $\tilde T$ we conclude that 
$Z(R((t))\subset \tilde P_{\und a}(R[[t]])$, for each $\und a$.
By [B-TII], 1.7, the homomorphism $\ti G\to G$ induces $\ti P_{\und a}\to P_{\und a}$.
Hence, the exact sequence (\ref{longzg}) shows the existence 
of an immersion of ind-schemes
\begin{equation}\label{immersion}
L\tilde G/L^+\tilde P_{\und a}\hookrightarrow LG/L^+P_{\und a}\ .
\end{equation}

\begin{Remarknumb}\label{central}
{\rm Observe that (\ref{zconstant}) is not true when $p={\rm char}(k)$ divides the order of $Z$.
For example, if $R$ contains $a\neq 0$ with $a^p=0$, then
 $at^{-1}+1$ is in $\mu_p(R((t)))$ but not in $\mu_p(R[[t]])$.
 As a result, when $p$ divides the order of $Z$, $LZ(R)=Z(R((t)))$ is not necessarily
 contained in the  parahoric points $L^+\tilde P_{\und a}(R)=\tilde P_{\und a}(R[[t]])$.
 For example, if $p=2$ and $R$ contains $a\neq 0$ with $a^2=0$, the central element
 \begin{equation*}
 z=\left(\begin{matrix}
 1+at^{-1}& 0\cr
 0& 1+at^{-1}
 \end{matrix}\right)\in SL_2(R((t)))
 \end{equation*}
 does not belong to $SL_2(R[[t]])$.
}
\end{Remarknumb}

Now suppose also that
$R$ is Artinian with algebraically closed residue field $k'$. Then $R((t))$ is local henselian
with residue field $k'((t))$,
and we have
\begin{equation}
{\rm H}^1_{\rm fppf}(R((t)), Z)={\rm H}^1_{\rm et}(R((t)), Z)={\rm H}^1_{\rm et}(k'((t)), Z).
\end{equation}
Recall that if $B$ is a finite abelian group, and $k'$ an algebraically closed
field of characteristic prime to the order of $B$, then we have an isomorphism
\begin{equation}\label{kummer1}
{\rm H}^1_{\rm et}(k'((t)), B(1))=k'((t))^\times\otimes_\Z B\xrightarrow{\sim} B
\end{equation}
where the last isomorphism is given using the valuation ${\rm ord}: k'((t))^\times\to \Z$.
Similarly, for the finite \'etale commutative group  scheme $A(1)$ over $k((t))$ (which
we may think of as an $I$-module with ${\rm Gal}(K^{\rm sep}/K)=I$), we have   
\begin{equation}\label{kummer2}
{\rm H}^1_{\rm et}(k'((t)), A(1))\xrightarrow{\sim} A_I.
\end{equation}
Hence, under our assumptions, 
the exact sequence (\ref{longzg2}) becomes
\begin{equation}\label{longzg3}
1\to (\pi_1(G)(1))^I\to \tilde G(R((t)))\to G(R((t)))\xrightarrow{\kappa}  \pi_1(G)_I\to 1  .
\end{equation}
(One can show that the homomorphism $\kappa$ is given by the Kottwitz homomorphism $\kappa_G$.)
Recall that by Theorem \ref{components}, $L\ti G/L^+\ti P_{\und a}$ is connected.
The exact sequence (\ref{longzg3}) %together with Theorem \ref{components}  
now implies that the immersion
(\ref{immersion}) identifies $L\ti G/L^+\ti P_{\und a}$ with the neutral component 
of $L  G/L^+  P_{\und a}$. Furthermore, $L  G/L^+  P_{\und a}$ is reduced if and only if
$L\ti G/L^+\ti P_{\und a}$ is reduced. This concludes the reduction of the proof of Theorem \ref{redu} to the simply connected case. 

\begin{Remarknumb}\label{pgl2}
{\rm To see that the condition $p\nmid\# \pi_1(G)$ in the statement of the Theorem
\ref{redu} is necessary, suppose that $k$ is an algebraically closed field of characteristic $2$ and let
$G=PGL_2$ so that $\ti G=SL_2$ and $\pi_1(G)=\mu_2$.  Consider the fpqc sheaf $Q$ associated to the presheaf
$R\mapsto (R((t)))^\times/(R((t))^\times)^2R[[t]]^\times$, for a $k$-algebra $R$. We can see that $Q$ is represented by an 
ind-group scheme which is of ind-finite type and has two connected components: $Q=\Z/2\Z\times Q^0$.
The ind-scheme $Q^0$ is also given by the fpqc quotient $(1+t^{-1}R[t^{-1}])^\times/((1+t^{-1}R[t^{-1}])^\times)^2$;
it is not hard to see that it is  infinitesimal. Now suppose that $R$ is an Artin local
$k$-algebra with algebraically closed residue field. By comparing the exact sequences $(\ref{longzg2})$ for $G=\Gm$ and
$G=PGL_2$, we see that we have an exact sequence
\begin{equation}
1\to \mu_2(R((t)))\to SL_2(R((t)))\to PGL_2(R((t)))\to  (R((t)))^\times/(R((t))^\times)^2\to 1\ .
\end{equation}
Similarly, we have
\begin{equation}
1\to \mu_2(R[[t]])\to SL_2(R[[t]])\to PGL_2(R[[t]])\to  (R[[t]])^\times/(R[[t]]^\times)^2\to 1\ .
\end{equation}
This gives an exact sequence of pointed sets
\begin{equation}\label{temp2}
1\to \mu_2(R((t)))/\mu_2(R[[t]])\to \ti \F(R)\xrightarrow{\phi(R)} \F(R)\to (R((t)))^\times/(R((t))^\times)^2R[[t]]^\times\to 1
\end{equation}
where $\ti \F$, $\F$ are the affine Grassmanians for $SL_2$, $PGL_2$ respectively and $\phi$ is the
natural morphism (which is proper). When $R=k$, we have
\begin{equation}
1\to \ti \F(k)\to \F(k)\to \Z/2\Z\to 1\ .
\end{equation}
Therefore, since $\ti\F$, $\F$ are of ind-finite type and $\ti\F$ is reduced by
 [Fa1], the reduced locus of the neutral component $\F^0$ of $\F$
coincides with the scheme theoretic image $\phi(\ti \F)$. However, if
there is an element in
$ (R((t)))^\times/(R((t))^\times)^2R[[t]]^\times$   whose class in $Q^0(R)$
is nontrivial, then there is also an $R$-valued point of $\F^0$ which does not
belong to $\phi(\ti \F)(R)$. This shows that
$\F$ is not reduced. In fact, the point of $\F^0(R)$ represented by 
\begin{equation*}
 z=\left(\begin{matrix}
 1+at^{-1}& 0\cr
 0& 1
 \end{matrix}\right)\in PGL_2(R((t)))
 \end{equation*}
does not lie in $\F_{\rm red}$. (Here again $a\in R$ satisfies $a\neq 0$ and $a^2=0$.) The morphism $\ti\F\to \F$ induces a surjective radiciel morphism between $\ti\F=\F_{SL_2}$ and the
reduced locus of the neutral component of $\F=\F_{PGL_2}$.
We can see that this is ind-finite and so it is a universal homeomorphism (in this context of ind-schemes).
}\end{Remarknumb}

The assumption in Theorem
\ref{redu} that $G$ be semi-simple, is also necessary, as is shown by the next proposition.

\begin{prop}\label{nonsemi}
Let $G$ be a reductive group over $K$ which is not semi-simple. Then $LG$ and $\F_{\und{a}}$, for any facet $\und{a}$, are not reduced. In particular, if  $G = T$ is a non-trivial
torus, then $LT$ and $\F_{{\mathcal T}^0}$ are not reduced.
\end{prop}
\begin{Proof}
We first prove the last statement. 
Suppose that $T$ is a non-trivial torus over $K$.  Let $K'/K$ be a finite Galois splitting field for $T$ with Galois group $\Gamma$; denote by $\Gamma_t$ the (cyclic) tame quotient of $\Gamma$ by the wild inertia 
subgroup  $\Gamma_w\subset \Gamma$. We write $K'=k((u))$ and $K_t=(K')^{\Gamma_w}=k((v))$
with $v=t^{1/m}$, where $m=\#\Gamma_t$.
Then, for every $k$-algebra $R$,
we have
\begin{equation*}
LT(R)=T(R((t)))=(T(R((u))))^\Gamma=(X_*(T)\otimes_{\Z}R((u))^\times)^\Gamma
\end{equation*}
where $X_*(T)$ is   the $\Gamma$-module 
given
 by the cocharacter lattice. Consider the $k$-algebra $R=k[x]/(x^a)$, for  $a>1$
 and set
 \begin{equation*}
 M=\{1+\sum_{i=1}^{a-1}a_ix^iv^{-i}\ |\ a_i\in k\}\subset R((v))^\times\ .
 \end{equation*}
This is subgroup of $R((u))^\times $ which is fixed by $\Gamma$ and is killed by $p^b$ for some $b>>0$.
 Now find a non-zero
element $w_0\in X_*(T)\otimes_{\Z}\Z/p^b\Z$ which is fixed by the $p$-group $\Gamma_w$
and consider the subgroup $V$ generated by the ``conjugates"
$\gamma\cdot w_0$, $\gamma\in\Gamma$. The action of $\Gamma$ on 
$$
V\otimes_{\Z} M\subset (X_*(T)\otimes_{\Z}R((u))^\times)^{\Gamma_w}
$$
is inflated from an action of the cyclic group $\Gamma_t=\Gamma/\Gamma_w$ and there is a 
non-trivial element $z_0$  such that $\gamma_t\cdot z_0=\chi(\gamma_t) z_0$
for a character $\chi: \Gamma_t\to k^\times$. Fix $n_\chi>0$ such that $\gamma_t(v^{n_\chi})=\chi(\gamma_t)v^{n_\chi}$. 
Then the $R$-valued point $z_0\cdot v^{-n_\chi}$
 is fixed by $\Gamma$; hence it gives an $R$-valued point of $LT$ and a non-trivial point of $\F_{{\mathcal T}^0}$.
We  conclude that $\F_{{\mathcal T}^0}$ is not discrete and hence not reduced.
This implies that $LT$ is also not reduced.

Now we consider the general case. Let $G' = G_{\rm der}$ and $D = G/G'$. If $G$ is not semi-simple, then $D$
is a non-trivial torus and the previous considerations apply to $D$. There exists a commutative diagram
$$\xymatrix@C=25mm{
G\ar[r]^\pi & D  \\
\ & D\ar[ul]\ar[u]_\lambda
}$$
where $\pi$ is the natural projection and where $\lambda$ is an isogeny. This induces a commutative diagram of
homomorphisms of ind-group schemes,
$$\xymatrix@C=25mm{
LG\ar[r] & LD  \\
\ & LD\ar[ul]\ar[u]_{L\lambda}
}$$
If $LG$ were reduced, then $L\lambda$ would factor through $(LD)_{\rm red}$. But then there would exist
$N>0$ such that for any Artin local ring $R$ and any $t\in LD(R)$ we have $t^N\in (LD)_{\rm red} (R)$. However,
we can easily see that there exists a point of the above form   for suitable 
$R = k[x]/(x^a)$ such that its $N$-th power is a
non-trivial infinitesimal $R$-valued point of $LD$. \endproof
\end{Proof}

\subsection{}\label{GG'} Suppose that $G$ is connected reductive but not necessarily semi-simple over $K=k((t))$ with $k$  algebraically closed. Let us set $G'=G_{\rm der}$
and $D=G/G'$. Suppose that $\und a$ is a facet in the Bruhat-Tits building of $G(K)$. Denote by $P'_{\und a}$, $P_{\und a}$
the corresponding parahoric group schemes, and $\F'_{\und a}$, $\F_{\und a}$ the
corresponding affine flag varieties for $G'$, $G$ respectively. For simplicity, in what follows, we
will omit the subscript $\und a$.
Consider the exact sequence
\begin{equation}\label{der1}
1\to G'\xrightarrow{\phi} G\xrightarrow{\ep}  D\to 1\ .
\end{equation}
By [B-TII], 1.7, there are homomorphisms of group schemes $P \to {\mathcal D}^0$,
$P'\to P$ over $k[[t]]$; the composition $P'\to P\to {\mathcal D}^0$
is trivial. %Clearly $P'(R[[t]])\to P(R[[t]])$ is injective.
We obtain   morphisms of ind-schemes
\begin{equation}
 \F'\xrightarrow{\phi} \F\xrightarrow{\ep} \F_{{\mathcal D}^0}\ .
\end{equation}
Using [B-TII], 5.2.11, and Theorem \ref{components} and its proof, 
we see that if $k'\supset k$ is algebraically closed, we have
\begin{equation}
(LG)^0(k')=\phi((LG')^0(k')){\mathcal T}^0(k'[[t]]) \ ,
\end{equation}
where ${\mathcal T}^0$ is the connected Neron model of any maximal torus $T$ of $G$
and $(LG)^0$, $(LG')^0$ denote the neutral components.
We can deduce that $\phi: (\F')^0(k')\to \F^0(k')$ is surjective; by our description
of $P'(k'[[t]])$, resp.  $P(k'[[t]])$ as the stabilizers of $\und a$ in $(LG')^0(k')$, resp.  $(LG)^0(k')$,
it is also injective, hence bijective. We conclude that
\begin{equation}
(\F')^0_{\rm red}\xrightarrow{\phi} (\F)^0_{\rm red}
\end{equation}
is a surjective radiciel morphism between the reduced loci of the neutral components
of $\F'$ and $\F$. 

We do not know whether this morphism is an isomorphism in general. In the tamely ramified case
there is the following result.

\begin{prop}\label{GG'prop}
Suppose that $G$ splits over a tamely ramified extension of $K$.
Then
\begin{equation}
(\F_{G_{\rm der}})^0_{\rm red}\xrightarrow{\phi} (\F_G)^0_{\rm red}
\end{equation} 
is an isomorphism.
\end{prop}

\begin{Proof} 
The essential ingredient is the following:

\begin{lemma}\label{NeronExact}
Let $1\to T'\to T\to T''\to 1$ be an exact sequence of tori
over $K=k((t))$ such that %$T' $  is induced and 
$T$ splits over a tamely ramified extension of $K$.
Then there is an exact sequence
of smooth group schemes over $\O=k[[t]]$
\begin{equation}\label{neronseq}
1\to  {\mathcal S}'\to  {\mathcal T}^0\xrightarrow{}  {\mathcal T}''^0\to 1\  
\end{equation}
which extends the above sequence of tori and is such that:
\begin{itemize}
\item ${\mathcal T}^0$, ${\mathcal T}''^0$ are the connected Neron models 
of $T$, $T''$ respectively.

\item The connected component of the group scheme ${\cal S}'$ is 
equal to ${\cal T}'^0$ and the group of connected components 
of the special fiber of ${\cal S}'$ is a subgroup of 
the torsion subgroup of  $X_*(T')_I$.
\end{itemize}
In fact, if $T'$ is induced, then $X_*(T')_I$ is torsion-free and so ${\cal S}'={\cal T}'^0$.
\end{lemma}

\begin{Proof}
Assume that $T$ splits over the tamely ramified (finite) Galois extension
$\ti K/K$ with Galois group $\Gamma$. Then the same is true for $T'$ and  $T$. 
Since we are then dealing with split tori, there is a
short exact sequence
\begin{equation}
1\to {\mathcal T}'_{\ti\O}\to {\mathcal T}_{\ti\O}\to {\mathcal T}''_{\ti\O}\to 1
\end{equation}
between lft N\'eron models  over $\ti\O$.
As in [Ed], proof of Theorem 4.2, we see that the fixed point schemes $({\mathcal T}'_{\ti\O})^\Gamma$,
$({\mathcal T}_{\ti\O})^\Gamma$, $({\mathcal T}''_{\ti\O})^\Gamma$
are the lft  N\'eron models ${\mathcal T}'$,
${\mathcal T}$, ${\mathcal T}''$ for $T'$, $T$, $T''$ over $\O$.
Hence, there is an exact sequence
\begin{equation}\label{neronseq2}
1\to {\mathcal T}' \to {\mathcal T} \to {\mathcal T}''\ .   
\end{equation}
In particular, ${\mathcal T}'$ is a closed subgroup scheme of ${\mathcal T}$.
 Hilbert's theorem 90 and the triviality of the Brauer group
of $K$ implies that
${\rm H}^1(K, T')=(0)$ (see e.g [C] Lemma 4.3) and so $1\to T'(K)\to T(K)\to T''(K)\to 1$
is exact. Since by the N\'eron extension property ${\mathcal T}(\O)=T(K)$, 
${\mathcal T}'(\O)=T'(K)$, ${\mathcal T}''(\O)=T''(K)$,
we conclude that 
\begin{equation}\label{neronseq4}
1\to {\mathcal T}'(\O) \to {\mathcal T}(\O) \to {\mathcal T}''(\O)\to 1 .   
\end{equation}
is also exact. On the other hand, the targets of the Kottwitz invariant homomorphisms
form the exact sequence
\begin{equation}
 X_*(T')_I\xrightarrow{\phi_*} X_*(T)_I\to X_*(T'')_I\to 0\ \ ,
\end{equation}
where ${\rm ker}(\phi_*)\subset (X_*(T')_I)_{\rm tor}$.
Now define ${\cal S}'$ to be the unique subgroup scheme of 
the lft N\'eron model
${\cal T}'$ which contains ${\cal T}'^0$
and with special fiber given by the connected components
of the special fiber of ${\cal T}'$ which are parametrized by 
 ${\rm ker}(\phi_*)$.  We can conclude by the snake lemma 
and Theorem \ref{components} that
\begin{equation}\label{neronseq3}
1\to {\mathcal S}'(\O) \to  {\mathcal T}^0(\O) \to  {\mathcal T}''^0(\O) \to 1\  
\end{equation}
is exact. (The exactness on the right follows also using [BLR], 9.6.2.)
By the above, ${\mathcal S}'$ is a closed subgroup scheme of ${\mathcal T}^0$.
The claim that (\ref{neronseq}) is exact everywhere except
at ${\mathcal T}''^0$ immediately follows from (\ref{neronseq2}) and (\ref{neronseq3}). 
The claim that  (\ref{neronseq})
is exact at ${\mathcal T}''^0$ also follows from the above and [B-TII], Prop. 1.7.6.
Indeed, both ${\mathcal T}''^0$ and the quotient ${\mathcal T}^0/{\mathcal S}'$ 
are smooth (affine) commutative group schemes with $T''$ as generic fiber and
with the same $\O$-valued points. \endproof
\end{Proof}

\smallskip

We will apply Lemma \ref{NeronExact} to the exact sequence $1\to T'\to T\to D\to 1$
where $T'$, $T$ are maximal tori in $G'$, $G$ respectively. Let $T'_{\rm sc}$ be the inverse image
of $T'$ in the simply connected cover of $G'$.
Consider
\begin{equation*}
0\to X_*(T'_{\rm sc})\to X_*(T')\to \pi_1(G')\to 0,\quad 0\to X_*(T'_{{\rm sc}})_I\to X_*(T')_I\to \pi_1(G')_I\to 0\ .
\end{equation*}
Here the first sequence is exact and $\pi_1(G')$ is a finite group.  Also $X_*(T'_{\rm sc})$ is an induced $I$-lattice
(the coroot lattice)
so $X_*(T'_{\rm sc})_I=X_*(T'_{\rm sc})^I$ is torsion-free.
Hence, the second sequence is also exact and we obtain
\begin{equation}\label{torpi1}
(X_*(T')_I)_{\rm tor}\hookrightarrow \pi_1(G')_I\ .
\end{equation}

To complete the proof of Proposition \ref{GG'prop}, suppose that $R$ is a  
strictly henselian local $k$-algebra; then
so is $R[[t]]$. Since ${\mathcal S}'$ is smooth over $k[[t]]$ the
homomorphism ${\mathcal T}^0(R[[t]])\to {\mathcal D}^0(R[[t]])$ is surjective.
Now suppose that we have $[g]\in LG^0(R)/L^+P(R)=\F^0(R)$ and that $R$ is in addition reduced. Then
by Corollary \ref{torusDiscrete} we have $\ep(g)\in {\mathcal D}^0(R[[t]])$.
By modifying $g$ by multiplying with the inverse of a lift of $\ep(g)$
in ${\mathcal T}^0(R[[t]])\subset P(R[[t]])$ we can assume that $\ep(g)=1$, therefore
$g\in LG'(R)$. In fact, $\kappa_{G'}(g)$ lies in the finite subgroup of $(X_*(T')_I)_{\rm tor}\subset 
\pi_1(G')_I$
given by Lemma \ref{NeronExact} and (\ref{torpi1}). By modifying $g$ further via multiplication
by an element $s\in {\cal S}'(R[[t]])\subset T'(R((t)))$ with $\kappa_{G'}(s)=-\kappa_{G'}(g)$ we can arrange
so that the new $g$ is also in $(LG')^0(R)$.
This shows that, when $G$ splits over a tame extension,
then for every {\sl reduced} strictly local $k$-algebra $R$,
the map $\phi: \F'^0(R)\to \F^0(R)$ is surjective. 
Now $(\F')^0_{\rm red}$ and $(\F)^0_{\rm red}$ are of ind-finite type.
If $\phi: X \to Y$ is a morphism between reduced $k$-schemes of finite type, which is a universal
homeomorphism and such that $ X(R) \to Y(R)$ is surjective for every reduced strictly local $k$-algebra $R$,
then $\phi$ is an isomorphism. This, combined with the above, allows us to conclude the proof.
\endproof
\end{Proof}

 \bigskip
\bigskip

\section{Tame forms and liftings}\label{tamepar}

\setcounter{equation}{0}

\subsection{}\label{tame1} In this section, we assume that $k$ is algebraically closed.
Suppose that $G$ is an absolutely simple, simply connected,  quasi-split
group over $K=k((t))$ that splits over the (totally ramified) separable extension $K'/K$ of degree $e$. This degree can take the values
$e=1$ (when $G$ is split over $K$), $2$ or $3$ (see [T1], \S4, or [B-TII], \S 4.1-4.2). We have
$e=2$ when the local Dynkin diagram of $G$ is of type ${\rm B}-{\rm C}_n$, (for $n\geq 3$),
${\rm C}-{\rm B}_n$, (for $n\geq 2$), ${\rm C}-{\rm BC}_n$ (for $n\geq 1$), ${\rm F}^I_4$. We have $e=3$ only in
one case, when the Dynkin diagram is of the form ${\rm G}^I_2$. (The notations
for the Dynkin diagrams are as in [T1], Table 4.2, p. 60).

We assume that $(e, {\rm char}(k))=1$ so that $K'/K$ is tame and we  set $K'=k((u))$ with $u^e=t$.  Fix
a generator $\tau$ of ${\rm Gal}(K'/K)$.
Let us denote by $H$ the {\sl split} simple, simply connected,  algebraic group over $K$ such that
 $G_{K'} =H_{K'}$ (we fix a choice of this isomorphism).
Consider the Weil restriction of scalars $H'={\rm Res}_{K'/K}(H_{K'})$; the automorphism
${\rm id}\otimes\tau$ of ${\rm Res}_{K'/K}(G\otimes_KK')= H'$
induces an automorphism $\sigma$ of $H'$ which can be written
on $H_{K'}=H\otimes_KK'$ in the form $\sigma=\sigma_0\otimes\tau$, where $\sigma_0$ is a Dynkin diagram
automorphism of order $e$ of the group $H$. In the cases ${\rm B}-{\rm C}_n$, ${\rm C}-{\rm BC}_n$ (these correspond
to $G$ a special ``ramified" unitary group as in \S \ref{unitary}), we have $H=SL_{2n}$,
or $H=SL_{2n+1}$ respectively; the involution $\sigma_0$ is the usual transpose inverse involution.
In the case ${\rm C}-{\rm B}_n$, we have $H={\rm Spin}_{2n+2}$. In the case ${\rm F}^I_4$, we have $H={\rm E}_6$.
Finally, in the case ${\rm G}^I_2$ (then $e=3$), we have $H={\rm Spin}_8$ (``triality").

The automorphism $\sigma$ gives an action on the Bruhat-Tits building
${\cal B}(H', K)={\cal B}(H, K')$. By the main result of [P-Y], the set of
points ${\cal B}(H, K')^\sigma $ of ${\cal B}(H, K')$ fixed by $\sigma$ can be identified ($G(K)$-equivariantly)
with the building ${\cal B}(G, K)$. Now let $x$ be a point in ${\cal B}(G, K)\subset {\cal B}(H, K')$.
Since $H$ is simply connected, the stabilizer $H(K')^x$ of $x$ in $H(K')$ gives the $k[[u]]$-valued points
of the corresponding parahoric subgroup scheme $P^{H}_x$ over $k[[u]]$. One can see
 that
the automorphism $\sigma$ extends to ${\rm Res}_{k[[u]]/k[[t]]}(P^{H}_x)$ (e.g using [B-TII], 1.7). Now let us consider the scheme of fixed points
$({\rm Res}_{k[[u]]/k[[t]]}(P^{H}_x))^\sigma$. Since $\sigma$ has order prime to the characteristic $p$, the scheme
$({\rm Res}_{k[[u]]/k[[t]]}(P^{H}_x))^\sigma$ is a group scheme which is smooth over $k[[t]]$. In fact, since its
$k[[t]]$-valued points give the stabilizer $G(K)^x$ and $G$ is simply connected, we have
\begin{equation}
P^G_x=({\rm Res}_{k[[u]]/k[[t]]}(P^{H}_x))^\sigma
\end{equation}
with $P^G_x$ the parahoric group scheme over $k[[t]]$ associated to $x\in {\cal B}(G, K)$.

Suppose now that in addition ${\rm char}(k)=p>0$; we denote by $W=W(k)$ the ring of Witt vectors of $k$.
Our goal is to extend this picture to a mixed characteristic situation (i.e., over $W$).
We lift $k[[u]]/k[[t]]$ to the ring extension $W[[u]]/W[[t]]$ with $u^e=t$; we will continue to denote
by $\tau$ the automorphism of $W[[u]]$ with lifts our choice of Galois generator
with the same notation. First we observe that the group scheme $P^H_x$ over $k[[u]]$
has a natural lift to a group scheme $\und P^H_x$ over $W[[u]]$. This can be obtained
by applying [B-TII], \S 3.2 to the base ring $A=W[[u]]$ and to   schematic root data
for the split group $H$ that correspond to $x$. More precisely, the needed ideals $f(a)\subset W[[u]]$ in
loc. cit. Prop. 3.2.7, are given by $(u^{m(a)})$, where $m(a)\in \Z$ and where $(u^{m(a)})\subset k[[u]]$
are the corresponding ideals describing the  schematic root data for $P^{H}_x$ that
correspond to $x\in {\cal B}(H, K')$; this choice guarantees that our lift is ``horizontal in the $W$-direction".
We can see again that we have an automorphism $\sigma$ on $({\rm Res}_{W[[u]]/W[[t]]}(\und P^{H}_x))$
that covers $\tau: W[[u]]\to W[[u]]$. We set
\begin{equation}
\und P^G_x=({\rm Res}_{W[[u]]/W[[t]]}(\und P^{H}_x))^\sigma\ ;
\end{equation}
this is a smooth group scheme over $W[[t]]$ that lifts the parahoric group scheme $P^G_x$ over $k[[t]]$.
In fact, one could also construct $\und P^G_x$ directly by applying [B-TII], 3.9.4 to $A=W[[t]]$
and to suitable ``$W$-horizontal" schematic root data given by $x\in {\cal B}(G, K)$.

Denote by $H_0$ the Chevalley group scheme over $W$ that corresponds to
the split group $H$. Let us set $\und H=H_0\otimes_W{W((u))}$.
This is a smooth affine group
scheme over $W((u))$ that lifts $H$ over $K'=k((u))$. The automorphism $\sigma$ lifts to an automorphism
of $\und H'={\rm Res}_{W((u))/W((t))}(\und H)$ over $W((t))$ and we can set
\begin{equation}
\und G:=({\rm Res}_{W((u))/W((t))}(\und H))^\sigma\ .
\end{equation}
(The group scheme $\und G$  over $W((t))$ can also be obtained as in [T2], p. 217).

\subsection{} We continue with the above assumptions and notations.
In particular, we fix $x\in {\cal B}(G, K)\subset {\cal B}(H, K')$.
Let us consider the following functors
of $W$-algebras:
\begin{equation*}
L\und H(R):= \und H(R((u))), \quad L^+ \und P^H_x(R): =\und P_x^H(R[[u]])\ .
\end{equation*}
Also we have
\begin{equation*}
L\und G(R)=(\und H(R((u)))^\sigma,\quad
L^+\und P^G_x(R) =(\und P_x^H(R[[u]]))^\sigma\subset L\und G(R)\ .
\end{equation*}

Consider the fpqc sheaf  $L\und H /L^+ \und P^H_x$ over $W$ which is associated to the presheaf
$R\mapsto L\und H(R)/L^+ \und P^H_x(R)$. Since $\und H$ is split, we may refer to [Fa1], section 1, for the fact that this sheaf  is represented by an ind-scheme
over $W$ (which is actually ind-proper); it is a (partial) affine flag variety  for the (untwisted) loop group scheme $\und H$ over $W$.

\begin{prop}\label{fixed}
Consider the fpqc quotient sheaf $L\und G/L^+\und P^G_x$ over $W$.

a)The fpqc sheaf $L\und G/L^+\und P^G_x$ can be identified with the neutral component of
the fixed point ind-scheme $(L\und H /L^+ \und P^H_x)^\sigma$.

b) $L\und G/L^+\und P^G_x$ is representable by an ind-scheme which is ind-proper over $W$.

c) The fiber of $L\und G/L^+\und P^G_x$ over $k$ is isomorphic to the affine flag variety $LG/L^+P^G_x$
for $G$ and the parahoric subgroup scheme $P^G_x$.
\end{prop}

\begin{Proof}
Let us first discuss part (a). We first assume that $P^H_x$ is an Iwahori
group scheme for $H$; recall that $x$ is fixed by $\sigma$. Then the maximal reductive quotient of the reduction of $\und P^H_x$
over $(u)$ is a split torus $T_0$ over $W$ which is  $\sigma$-stable.
Let $R$ be a local $W$-algebra. Recall that by [Fa1], the quotient morphism $L\und H\to L\und H/L^+\und P^H_x$
splits locally for the Zariski topology of the ind-scheme $L\und H/L^+\und P^H_x$
and so we have
\begin{equation*}
(L\und H/\und P^H_x)(R)=L\und H(R)/L^+\und P^H_x(R).
\end{equation*}
There is a natural ``exact" sequence of pointed sets
\begin{equation}\label{pointed}
1\to (L\und H(R)^\sigma/(L^+\und P^H_x(R))^\sigma\to (L\und H(R)/L^+\und P^H_x(R))^\sigma\to
{\rm H}^1(\langle\sigma\rangle, \und P^H_x(R[[u]]))
\end{equation}
Now denote by $P^+(R)$ the kernel of the (surjective) reduction homomorphism
\begin{equation*}
\und P^H_x(R[[u]])\to  T_0(R)\ .
\end{equation*}
We can see that ${\rm H}^1(\langle\sigma\rangle, \und P^+(R))=(0)$
since $\langle\sigma\rangle$ has order prime to $p$ and $P^+(R)$ is a projective
limit of a system of groups which are iterated extensions of $W$-modules.
Observe that the action of $\tau$ on $T_0(R)$ is trivial. Also recall that $\sigma_0$
 is a diagram automorphism and hence acts on the cocharacter lattice $X_*(T_0)$
 of the torus $T_0$ by a permutation representation. Therefore, we have
\begin{equation}
{\rm H}^1(\langle\sigma\rangle, \und P^H_x(R[[u]]))={\rm H}^1(\langle\sigma_0\rangle, T_0(R))={\rm H}^1(\langle\sigma_0\rangle,  R^\times\otimes_\Z X_*(T_0)).
\end{equation}
We can now see, (the group $\langle \sigma_0\rangle$ is cyclic of order $e$), that
${\rm H}^1(\langle\sigma_0\rangle,  R^\times\otimes_\Z X_*(T_0))=
{\rm H}^1(\langle\sigma_0\rangle,  k^\times\otimes_\Z X_*(T_0))$, i.e.,
it is independent of  $R$, and is a finite torsion group which we denote by $C(\langle\sigma_0\rangle, H)$.
We obtain a map
\begin{equation}
(L\und H/L^+\und P^H_x)^\sigma\to
C(\langle\sigma_0\rangle, H)
\end{equation}
which is locally constant for the Zariski topology; this fact together with (\ref{pointed})
completes the proof of (a) when $P^H_x$ is Iwahori. To deal with the general case,
suppose that $x$ is a fixed point of $\sigma$ in the building ${\cal B}(H, K')$. We claim
that there is a $\sigma$-fixed point $z\in {\cal B}(H, K')$ such that $P^H_{z}$ is Iwahori and
$\und P^H_z(R[[u]])\subset \und P^H_x(R[[u]])$: Indeed, consider the maximal reductive quotient
of the reduction of $\und P^H_x$ modulo $(u)$. (This is a reductive
group $\und Q$ over $W$). We can see that this supports an action of
$\sigma_0$ which is given by a diagram automorphism. Hence we can find a
$\sigma_0$-stable Borel subgroup $\und B_Q\subset \und Q$. This gives
a $\sigma$-stable Iwahori subgroup $P^H_z$ whose lifting $\und P^H_z$
has the desired property.  Now observe that
\begin{equation*}
L\und H/L^+\und P^H_z\to L\und H/L^+\und P^H_x
\end{equation*}
is a $\sigma$-equivariant (Zariski locally trivial) smooth fibration with fiber over the origin the homogeneous
space $\und P^H_x/\und P^H_z=\und Q/\und B_Q$. Note that all 
the fibers
of the fixed point set
$(\und Q/\und B_Q)^\sigma$ over $W$ are non-empty. 
Indeed, it suffices by the Lefschetz fixed point formula to see that the trace of the automorphism $\sigma_0$ on the $\ell$-adic cohomology of the flag variety $Q/B_Q$ is non-zero. But this cohomology has a natural basis enumerated by the elements of the Weyl group of $Q$, and $\sigma_0$ acts by permuting the basis vectors. Consider the representation of the cyclic group generated by $\sigma_0$. On an irreducible constitutent associated to a fixed vector, the trace of $\sigma_0$ is $1$, on all other irreducible constituents it is $0$. Since there are fixed vectors (e.~g., the vector associated to the unit element in the Weyl group), the trace is non-zero. 
We conclude that
\begin{equation*}
(L\und H/L^+\und P^H_z)^\sigma\to (L\und H/L^+\und P^H_x)^\sigma
\end{equation*}
is surjective. The proof of (a) in the general case now follows from these considerations and
Theorem \ref{components}.

Part (b) now follows from (a) and the corresponding properties of $L\und H/L^+\und P^H_x$
since the action of $\sigma$ respects the ind-scheme structure of $L\und H/L^+\und P^H_x$.
In fact, we also see from (a) and its proof that for any  local $W$-algebra $R$, we have
\begin{equation}
L\und G(R)/L^+\und P^G_x(R)=(L\und G/L^+\und P^G_x)(R) .
\end{equation}
Part (c) follows now from the above and our construction.
\end{Proof}

\begin{Remarknumb}
{\rm   Let us consider the example of the ``ramified" special unitary group as in \S \ref{unitary}.
Take $H=SL_n$, with $\sigma_0$  the involution
\begin{equation*}
g\mapsto L_n\cdot (g^{\rm tr})^{-1}\cdot L^{-1}_n
\end{equation*}
where $L_n$ is the antidiagonal matrix of size $n$
and $\tau(u)=\bar u=-u$. Then, $G=SU_{n}$ as in \S \ref{unitary}.
In this case, we can obtain the group $\und G$
over $W((t))$ as follows: Consider the $W[[u]]$-lattice
\begin{equation*}
\Lambda_0={\rm span}_{W[[u]]}\{e_1,\ldots ,e_n\}\
\end{equation*}
with the (perfect) $W[[u]]$-valued hermitian form $\Phi$ defined by $\Phi(e_i, e_j)=\delta_{i, n+1-j}$.
First consider the group scheme $\und P_0$ over $W[[t]]$ whose $S$-valued points
are the $W[[u]]\otimes_{W[[t]]}\O_S$-linear automorphisms $A$ of $\Lambda_0\otimes_{W[[t]]}\O_S$ with determinant $1$
which respect the form $\Phi$, i.e., such that $\Phi(A\cdot v, A\cdot w)=\Phi(v, w)$ for all $v$, $w\in \Lambda_0\otimes_{W[[t]]}\O_S$.
Then set $\und G=\und P_0\otimes_{W[[t]]}W((t))$. In fact, we can  similarly also give $W[[u]]$-lattices $\Lambda_i$,  $\Lambda_{m'}$ by lifting the construction of the lattices $\lambda_i$, $\lambda_{m'}$ in \S \ref{unitary}
and consider group schemes $\und P'_I$ over $W[[t]]$ defined via $\Phi$-unitary automorphisms
that stabilize suitable chains of the $\Lambda$'s.
The group schemes $\und P'_I$ give explicit realizations of the
group schemes $\und P^G_x$ constructed above. They
lift the parahoric group schemes $P'_I$ over $k[[t]]$ which
were constructed in \S \ref{unitary}.

Let us comment on the statement  of Proposition \ref{fixed}  (a) in this special case.
For simplicity, we just restrict our consideration to the special fiber over $k$
and to the case of an Iwahori subgroup.

i) $n=2m+1\geq 3$ odd. Take $I=\{0,\ldots, m\}$. Then
$P'_I$ is an Iwahori subgroup of $SU_n(K)$. In this case, by Remark \ref{uniRem} (B), which identifies 
$LSU_n/L^+P'_I$ with the functor of complete selfdual lattice chains of the correct determinant,
we have
\begin{equation}
LSU_n/L^+P'_I=(LSL_n/L^+I_{SL_n})^\sigma
\end{equation}
where $I_{SL_n}$ is a corresponding Iwahori of $SL_n$.

ii) $n=2m\geq 4$ even. Take $I=\{0,\ldots , m-2, m, m'\}$
so that $P'_I$ is an Iwahori subgroup of $SU_n(K)$. By Remark
\ref{uniRem} (B) we have an ``exact" sequence of pointed sets
\begin{equation}
1\to LSU_n/L^+P'_I\to (LSL_n/L^+I_{SL_n})^\sigma\xrightarrow{\iota} \Z/2\Z\to 0\ .
\end{equation}
Again $(LSL_n/L^+I_{SL_n})^\sigma$ classifies complete self-dual lattice chains of the correct determinant
and the map $\iota$  is given by
the sign in the identity $a=\pm \bar b\cdot b^{-1}$ in loc. cit.}
\end{Remarknumb}

\subsection{} \label{KacM1}
 We continue with the above assumptions and notations, in particular $G$ splits over the 
 totally tamely ramified extension $K'/K$ of degree $e=1, 2$ or $3$.
 Let $\Delta_G$ be the local Dynkin diagram of the quasi-split group $G$ over the local field $K=k((t))$
(see for example [T1] 1.8 and table 4.2). Each of these local Dynkin diagrams
is (up to numbering) the Dynkin diagram of a unique  generalized Cartan matrix of 
affine type in the sense of Kac 
(see [Kac],  4.8, in particular the tables Aff 1,2,3 on p. 54-55). 

Now suppose that $F$ is a field that contains $W=W(k)$. We will denote by 
$\gfr_{KM}(\Delta_G)_F$ or simply $\gfr_{KM}(G)_F$ the {\sl affine} Kac-Moody Lie algebra over $F$
given by the corresponding generalized Cartan matrix $A(G)=\{a_{ij}\}_{1\leq i,j \leq l+1}$ ([Kac] Chapter 1).
Here $l$ is the rank of the split group $G_{K'}$. 
In fact, we will only need to consider the 
factor algebra $\gfr'_{KM}(G)_F$ of $\gfr_{KM}(G)_F$ by the one 
dimensional subspace generated by the derivation element.
The Lie algebra $\gfr'_{KM}(G)_F$ is generated over $F$ by $3(l+1)$ generators
$e_i$, $f_i$, $h_i$ with the relations:
$$
[h_i, h_j]=0,\quad [e_i, f_i]=h_i,\quad  [e_i, f_j]=0, \ \ \hbox{\rm if\ } i\neq j,
$$
$$
[h_i, e_j]=a_{ij}e_j,\quad  [h_i, f_j]=-a_{ij}f_j,
$$
$$
({\rm ad}\, e_i)^{1-a_{ij}}e_j=0,\quad ({\rm ad}\, f_i)^{1-a_{ij}}f_j=0, \ \ \hbox{\rm if\ } i\neq j.
$$

We also have the following description of $\gfr'_{KM}(G)_F$:

If the group $G$ is split, then $\und G=H_0\otimes_W W((t))$ and $\sigma={\rm id}$. Then
\begin{equation*}
 \gfr'_{KM}(G)_F\simeq F\cdot c\oplus  (Lie(H_0)\otimes_W F[t, t^{-1}])\ .
\end{equation*}
with bracket described in [Kac], 7.1. (The element $c$ is central, $(Lie(H_0)\otimes_WF[t, t^{-1}])$
is the loop Lie algebra of $H_0$.) This is an {\sl untwisted} affine Kac-Moody algebra
(its generalized Cartan matrix is from table Aff 1 of [Kac], 4.8). To explain the isomorphism, 
denote by $\Phi'$ the set of roots of the split group $G_{K'}=H_{K'}$ and fix a  choice of 
a set of simple roots $a'_1, \ldots, a'_l$ and hence a set of positive roots $\Phi'^+$. Denote by $\theta$ the highest root of $\Phi$.
Fix a Chevalley basis $E_a$, $a\in \Phi'$,  $H_i=H_{ a'_i}$, $1\leq i\leq l$, for the Lie algebra
of the group $H_0$ over $W$. The elements $e_i=E_{a'_i}$, $f_i=E_{-a'_i}$, $h_i=H_i$ for $1\leq i\leq l$,
$e_0=t\otimes E_{-\theta}$, $f_0=t^{-1}\otimes E_{\theta}$ and $h_0=\frac{2}{(\theta|\theta)}\cdot c-\theta^\vee$
in $F\cdot c\oplus  (Lie(H_0)\otimes_W F[t, t^{-1}])$ satisfy the defining relations for $\gfr'_{KM}(H)_F$
above (the notation is as in [Kac] \S 7.4.)

If the group $G$ is not split, then $\sigma=\sigma_0\otimes\tau$ is not trivial. Define
an automorphism $\ti \sigma$ of
$$
\gfr'_{KM}(H)_F=F\cdot c'\oplus (Lie(H_0)\otimes_WF[u, u^{-1}])
$$
by $\ti\sigma(c')=c'$,  while for $x\in Lie(H_0)$, $f(u)\in F[u, u^{-1}]$,
$\ti\sigma (x\otimes f(u))=\sigma_0(x)\otimes f(\tau(u))$. Then
\begin{equation*}
\gfr'_{KM}(G)_F\simeq (\gfr'_{KM}(H)_F)^{\ti \sigma},
\end{equation*}
the fixed point set of $\ti\sigma$ on $\gfr'_{KM}(H)_F$ (see [Kac], Theorem 8.3).
This is then a {\sl twisted} affine Kac-Moody algebra
(its generalized Cartan matrix is from table Aff 2 or Aff 3 of [Kac], 4.8;
in each case $2$ or $3$ is the order of $\sigma$.) 
In this case, the elements $e_i$, $f_i$, $h_i$ of $\gfr'_{KM}(H)_F$
can also be given explicitly from the Chevalley
generators $E_a$, $H_1,\ldots , H_l$, of the group $H_0$ over $W$ using the formulas in [Kac] \S 8.3.
 
Note that  there  is a Lie algebra embedding
\begin{equation*}
(Lie(H_0)\otimes_WF[u, u^{-1}])^{\ti \sigma}\hookrightarrow Lie(L(\und G\otimes_W F))\ .
\end{equation*}
(The Lie algebra $Lie(L(\und G\otimes_W F))$ also includes ``power series" in   $t$.)

\bigskip
\bigskip

\section{Weyl groups and Schubert varieties} \label{weyl}

\setcounter{equation}{0}

In this section, we assume that $k$ is algebraically closed and that
  $G$ is a connected reductive group over $K=k((t))$.

\subsection{}\label{affineWeylgrp} We start by recalling some  facts on affine Weyl groups (cf. [T1], [B-TI], [B-TII],  [R], and especially the appendix [H-R] to this paper).
Let $S$ be a maximal split torus in $G$ and let $T$ be its centralizer.
Since $k$ is algebraically closed, $G$ is quasi-split and so $T$ is a
maximal torus in $G$. Let $N=N(T)$ be the normalizer of $T$; denote by
$T(K)_1 $ the kernel of the Kottwitz homomorphism
$\kappa_T: T(K)\to X_*(T)_I$; then $T(K)_1={\mathcal T}^0(k[[t]])$
where ${\mathcal T}^0$ the connected Neron model of $T$ over $k[[t]]$.
By definition, the {\sl Iwahori-Weyl group associated to $S$} is the
quotient group
\begin{equation*}
\ti W=N(K)/T(K)_1\ ,
\end{equation*}
cf. [H-R], Def. 7. Since $\kappa_T$ is surjective, the Iwahori-Weyl group $\ti W$
is an extension of the relative Weyl group $W_0=N(K)/T(K)$ by $X_*(T)_I$:
\begin{equation}\label{exactWeyl}
0\to X_*(T)_I\to \ti W\to W_0\to 1 .
\end{equation}
We have ([H-R], Prop. 8, comp.~also [R])

\begin{prop}
Let $B_0$ be the Iwahori subgroup of $G(K)$ associated to an alcove contained
in the apartment associated to the maximal split torus $S$. Then $G(K)=B_0\cdot N(K)\cdot B_0$
and the map $B_0\cdot n\cdot B_0\mapsto n\in \ti W$ induces a bijection
\begin{equation}
B_0\backslash G(K)/B_0\xrightarrow{\sim} \ti W\ .
\end{equation}
If $P$ is the parahoric subgroup of $G(K)$ associated to a facet contained in the apartment corresponding to $S$, then
\begin{equation}
P \backslash G(K)/P \xrightarrow{\sim} \ti W^P\backslash \ti W/\ti W^P, \ \ \hbox{\rm where\ \ \ } \ti W^P:=(N(K)\cap P)/T(K)_1\ .
\end{equation}
\end{prop}

In fact, if $P$ is the (special) parahoric subgroup $P_{\und x}$ that corresponds to a special vertex $\und x$ in
the apartment corresponding to $S$, then the subgroup $\ti W^P\subset \ti W$ maps isomorphically to $W_0$ under
the quotient $\ti W\to W_0$ and the exact sequence (\ref{exactWeyl}) represents the Iwahori-Weyl group as
a semidirect product
\begin{equation}\label{semidirect}
\ti W=W_0\ltimes X_*(T)_I\ ,
\end{equation}
see [H-R], Prop. 13.

Now let $S_{{\rm sc}}$,   $T_{{\rm sc}}$, resp. $N_{{\rm sc}}$ be the inverse images
of $S\cap G_{{\rm der}}$, $T\cap G_{{\rm der}}$, resp. $N\cap G_{{\rm der}}$ in the simply connected covering
$G_{{\rm sc}}$ of the derived group $G_{{\rm der}}$. Then $S_{{\rm sc}}$ is a maximal split torus of $G_{{\rm sc}}$  and $T_{{\rm sc}}$, resp. $N_{{\rm sc}}$ is its centralizer, resp. normalizer. Hence
\begin{equation*}
W_a:=N_{{\rm sc}}(K)/T_{{\rm sc}}(K)_1
\end{equation*}
is the Iwahori-Weyl group of $G_{{\rm sc}}$. This group is also called the {\sl affine Weyl group
associated to $S$} and is a Coxeter group.
Indeed, we can  recover $W_a$ in the following way:
Let $N(K)_1$ be the intersection of $N(K)$ with the kernel $G(K)_1$ of the Kottwitz homomorphism
$\kappa_G: G(K)\to \pi_1(G)_I$. Then  the natural homomorphism
\begin{equation}
W_a=N_{{\rm sc}}(K)/T_{{\rm sc}}(K)_1\xrightarrow{\sim} N(K)_1/T(K)_1
\end{equation}
is an isomorphism ([H-R], Lemma 17) and  there is  an exact sequence
\begin{equation}
1\to W_a\to \ti W\xrightarrow{\kappa_G}  \pi_1(G)_I=X_*(T)_I/X_*(T_{{\rm sc}})_I\to 1\ ,
\end{equation}
cf. [H-R], Lemma 14. 
Now let $B_0$ be the Iwahori subgroup of $G(K)$ associated to an alcove  $C$ in the apartment corresponding to $S$
and let ${\bf S}$ be the set of reflections about the walls of $C$. Then by [B-TII], 5.2.12, the quadruple
$(G(K)_1, B_0, N(K)_1, {\bf S})$ is a double Tits system and $W_a=N(K)_1/T(K)_1$
is the affine Weyl group of the affine root system $\Phi_a$ of $S$. The affine Weyl group $W_a$
acts simply transitively on the set of alcoves in the apartment of $S$. Since $\ti W$ acts transitively
on the set of these chambers, $\ti W$ is the semi-direct product of $W_a$ with the normalizer $\Omega$ of
the base alcove $C$, i.e., the subgroup of $\ti W$ which preserves the alcove,
\begin{equation}
\ti W=W_a\rtimes \Omega\ .
\end{equation}
We can identify $\Omega$ with $X_*(T)_I/X_*(T_{{\rm sc}})_I=\pi_1(G)_I$, cf. [H-R], Lemma 14.

Let us write
\footnote{Here we follow the tradition of parametrizing the set of simple affine roots 
by $I$. This should not be confused with the notation for the inertia group.} 
${\bf S}=\{s_i\}_{i\in I}\subset W_a$ for the finite set of reflections about the walls of
$C$ that generate the Coxeter group $W_a$.
For each $w\in W_a$
its length $l(w)$ is the minimal number of factors in a product of $s_i$'s representing $w$.
Any such product realizing the minimum is called a reduced decomposition of $w$.
We will denote by $\leq $ the corresponding Bruhat order. Recall its definition:  Fix a reduced decomposition
of $w\in W_a$. The elements $w'\leq w$
are obtained by replacing some factors in it by $1$. (This set of such $w'$'s
is independent of the choice of the reduced decomposition of $w$.)

Let us denote by $\alpha_i\in \Phi_a$ the
unique affine root with corresponding affine reflection
equal to $s_i$. (Since the residue field $k$ is algebraically closed, $\frac{1}{2}\alpha_i\not\in \Phi_a$;
see [T1], 1.8.)
 We will denote by $\Delta=\Delta_G$ the (local) Dynkin diagram of the affine root system
$\Phi_a$. (This can be obtained from the subset $\{\alpha_i\}_{i\in I}$;
see [B-TI], 1.4 and [T1], 1.8.) Suppose that the Coxeter system $(W_a, {\bf S})$ is a product
of the irreducible systems $(W^k_a, {\bf S}^k)$ which correspond to
the  simple factors of the derived group $G_{\rm der}$; we have ${\bf S}=\cup_{k}{\bf S}^k$
([B-TI], 1.3).
 For  a  subset $Y\subset {\bf S}$ such that $Y\cap {\bf S}^k\neq {\bf S}^k$
 for all $k$ , we denote by $W_Y\subset W_a$ the
subgroup generated by $s_i$ with $i\in Y$; we set $P_Y=B_0\cdot W_Y\cdot B_0 $.
By  general properties of Tits systems these are subgroups of $G(K)_1\subset G(K)$; by [B-TII], 5.2.12 (i)
they are the parahoric subgroups of $G(K)$ that contain $B_0$. Using [B-TI], 1.3.5 we see that
we can identify $P_Y$ with the parahoric subgroup $P_{C_Y}$ where $C_Y$
is the facet consisting of $a$ in the closure $ \overline C$ of $C$ for which $Y$ is exactly
the set of reflections $s\in {\bf S}$ which fix $a$.

Now let us fix a special vertex $\und x$ in the apartment
corresponding to $S$. We may assume that $\und x$ is in  $\overline C$.
Then there exists a reduced root system ${}^x\Sigma$ such that the semi-direct
product (\ref{semidirect}) (for $G_{{\rm sc}}$ instead of $G$) presents $W_a$ as the affine Weyl group
associated (in the sense of Bourbaki) to ${}^x\Sigma$,
\begin{equation}\label{semidirect2}
W_a=W({}^x\Sigma)\ltimes Q^\vee ({}^x\Sigma) ,
\end{equation}
(cf. [B-TI],~1.3.8,  [T1],~1.7, 1.9). In other words, we have   identifications $W_0\simeq W({}^x\Sigma)$,
$X_*(T_{{\rm sc}})_I\simeq Q^\vee({}^x\Sigma)$ compatible with the semidirect product
decompositions (\ref{semidirect}) and (\ref{semidirect2}), cf. [H-R], Prop. 13.

\begin{Remarknumb}
{\rm
The following variant of the Iwahori-Weyl group appears in [T1], p.~34. Set $W^\flat=N(K)/T(K)_b$
where $T(K)_b$ is the maximal bounded subgroup of $T(K)$
(this group is denoted by $\ti W$ in loc. cit.).  We have an exact sequence
\begin{equation}
1\to T(K)_b/T(K)_1\to \ti W\to W^\flat\to 1\ ,
\end{equation}
cf. [H-R], Rem. 10.
}
\end{Remarknumb}

\subsection{} Let us now discuss Schubert cells and varieties.

\begin{Definition}
Let $P$ be a parahoric subgroup of $G$ which corresponds
to a facet in the apartment of the maximal split torus $S$ 
and let $w\in \ti W^P\backslash \ti W/\ti W^P$.

 a) The {\sl Schubert cell} $C_w$ is the reduced subscheme
\begin{equation*}
L^+P\cdot n_w\subset \F_P=LG/L^+P,
\end{equation*}
 with $n_w\in N(K) $ a representative
of $w\in \ti W^P\backslash \ti W/\ti W^P$. \footnote{Note that $C_w$ is not a topological cell when $P$ is not an Iwahori subgroup.}

b) The {\sl Schubert variety} $S_w$ is the reduced scheme with underlying set the
Zariski closure $\overline{C_w}$ of the Schubert cell $C_w$ in $\F_P$.
It is a projective variety over $k$.
\end{Definition}

Our main result is the following:

\begin{thm}\label{normalmain}
Suppose that $G$ splits over a tamely ramified extension of $K$
and that the order of
the fundamental group of the derived group $\pi_1(G_{\rm der})$ is prime to the characteristic of $k$.
Then for each parahoric subgroup
$P$ of $G$ and $w\in   \ti W^P\backslash \ti W/\ti W^P$ the Schubert variety $S_w$
is normal, Frobenius split (when ${\rm char}(k)>0$), and has rational singularities.
\end{thm}

\begin{Remarknumb}
{\rm The normality of Schubert varieties has been proved in the
``Kac-Moody setting" by Kumar [Ku1],  by Mathieu [Ma1], and by Littelmann [Li]. That these Schubert varieties coincide in the case
of ${\rm SL}_n$
with the ones defined in our context follows in characteristic $0$ by an  integrability  result of Faltings. (See [B-L], App. ~to Sect.7.
The general split case is similar, see also \S \ref{KacM}). 
In  characteristic $p>0$, the comparison of the Schubert varieties in this article with the Kac-Moody Schubert varieties 
is harder and follows from the normality theorem (see Section \ref{CompAtp}). The normality of Schubert varieties 
in our context is proved in the case of ${\rm SL}_n$ in [P-R1]. Faltings [Fa1] proved the above result in the  case of an arbitrary  split
semisimple simply connected group. More details about this proof are given in [Go2] and [Fa2]. Our proof 
is an extension of Faltings' proof. 
}
\end{Remarknumb}

\begin{Remarknumb} \label{Borbit}
{\rm In fact, the conclusions of Theorem \ref{normalmain} also hold for 
the reduced Zariski closure of any (left) $L^+B$-orbit of the form $L^+B\cdot n_w$
in $\F_P$.
This follows immediately from the proof, see \ref{r1}.}
\end{Remarknumb}

\subsection{}\label{demazure}
 Before  dealing  with the proof we define the ``Demazure varieties" which are
of central importance in our arguments. Let us fix an alcove $C$
in the apartment corresponding to $S$ and consider the corresponding Iwahori subgroup $B=B_0$. We would like to fix ideas
and restrict our considerations to the Schubert
cells (varieties) for $P=B $ and $w\in W_a\subset \ti W$; by the above discussion and Theorem \ref{components},
the union of these is the (reduced) neutral component $\F^0_B=(LG)^0/L^+B $ of $\F_B$.
The advantage is that, at least set theoretically, these are the Schubert cells (varieties)
for the double Tits systems $(G(K)_1, B, N(K)_1, {\bf S})$.

\begin{prop}\label{QY}
a) The group scheme $L^+B$ over $k$ is a closed subgroup scheme of $L^+P_Y$ for each
$Y\subset {\bf S}$ as in \S \ref{affineWeylgrp}.

b)The fpqc quotient $L^+P_Y/L^+B$ is represented by a smooth projective
homogeneous space $Q_Y$ for the maximal reductive quotient $\overline P^{\rm red}_Y$
of the special fiber $\overline P_Y$;
the morphism $L^+P_{Y}\to L^+P_{Y}/L^+B\simeq Q_Y$ splits locally for the Zariski
topology on $Q_Y$. When $Y=\{i\}$ then $Q_Y\simeq {\bf P}^1_k$.
\end{prop}

\begin{Proof}
For every $k$-algebra $R$, we have $B(R[[t]])\subset P_Y(R[[t]])\subset G(R((t)))$.
Since both $B$, $P_Y$ are affine over $k[[t]]$,
 both $L^+B$, $L^+P_Y$ are closed $k$-subschemes of the ind-scheme $LG=LP_Y=LB$.
Part (a) follows.

To show part (b) we have to appeal to the construction of the parahoric group schemes $P=P_Y$
 in [B-TII], \S 4.6. (Note that since $k$ is algebraically closed,
$G$ is quasi-split and every unipotent group of finite type over $k$ is split;
in particular, ``quasi-reductive$=$reductive".) We consider the reduction homomorphism,
modulo $(t)$ which defines an exact sequence
\begin{equation}\label{reduct}
1\to \tilde U\to L^+P\xrightarrow{\ q\ }   \overline P\to 1.
\end{equation}
Here $\tilde U$ is a pro-unipotent pro-algebraic group over $k$.
Since $\overline P$ is affine and $\tilde U$ has a filtration by commutative
pro-unipotent group schemes
the $\ti U$-torsor $q$ is trivial. Let
$ \overline P^{\rm red}$ be the maximal reductive quotient
of $\overline P$. Then (\ref{reduct}) defines an exact sequence
\begin{equation}
1\to U\to L^+P\xrightarrow{\ q\ }   \overline P^{\rm red}\to 1
\end{equation}
where $U$ is again a pro-unipotent pro-algebraic group over $k$. (We can apply this to either
$P=P_Y$ or $P=P_{\emptyset}=B$.)
By loc.~cit.,~Thm.~4.6.33, applied to $P=P_{Y}$
we see that $L^+B\subset L^+P_Y$ is the inverse image
$q^{-1}(\overline B)$, where $\overline B$ is a  Borel subgroup of
$\overline P^{\rm red}_{Y}$. This implies that
$L^+P_Y/L^+B= \overline P^{\rm red}_Y/\overline B=Q_Y$
and that the quotient splits locally for the Zariski topology.
In fact, when $Y=\{i\}$, then the derived group of $\overline P^{\rm red}_{Y}$ is isomorphic to either
$SL_2$ or $PSL_2$  and $Q_Y\simeq {\bf P}^1_k$.  \endproof
\end{Proof}

\begin{prop}
Consider $w\in W_a$ and fix a reduced decomposition $w=s_{i_1}\cdots s_{i_r}$.
The Demazure variety $D(\ti w)$ is by definition the multiple ``contracted product"
\begin{equation*}
D(\ti w):=L^+P_{i_1} \times^{L^+B}\cdots \times^{L^+B} L^+P_{i_r}/L^+B\ .
\end{equation*}
(This depends on the choice of reduced decomposition; we write $\ti w$
instead of $w$ to indicate this dependence.) It is a smooth projective variety over $k$
of dimension equal to the length $l(w)=r$ of $w$ and affords a surjective morphism
\begin{equation}
D(\ti w)\xrightarrow{\pi_w} S_w\subset \F_B\ .
\end{equation}
\end{prop}

\begin{Proof}
By definition, the contracted product in the statement is the quotient
\begin{equation}\label{demaprod}
(L^+P_{i_1}\times\cdots \times L^+P_{i_r})/(L^+B)^r,
\end{equation}
where the product group scheme $(L^+B)^r$ acts on $L^+P_{i_1}\times\cdots \times L^+P_{i_r}$
from the right via
\begin{equation}\label{814}
(p_1,\ldots, p_r)\cdot (b_1,\ldots , b_r)=(p_1b_1, b_1^{-1}p_2b_2, \ldots , b_{r-1}^{-1}p_r b_r)\ .
\end{equation}
Recall that, by Proposition \ref{QY}, the quotient $L^+P_{i }/L^+B$
is isomorphic to the projective line ${\bf P}^1_k$. Let us write $w=w'\cdot s_{i_r}$. Note that
forgetting the last coordinate gives a morphism $D(\ti w)\to D(\ti w')$ which is a locally trivial fibration
with fiber $L^+P_{i_r}/L^+B\simeq {\bf P}^1_k$, comp.~[Go1], 3.3.1. The first statement of the Proposition
now follows by induction on the length of $w$. The morphism $\pi_w: D(\ti w)\to \F_B=LG/L^+B$
is given using the product in $LG$.
\endproof
\end{Proof}

\subsection{} \label{compMathieu}

Suppose here that $G$ is absolutely simple, simply connected 
and splits over a ramified extension $K'/K$ so that the setup 
and notations of \S \ref{tamepar}, \S \ref{KacM1} applies. 
For simplicity, set $\mathfrak g$ for the affine Kac-Moody algebra ${\mathfrak g}'_{KM}(G)_F$
that corresponds to $G$, where $F=K_0$ is the fraction field of $W$. For $w\in W_a$
(which then can be identified with the affine Weyl group in the Kac-Moody theory, [Kac], Ch. 3)
with reduced decomposition $\tilde w$, there is a Demazure variety $D_{\mathfrak g}(\ti w)$ over $k$
in the Kac-Moody setting (see [Ma1], p. 51). The variety $D_{\mathfrak g}(\ti w)$ is given by a contracted product similar to 
the one above but with the groups $L^+P_i$, $L^+B$, replaced by the ``parahoric" group schemes $\cal P_i$, $\cal B$
defined in [Ma1] Ch. I. The group schemes in [Ma1] are defined using the Chevalley 
integral form $U({\mathfrak g})^W$ of the universal enveloping algebra $U({\mathfrak g})$ of $\mathfrak g$;
by definition, this is the subalgebra of $U({\mathfrak g})$ which is generated over $W$ by the divided powers
$\displaystyle{\frac{e^{m}_i}{m!}}$, $\displaystyle{\frac{f^{m}_i}{m!}}$ and $\displaystyle{h_i\choose m}$
where $e_i, f_i, h_i$ are as in \S \ref{KacM1}. In particular, the group schemes $\cal P_i$, $\cal B$,
depend on this presentation of ${\mathfrak g}'_{KM}(G)_F$.
In what follows, we will see that the Demazure varieties $D(\ti w)$ above are naturally isomorphic
to their analogues $D_{\mathfrak g}(\ti w)$ in the Kac-Moody theory. 

It is enough to show this for the Iwahori subgroup scheme $B$ which 
is given as the $\sigma$-fixed points of the inverse image under reduction modulo $u$
of the Borel subgroup of $H_0$ that corresponds to $\Phi^+$.
Let also $P$ be a parahoric subgroup scheme
$P=P^G$ for $G$ which we can write as the $\sigma$-fixed point scheme of ${\rm Res}_{k[[u]]/k[[t]]}(P^H)$
as in  \S \ref{tamepar}; we can assume $B\subset P$. We have a semi-direct product decomposition which is stable 
under the $\sigma$-action
\begin{equation*}
L^+P^H=M^H\ltimes U^H\ .
\end{equation*}
Here $M^H  \simeq \overline{(P^H)}^{\rm red}$ is isomorphic to the maximal reductive quotient of the reduction
of $P^H$ modulo $u$ and $U^H$ is a pro-unipotent subgroup scheme over $k$. 
After taking $\sigma$-fixed points, this gives
\begin{equation}
L^+P=M\ltimes U(P)
\end{equation}
with $M\simeq \overline{P}^{\rm red}$ and $U$ a pro-unipotent subgroup scheme over $k$.
Denote by $U^H_m$
the ($\sigma$-stable) normal subgroup scheme of $U^H$ which is obtained as the kernel
of reduction modulo $u^m$
\begin{equation*}
P^H(m)(R)=\{ u\in P^H(R[[u]])\ |\ u\equiv 1 \in P^H(R[[u]]/(u^m))\}
\end{equation*}
We also set $P(m):=(U^H_m)^\sigma$ a normal pro-unipotent subgroup scheme of $U$. 
Given a reduced decomposition $ w=s_{i_1}\cdots s_{i_r}$, we can find integers $m_1,\ldots , m_r$
such that 
\begin{equation}
 P_{i_1}(m_1)\subset \cdots \subset  P_{i_r}(m_r)\subset L^+B\ .
\end{equation}
Consider the quotients $Q_{i_1}^{m_1}= L^+P_{i_1} / P_{i_1}(m_1), \ldots , Q_{i_r}^{m_r}=L^+P_{i_r}/ P_{i_r}(m_r)$, and
$L^+B^{m_1}=L^+B/ P_{i_1}(m_1),\ldots , L^+B^{m_r}=L^+B/ P_{i_r}(m_r)$; these are all smooth algebraic groups over $k$. 
It is not hard to see that the Demazure variety $D(\ti w)$ can be identified with the quotient 
\begin{equation} 
D(\ti w)=(Q_{i_1}^{m_1}\times\cdots \times Q_{i_r}^{m_r})/(L^+B^{m_1}\times \cdots \times L^+B^{m_r})
\end{equation}
where the  action is induced by (\ref{814}). A similar construction can be performed
in the Kac-Moody set-up. We obtain:
\begin{equation}
D_{\mathfrak g}(\ti w)=({\cal Q}_{i_1}^{m_1}\times\cdots \times {\cal Q}_{i_r}^{m_r})/({\cal B}^{m_1}\times \cdots \times {\cal B}^{m_r})
\end{equation}
with ${\cal Q}_{i}^{m}={\cal P}_i/ {\cal P}_i(m)$, similarly for ${\cal B}^{m_i}$. The algebraic groups
${\cal Q}_{i}^{m_1}$ and ${\cal B}^{m_i}$ both contain a central torus $\Gm$ that corresponds to the
central element $c$. 
To show our claim $D(\ti w)\simeq D_{\mathfrak g}(\ti w)$ is enough to show that there is an isomorphism $L^+P_{i}/ P_{i}(m)\simeq {\cal P}_{i}/(\Gm\times 
{\cal P}_i(m))$
which restricts to give
$L^+B/P_{i}(m)\simeq {\cal B}/(\Gm\times  {\cal P}_i(m))$. 
In fact, these algebraic groups lift to smooth group schemes over $W$ (using \S \ref{tamepar} and [Ma1]). We can easily see
from \ref{KacM1} and [Ma1] that such isomorphisms exist on the level of Lie algebras and hence between the generic fibers
over the fraction field $K_0$ of $W$. Denote by ${\mathfrak q}^{m}_i$ the Lie algebra of 
$L^+P_{i}/ P_{i}(m)_{K_0}= {\cal P}_{i}/(\Gm\times  {\cal P}_{i}(m))_{K_0}$ over $K_0$. 
Consider the $W$-module of distributions of these group schemes with support on the identity section (see [B-TII] 1.3).
By the above, these can be identified with $W$-lattices in the $K_0$-algebra $U({\mathfrak q}^m_i)$.
Using [B-TII] 3.5.1, we see that it is enough to show that the distribution lattices that correspond to the 
two group schemes coincide. Using the definition in [Ma1],
we can see that this amounts to checking that the image of $U({\mathfrak g})^W$ under
the natural map $U({\mathfrak g})\to U({\mathfrak q}^m_i)$
coincides with the lattice of distributions of $L^+P_{i}/ P_{i}(m))_{K_0}$.
This lattice can now be calculated in terms of the Chevalley basis $E_a$, $H_i$
(cf. [Kon], [Lu]). The expressions of $e_i$, $f_i$, $h_i$ in terms of  $E_a$, $H_i$
in \S \ref{KacM1} now allows us to show that the image of $U({\mathfrak g})^W$ above
agrees with this lattice.

We will return to a study of the Demazure varieties in the next section.

\subsection{}\label{reduceproof}
In the following few paragraphs we assume  that   the statement of Theorem \ref{normalmain}
is true under the following
%(additional)
 conditions:

\begin{itemize}
\item   $G$ is  absolutely simple,
simply connected and splits over a tamely ramified extension of $K$,

 \item  $P$ is an Iwahori subgroup.
\end{itemize}

We will show how to deduce the full statement of the Theorem
from this assumption.

\subsubsection{}\label{r1}
First of all observe that by Theorem \ref{represent} and Proposition \ref{QY} the natural morphism \begin{equation*}
\pi: \F_B=LG/L^+B\to \F_P=LG/L^+P_Y
\end{equation*}
is a proper and an \'etale locally trivial fibration
with   fibers  $L^+P_Y/L^+B\simeq Q_Y$. The inverse image under $\pi$ of
a Schubert variety  in $\F_P$ is a Schubert variety in $\F_B$. Hence, if all Schubert varieties in $\F_B$
are normal, Frobenius split and with rational singularities, the same is true for all Schubert varieties
in $\F_P$ (for the transfer of the $F$-splitting of $\F_B$ to the $F$-splitting of $\F_P$, one uses that the direct image
of ${\mathcal O}_{\F_B}$ is  ${\mathcal O}_{\F_P}$, cf. [Go1], Prop.~2.4.). From here  on we will assume that $P_Y=P=B$, i.e., is an Iwahori subgroup.

\subsubsection{}\label{r2} Suppose that $G$ is semi-simple, simply connected, splits over a tame extension and that
$P=B$. Then we can write
\begin{equation}
G=\prod_i{\rm Res}_{K_i/K}G_i,\qquad B=\prod_i{\rm Res}_{\O_{K_i}/\O_K}B_i,
\end{equation}
where $G_i$ are  absolutely simple and simply connected algebraic
groups over $K_i$ and $B_i$ an Iwahori subgroup scheme of $G_i$.
Each $G_i$ splits over a tamely ramified extension of $K_i$
(and $K_i/K$ itself is tamely ramified). We obtain
\begin{equation}
\F_B\simeq \prod_i L{\rm Res}_{K_i/K}(G_i)/L^+{\rm Res}_{\O_{K_i}/\O_K}(B_i)
\end{equation}
and the Schubert varieties in $\F_B$ are products of Schubert varieties
in the affine flag varieties $L{\rm Res}_{K_i/K}(G_i)/L^+{\rm Res}_{\O_{K_i}/\O_K}(B_i)$.
Note however that $K_i=k((u))\simeq K=k((t))$, $\O_{K_i}=k[[u]]\simeq \O_{K}=k[[t]]$ for each $i$.
These induce isomorphisms
$L{\rm Res}_{K_i/K}(G_i)\simeq LG_i$, $L^+{\rm Res}_{\O_{K_i}/\O_K}(B_i)\simeq L^+B_i$ as in (\ref{weilre2}). Therefore the
truth of the assertion of Theorem \ref{normalmain} for each $G_i $ implies the assertion for $G$. 

\subsubsection{}\label{r3} Suppose now that $G$ is semi-simple,
$p\nmid\# \pi_1(G)$,
and that  $G$ splits over a tamely ramified extension. Let $\ti G$ be the simply connected cover of $G$.
Then $\ti G$ satisfies the assumptions of the previous paragraph and by
the discussion above, Theorem \ref{normalmain} is true for all Schubert varieties
in $L\ti G/L^+\ti B$. By the work in \S \ref{Parcentral2}, these Schubert varieties
can be identified with the Schubert varieties $S_w$ in the neutral component
of $LG/L^+B$, i.e., the ones with $w\in W_a\subset \ti W$. Since $\ti W=W_a\rtimes \Omega$,
translating by elements of $\Omega$ allows us to conclude that Theorem \ref{normalmain} is true
for all Schubert varieties in $LG/L^+B$.

\subsubsection{}\label{r4} Finally suppose that $G$ satisfies the assumptions of Theorem \ref{normalmain}.
Then it follows from Proposition \ref{GG'prop} and  \S \ref{r3} (using a similar  argument
to extend ``outside the neutral component")
that Theorem \ref{normalmain} is true for all Schubert varieties in $LG/L^+B$.

\bigskip

\section{The proof of Theorem \ref{normalmain}} \label{proofmain}

\setcounter{equation}{0}

By \S \ref{reduceproof}, it is enough to deal with the case that $k$ is algebraically closed
and $G$
is absolutely simple, simply connected and splits
over a tamely ramified extension of $K=k((t))$. We can also assume
that $P=B$ is an Iwahori subgroup and for simplicity set $\F=\F_B=LG/L^+B$.
 Recall that then $G$ is quasi-split and the set-up of \S \ref{tamepar} applies.
We will continue with these assumptions for the rest of the section.
Recall that we then have $W_a=\ti W$, $G(K)_1=G(K)$.

\subsection{}\label{9a} Choose a maximal torus $T_0$ of the Chevalley group   $H=H_0\otimes_k K$ over $K$
which is defined over the field of constants $k$; we can assume that there is an action
of the diagram automorphism $\sigma_0$  on $T_0$ also defined over $k$. Then $T:=T_0\otimes_KK'$ is
a maximal torus of $G_{K'}=H_{K'}$ with an action of $\sigma=\sigma_0\otimes \tau$;
the fixed points  $S=T^{\sigma}$ form a maximal split torus of $G$ over $K$.
We will denote by $\Phi$ the set of (relative) roots of $G$ over $K$ with respect to $S$, and by $\Phi_{\aff}$ the
set of affine roots. Also let $ \Phi'$ be the (absolute) roots of $G_{K'}=H_{K'}$
with respect to $T$;  the Galois group
${\rm Gal}(K'/K)$ acts on $\Phi'$.   The set of orbits of $\Phi'$ under ${\rm Gal}(K'/K)$
can be identified with $\Phi$. Let us fix  a choice of Chevalley generators for the Lie algebra of $H_0$;
this provides us with Chevalley generators for the Lie algebra of $H=H_0\otimes_k K$ and therefore
with a Chevalley-Steinberg system  $(x_a)_{a\in \Phi}$  
([B-TII], 3.2, 4.1). We also fix a choice of a system of positive roots $\Phi^+$
and hence   a choice of an Iwahori subgroup $B$ of $G(K)$. It is enough to show the 
Theorem for such a $B$.

For every affine root $\alpha\in \Phi_{\aff}$, there is a corresponding
root subgroup scheme $U_\alpha\subset LG$ over $k$. In what follows, we will explain the construction
of these subgroups in some detail. We will refer to [Fa1] p. 46 for the construction of the affine root subgroups
$U'_{\alpha'}\subset LH_{K'}$ for the split group $H_{K'}$; this construction works over a general base and in particular 
over the ring of Witt vectors $W$.

We first have to recall the shape of
the root subgroup $U_a\subset G$ for $a\in \Phi$.
This is a subgroup scheme over $K$.  We distinguish the following  cases (see [B-TII], 4.1):

(I) $\frac{1}{2}a, 2a\not\in \Phi$. Then $x_a: U_a\simeq {\rm Res}_{L/K}{\bf A}^1_{L}$. In this case, 
$U_a\times_KK'$ is the direct product $\prod_{a'}U'_{a'}$ of  root
subgroups of $G_{K'}=H_{K'}$ for $a'$ ranging over the ${\rm Gal}(K'/K)$-orbit in $\Phi'$ that corresponds to
$a$ and $L$ is the fixed field of the stabilizer of $a'$, in particular $L=K$ or $K'$; we have
\begin{equation}
U_a=(\prod_{a'}U'_{a'})^\sigma\ .
\end{equation}
There is  
a group scheme homomorphism $\pi_a: G^a:={\rm Res}_{L/K}{ SL}_2\xrightarrow{} G$
such that $\pi_a$
induces an isomorphism between the maximal unipotent subgroups $U_+$ and $U_-$
of $G^a$ normalized by $\pi^{-1}_a(S)$ and
the subgroups $U_a$, $U_{-a}$ of $G$ respectively (see [B-TII], 4.1.4).
The homomorphism $\pi_a$ identifies the corrot lattice of $G^a$ with a direct summand
of the coroot lattice of $G$. Under our assumption that $G$ is simply connected 
this implies that $\pi_a$ is injective. 
The homomorphism $\pi_a$ of course depends 
on our choice of the Chevalley-Steinberg system.

(II) $a, 2a\in \Phi$. Then $K'/K$ is a separable (ramified) quadratic extension 
with Galois automorphism $\sigma(c)=\bar c$ and
$x_a: U_a\simeq  H(K',K)$ 
 where $H (K', K)$ is the group scheme over $K$ whose $R$-valued
points (for $R$ an $K$-algebra) are
\begin{equation*}
H (K',K)(R)=\{(c,d)\in (K'\otimes_KR)\times (K'\otimes_KR)\ |\   c\bar c=d+\bar d\}
\end{equation*}
with group law
\begin{equation}
(c,d)+(c', d')=(c+c', d+d'+\bar c c')\ .
\end{equation}
One can see that $U_{2a}$ is the subgroup scheme of $U_a$ which corresponds, under $x_a$,
to  $\{(0,d)\ |\ d+\bar d=0\}\subset H(K',K)(R)$. In this case, we can see that 
$U_a\times_KK'$ is the direct product 
$$
U'_{a'_1}U'_{a'_1+a'_2}U'_{a'_2}\subset G_{K'}=H_{K'}
$$
where $a'_1$, $a'_2$ are the two roots in $\Phi'$ 
in the Galois orbit corresponding to $a\in \Phi$. (Then $a'_1+a'_2$ is also in $\Phi'$.)
We have
\begin{equation}
U_a=(U'_{a'_1}U'_{a'_1+a'_2}U'_{a'_2})^\sigma, \quad U_{2a}=(U'_{a'_1+a'_2})^\sigma\ .
\end{equation}
As above we can see that there is an  injective
group scheme homomorphism $\pi_a: G^a:={SU}_3\xrightarrow{} G$
such that $\pi_a$
induces an isomorphism between the maximal unipotent subgroups $U_+$ and $U_-$
of $G^a$ normalized by $\pi^{-1}_a(S)$ and
the subgroups $U_a$, $U_{-a}$ of $G$ respectively. Here again $SU_3$ is the quasi-split special unitary group for $K'/K$, and the form with matrix equal 
to the antidiagonal unit matrix. 
In case (II) we can also set $G^{2a}:=G^{a }$, $\pi_{2a}:=\pi_{a}$. 
(Then $\pi_a$ and $G^a$ make sense for any root $a$, but they   really 
only depend on the ``root ray" given by $a$.)
Note that case II can only occur when ${\rm char}(k)\neq 2$.

\subsection{} \label{9b}
Recall that we denote by $t$, resp. $u$, uniformizers of
$K$, resp. $K'$, with $u^{[K':K]}=t$. 
Now consider the (hyperspecial) vertex $x_0$ in the Bruhat-Tits building ${\cal B}(H, K')$ that corresponds to 
the parahoric subgroup $H_0(k[[u]])$; since $\sigma$ acts on $H_0(k[[u]])$, the vertex $x_0$ is fixed by $\sigma$
and so it corresponds to a vertex $x_0$ in ${\cal B}(G, K)$. Affine roots $\alpha$, $\alpha'$ for $G$, resp. $H_{K'}$,  
can then be written in the form $\alpha(x)=\nu(\alpha)(x-x_0)+m$, $\alpha'(x)=\nu(\alpha')(x-x_0)+m$
with $m\in \Q$ and with vector parts $a=\nu(\alpha)\in \Phi$, $a'=\nu(\alpha')\in \Phi'$. Here, we are using the valuation 
of $K'$ with ${\rm val}(u)=1/e$. For simplicity, we write $\alpha=a+m$, $\alpha'=a'+m$.
The affine root subgroup $U'_{a'+m}$ in $H_{K'}=H_0\otimes_k k((u))$ corresponds to 
$\{r\cdot X_{a'}\cdot u^{m}|r\in k\}$ in the Lie algebra, where $X_{a'}$ is our choice (fixed above)
of a Chevalley generator 
for $H_0$ and the root $a'$.

If
$a=\nu(\alpha)$ falls in case I,  
we set
\begin{equation}\label{94}
U_\alpha=(\prod_{a'}U'_{a'+m})^\sigma\ ,
\end{equation}
where $U'_{a'+m}$ is the affine root subgroup
in $G_{K'}=H_{K'}$, and in the product $a'$ 
ranges over the ${\rm Gal}(K'/K)$-orbit in $\Phi'$ that corresponds to
$a$. In this case, we have $[L:K]\cdot m\in \Z$ and $U_\alpha$ is the $k$-subgroup scheme of $LU_a$ such that
\begin{equation*}
x_a(U_\alpha(R))=R\cdot t^{m}=R\cdot u^{em}.
\end{equation*}
If either $\frac{1}{2}\nu(\alpha)\in \Phi$, or $2\nu(\alpha)\in \Phi$, (case II), then there are two subcases:

1(II)  $ \nu(\alpha)=a$ with $2a\in \Phi$. Then $m\in \Z$ and 
we set
\begin{equation}\label{95}
U_\alpha=(U'_{a'_1+m}U'_{a'_1+a'_2+2m}U'_{a'_2+m})^\sigma\ .
\end{equation}
Then we have
\begin{equation*}
x_a(U_\alpha(R))= \left\{ ( r\cdot t^{m/2}, \frac{ (-1)^{m}}{2} r^2 t^{m}), \ \ r\in R\ \right\}.
\end{equation*}

2(II)  $\nu(\alpha)=2a$ with $a\in \Phi$. Then $2 m$ is an odd integer and 
we set
\begin{equation}\label{96}
\quad U_{\alpha}=(U'_{a'_1+a'_2+2m})^\sigma\ .
\end{equation}
Then we have
\begin{equation*}
x_{2a}(U_\alpha(R))=\{ ( 0, r\cdot t^m), \ \ r\in R\ \}.
\end{equation*}

(See also  [T1], Example 1.15, which explains the odd ramified unitary groups.
Observe
of course, that in the notation of Bruhat-Tits the affine root subgroups are given as groups of
field elements whose valuations are bounded from below, see Remark \ref{canonical}.
Some more details 
about the calculation of the affine root system $\Phi_{\rm aff}$
can be found in   [P-R3] \S 2.c .)

Note that for $(c,d)\in H(K',K)(R)$, we can write
\begin{equation}\label{sumRoot}
(c, d)=\left(c, \frac{1}{2}c\bar c\right)+\left(0, d-\frac{1}{2}c\bar c\right)
\end{equation}
(the sum in $H(K',K)(R)$).

For each affine root $\alpha=a+m$ there is a corresponding 
homomorphism 
\begin{equation}\label{rootsl2}
\phi_\alpha: SL_2\xrightarrow{} LG
\end{equation}
which has the property that it induces an isomorphism 
between
the unipotent upper triangular  matrices resp. the unipotent lower triangular matrices and 
$U_{\alpha_i}$, resp.~$U_{-\alpha_i}$. This homomorphism can be constructed as follows:

First assume that $a/2$ is not a root. Recall the homomorphism $\pi_a: G^a\to G$
over $K$ given in \S \ref{9a}. It is enough to give a homomorphism 
$SL_2\to LG^a$ corresponding to $\alpha=a+m$; then (\ref{rootsl2}) is obtained by composing with 
$\pi_a: LG^a\to LG$. If $2a\not\in \Phi$,  $G^a$ is either $SL_2$ or ${\rm Res}_{K'/K}SL_2$. Then 
we can reduce to the case $G=SL_2$. In this case, the homomorphism is given by
\begin{equation*}
\left( \begin{matrix}
 a&b\cr c& d 
 \end{matrix}\right)
 \mapsto 
\left( \begin{matrix}
 a& b\, t^{m}\cr c\, t^{-m}& d 
 \end{matrix}\right)\ .
\end{equation*}
If $2a$ is a root, then $G^a=SU_3$ and we can reduce to the case $G=SU_3$. 
Then the homomorphism is given by
\begin{equation*}
\left( \begin{matrix}
 a&b\cr c& d 
 \end{matrix}\right)
 \mapsto 
\left( \begin{matrix}
 a^2&   -u^{m}  ab & -(-1)^m t^{m} b^2/2   \cr
 -  u^{-m}2 ac &ad+bc &  (-1)^m u^m  bd \cr
 -(-1)^m t^{-m} 2c^2& (-1)^m u^{-m} 2 cd & d^2
 \end{matrix}\right)\ .
\end{equation*}

Now assume that $a/2=a'$ is a root. Then by considering $\pi_{a'}: G^{a'}\to G$
we can reduce to the case that $G=SU_3$. In this case, the homomorphism is given by
\begin{equation*}
\left( \begin{matrix}
 a&b\cr c& d 
 \end{matrix}\right)
 \mapsto 
\left( \begin{matrix}
 a& 0 & -b\, t^m\cr
 0&1&0\cr
 -c\, t^{-m}&0& d\cr
 \end{matrix}\right)\ .
\end{equation*}
(Note that in this case $m$ is a half integer and $t^m=u^{2m}$.)

For each simple affine root $\alpha_i=a_i+m_i$, the product
$U_{-\alpha_i}U_{\alpha_i}U_{-\alpha_i}$ contains
a representative of the reflection $s_i$. 
This can be seen as follows: By [T1] 1.4 we see that $s_i$ is the image of the unique element 
in the intersection $U_{-a_i}u\, U_{-a_i}\cap N$
for $u\in U_{\alpha_i}-\{1\}$ where $N$ is the normalizer of the maximal torus $T$
of $G$. Recall the homomorphism $\phi_{\alpha_i}: SL_2\to LG$
constructed above. 
Since 
\begin{equation}\label{uns}
U_{-}\left(\begin{matrix} 1&1\cr 0& 1\end{matrix}\right) U_{-}\cap N_{SL_2}=
\left\{\left(\begin{matrix} 0&1\cr -1& 0\end{matrix}\right)\right\},
\end{equation}
we can produce the element $s_i$ 
by taking the image under $\phi_{\alpha_i}: SL_2\to LG$ of 
the matrix on the right hand side of the equation (\ref{uns}) above.

For $w=s_{i_1}\cdots s_{i_r}\in W_a$, $\alpha\in \Phi_{\aff}$, we see that we have
\begin{equation}\label{conj}
w^{-1}\cdot U_{\alpha}\cdot w=U_{\alpha\cdot w}
\end{equation}
where we view $w$ as acting via an affine transformation on the apartment of $S$.

\begin{Remark}\label{canonical}
{\rm The group schemes $U_\alpha\subset LG$ can be understood
in the context of Bruhat-Tits theory as follows: Our choice of Chevalley basis
gives a root datum valuation $(\phi_a)_{a\in\Phi}$. Given an affine root $\alpha=a+m$ as above,
Bruhat and Tits [B-TII] construct a smooth affine group scheme
${\mathfrak U}_\alpha$ over $k[[t]]$. This has generic fiber the root subgroup $U_a$ 
and is such that ${\mathfrak U}_{\alpha}(k[[t]])=\{u\in U_a(k((t)))\ |\  u=1,\ \hbox{\rm or\ }
\phi_a(u)\geq m\}$.  Similarly, we have the subgroup scheme ${\mathfrak U}_{\alpha+}$ with 
${\mathfrak U}_{\alpha+}(k[[t]])=\{u\in U_a(k((t)))\ |\  u=1,\ \hbox{\rm or\ }
\phi_a(u)>m\}$.
We can see that $U_\alpha\subset L^+{\mathfrak U}_\alpha$ becomes isomorphic
to the
 quotient
$\overline {\mathfrak U}_{\alpha}:=L^+{\mathfrak U}_\alpha/L^+{\mathfrak U}_{\alpha+}$
with the obvious map.
}
\end{Remark}

\begin{Remark}\label{overW}
{\rm Our construction also works to 
define affine root subgroups $\und{U}_{\alpha}$, $\alpha\in \Phi_{\rm aff}$ over the ring 
of Witt vectors. These are obtained via taking fixed points of   products (as in (\ref{94}), (\ref{95}), 
(\ref{96})) of the affine root subgroups of [Fa1], p. 46, 
for the (non-twisted) loop group $L\und{H}$ over $W$. Using \S \ref{tamepar}
we see that $\und{U}_{\alpha}$ are subgroup schemes of $L\und G$ that satisfy (\ref{conj})
with special fiber the group schemes $U_{\alpha}$ that we have constructed above.
We also see that the group schemes $G^a$ and the homomorphisms $\pi_a: G^a\to G$ lift to group
schemes $\und{G}^a$ and homomorphisms $\und{\pi}_a: \und{G}^a\to \und{G}$ over $W$. This can be seen as follows.
If $G$ is split, $G^a=SL_2$ and the statement follows from the construction of the Chevalley group
scheme $\und G$ over $W[[t]]$. In general, there are two cases as in \S \ref{9a}: Recall the group scheme $\und{H}$ 
over $W[[u]]$ of \S \ref{tamepar} which is such that $\und{G}=({\rm Res}_{W[[u]]/W[[t]]}\und {H})^\sigma$. In case (I), we can construct
$\und{G}^a$ as the $\sigma$-fixed point scheme of $\prod_{\{a'\}}\und {H}^{a'}$.
In case (II), there is group scheme embedding $\und{SL}_3\to \und{H}$ over $W[[u]]$ corresponding 
to the root subset $\{a'_1, a'_2, a'_1+a'_2\}$. Then $\und{G}^a=({\rm Res}_{W[[u]]/W[[t]]}\und{SL}_3)^{\sigma}$.}
\end{Remark}

\subsection{}

We now state a key ingredient of the proof:

\begin{prop}\label{key}{(\bf Key observation)}
The Lie algebra $Lie(LG)$ of $LG$ is spanned over $k$  by the Lie  algebras
of $L^+B$ and of $U_\alpha$ for all $\alpha\in \Phi_{\rm aff}$.
\end{prop}
\medskip

\begin{Proof}
The reductive group $G$ over $K=k((t))$ contains  the product
\begin{equation}
\prod_{a\in \Phi^+}U_{-a}\times T\times \prod_{a\in \Phi^+}U_a
\end{equation}
as a Zariski open neighborhood of the origin. Hence, it is enough to
show that the Lie algebras $Lie(U_a) $, $a\in \Phi$, and $Lie(T)$
are spanned over $k$ by elements from the Lie  algebras
of $L^+B$ and $U_\alpha$ for  $\alpha\in \Phi_{\rm aff}$.

Let us first consider the unipotent subgroups.
For each $a\in \Phi$, consider   $u\in LU_a(R)=U_a(R((t)))$.
Note that elements  of $ LU_a(R)$ of sufficiently high valuation
lie in $L^+B(R)=B(R[[t]])$. Hence, we can see using
(\ref{sumRoot}) that we can write $u=u'+u_B$ with
$u_B\in B(R[[t]])$ and $u'$ a finite sum of elements in
affine root subgroups $U_{\alpha}(R)$ with $\nu(\alpha)=a$.
In particular, the result on Lie algebras follows.

Now let us consider the maximal torus $T$. We would like to show that
each element of $T(R((t)))$ for $R$ a local Artinian $k$-algebra with $R/M=k$,
$M^2=(0)$, which reduces to the identity  modulo $M$, can be written in $G(R((t)))$ as a product of   elements  
in $U_a(R((t)))$ for $a\in \Phi$. Then the result will follow by the above.
Recall that $G$ is quasi-split, absolutely simple, semi-simple and simply connected, and splits over the  
field $K'/K$, which is tame. The cocharacter group $X_*(T)$ coincides with the coroot lattice
of the  split group $G_{K'}$.
Let $\Delta'$ be a basis for the roots $\Phi'$ of  
$G_{K'}$ over $K'$ that corresponds to a Borel
subgroup defined over $K$. 
The torus $T$ splits
into a product of induced tori $T_b$ over the ${\rm Gal}(K'/K)$-orbits $b$ on
 $\Delta'$. This set of orbits  can be identified with a basis
of  roots of $\Phi$. For a root $a$ in this basis, recall the homomorphism $\pi_a: G^a\to G$
given in \S \ref{9a}. The torus $T_b$ is equal to $\pi_a(T^a)$ with $T^a=\pi^{-1}_a(T)$
a maximal torus of $G^a$. Since $\pi_a$ identifies the unipotent subgroups $U_+$ and $U_-$
of $G^a$ with $U_a$ and $U_{-a}$, these considerations show that it is enough to consider   the cases
(corresponding to (I) and (II) of \S \ref{9a}):

 (I) $G=SL_2$ (split), and (II) $G=SU_3$ (for the ramified quadratic extension
 $K'/K$).

 In case (I) our claim follows from the identity, cf. [Fa1], p. 53:
\begin{equation}
\left(\begin{matrix}
c& 0\\ 0& c^{-1}\\
\end{matrix}\right)=
\left(\begin{matrix}
1& c\\ 0& 1\\
\end{matrix}\right)
\left(\begin{matrix}
1& 0\\ -c^{-1}& 1\\
\end{matrix}\right)
\left(\begin{matrix}
1& c\\ 0& 1\\
\end{matrix}\right)
\left(\begin{matrix}
0& -1\\ 1& 0\\
\end{matrix}\right).
\end{equation}
This indeed shows that these elements are in the subgroup generated 
by the root groups.

 In case (II) let
 $$
 t(d)={\rm diag}\left(-d, \frac{\bar d}{d}, -\frac{1}{\bar d}\right),
 $$
  with $d\in (R((u)))^\times$,  be an element of
 the maximal torus $T\subset SU_3(R((t)))\subset SL_3(R((u)))$. If $d+\bar d=c\bar c$
 with $c\in R((u))$ we can use the identity, cf. [St], p.~537:
 \begin{equation}\label{identitySu3}
t(d)={\rm diag}\left(-d, \frac{\bar d}{d}, -\frac{1}{\bar d}\right)= x(c,d)\cdot y\left(-\frac{\bar c}{\bar d}, \frac{1}{\bar d}\right)\cdot x\left(c\frac{\bar d}{d}, d\right)\cdot w
 \end{equation}
where
\begin{equation}
x(c, d)=\left(\begin{matrix}1 &-c& -d\\ 0&1& \bar c\\ 0&0& 1\\
\end{matrix}\right),
\quad y(c, d)=
\left(\begin{matrix} 1 &0& 0\\ -c
&1&0 \\ -d&\bar c& 1\\
\end{matrix}\right), \quad
w=\left(\begin{matrix}0 &0& 1\\ 0&-1& 0\\ 1&0& 0\\
\end{matrix}\right).
\end{equation}
 Here $x(c, d)\in U_+$, $y(c,d)\in U_-$ for the special unitary group $SU_3$.
Note that $U_+\simeq U_-\simeq H(K', K)$. In general, given an element $t(d)\in T(R((t)))\simeq R((u))^\times$, which reduces to the identity modulo $M$, 
we can find $c$ in $R((t))$ with $c^2=c\bar c=d+\bar d$ (since char $k\neq 2$). The identity
(\ref{identitySu3}) now allows us to write $t(d)$ as a product of unipotent matrices
and the result follows.\endproof
\end{Proof}
\medskip

In fact, the following extension
of the key Proposition \ref{key} can be shown to hold by using Remark \ref{overW} and repeating the proof verbatim
(we can redo the argument by replacing  the local Artin $k$-algebra $R$ with a local  Artin $W$-algebra 
with residue field $k$).
 
\begin{prop}\label{prop911}
The Lie algebra of $ L\und G $ is generated as a $W$-module by the Lie algebras
of $L^+ \und B$ and $\und U_{\alpha}$, $\alpha\in \Phi_{\aff}$.\endproof
\end{prop}

\begin{cor}\label{generate}
 If $R$ is an Artinian local $k$-algebra, then $LG(R)$ is generated by
the subgroups $L^+B(R)$
and $U_{\alpha}(R)$ for $\alpha\in \Phi_{\rm aff}$.
\end{cor}

\begin{Proof}
When $R $ is a field, our claim follows from the Bruhat decomposition
and the fact that $U_{-\alpha_i}U_{\alpha_i}U_{-\alpha_i}$ contains
the reflection $s_i$ for every simple affine root $\alpha_i=a_i+m_i$. 

Now for a general Artinian $R$, choose  a minimal ideal $I\subset R$ such that the
assertion is true for $R/I$. Consider $g\in LG(R/I^2)$ and its reduction $\bar g\in LG(R/I)$.
By assumption we can write $\bar g$ as a product of
points in  $L^+B(R/I)$, $U_{\alpha}(R/I)$. By smoothness, these
points   can be lifted to corresponding points with values in $R/I^2$
whose product is an element $g'\in LG(R/I^2)$ such that $h=g'\cdot g^{-1}\equiv \ 1\ {\rm mod}\ I/I^2$.
Since $h$ is given by an element of $Lie(LG)\otimes _k I/I^2$ our claim follows
from the key Proposition \ref{key} above.\endproof
\end{Proof}

\subsection{}\label{demazureplus}
We now give some additional properties of the Demazure varieties $D(\ti w)$;
we continue with the notations of \S \ref{demazure}. The proofs of the following two 
propositions are similar to proofs of corresponding statements in [Fa1] (see also [Go2]). We will
sketch the arguments and refer the reader to [Fa1] and [Go1-2] for more details.

\begin{prop}\label{d1}
a) For any reduced  $\ti u\leq \ti w$  (i.e., a reduced decomposition
obtained by omitting reflections from $w=s_{i_1}\cdots s_{i_r}$),
there is a closed
immersion $\sigma_{\ti u, \ti w} : D(\ti u)\to D(\ti w)$ and
a  commutative
diagram
\begin{equation}
\begin{CD}
D(\ti u)@>\sigma_{\ti u, \ti w} >>D(\ti w)\\
@V\pi_uVV @VV\pi_wV\\
S_u@>>>S_w\\
\end{CD}
\end{equation}
in which the bottom horizontal arrow is the natural closed immersion.
For $\ti v\leq \ti u\leq \ti w$, we have $\sigma_{\ti v, \ti w} =\sigma_{\ti u, \ti w} \cdot \sigma_{\ti v, \ti u} $.

b) For every $\ti w$, the morphism  $\pi_w: D(\ti w)\to S_w$ is proper and birational
and has geometrically connected fibers.

c) When ${\rm char}(k)>0$, the varieties $D(\ti w)$ are Frobenius split compatibly with the closed
immersions $\sigma_{\ti u, \ti w} : D(\ti u)\to D(\ti w)$, for all reduced $\ti u\leq \ti w$.
\end{prop}

\begin{Proof}
Here the assertion (a) is clear. To see  the birationality of the morphism in (b), note that the open 
subset of $D(\tilde w)$ ``where no factor in (\ref{demaprod}) lies in $L^+B$" maps isomorphically onto the  Schubert
cell $C_w$. This open subset is equal to $U_{\alpha_{i_{1}}}\cdot s_{i_{1}}\cdot U_{\alpha_{i_{2}}}\cdot s_{i_2} \dots \cdot U_{\alpha_{i_l}}\cdot s_{i_r}$. 
The fact that the fibers are geometrically connected is shown by induction on $r=l(w)$.
Indeed, if we write $\ti w=s_i\ti u$ with
$l(u)=l(w)-1$ then we can factor $\ti \pi_w$ as   
\begin{equation*}
D(\ti w)=L^+P_i\times^{L^+ B}D(\ti u) \rightarrow L^+ P_i\times^{L^+ B}  S_u \rightarrow   S_w .
\end{equation*}
The Frobenius splitting of $D(\tilde w)$ is constructed by using the criterion of Mehta and Ramanathan, comp.~[Go1], Prop.~2.5. More precisely,
one calculates the canonical bundle $\omega_{D(\tilde w)}$ of $D(\tilde w)$ as in [Go1], Prop.~3.19, and shows that there
exists a global section $s$ of $\omega_{D(\tilde w)}^{-1}$ with divisor equal to the sum of the boundary divisor
of $D(\tilde w)$ and of an effective divisor which does not contain the origin of $D(\tilde w)$. 
(The corresponding calculation in the general symmetrizable Kac-Moody case is also  contained in [Ma1], Ch. 8, 18.) As in [Go1],
Cor.~3.23 it is checked that the resulting Frobenius splitting of $D(\tilde w)$ is compatible
with the closed immersions $\sigma_{\ti u, \ti w}$, for all reduced $\ti u\leq \ti w$.
\endproof
\end{Proof}

\smallskip
 Let $\psi_w: \ti S_w\to S_w$ be the normalization of $S_w$. We obtain a factorization of $\pi_w$
\begin{equation}\label{factorise}
D(\ti w)\xrightarrow{\ti \pi_w}\ti S_w\xrightarrow{\psi_w} S_w\ .
\end{equation}

\begin{prop}\label{d2}
a) The morphism $\psi_w$ is a universal homeomorphism.

b) For any  $ u\leq  w$, there is a closed immersion $\ti S_u\to \ti S_w$
which lifts the natural closed immersion $S_u\to S_w$. If $\ti u\leq \ti w$ are reduced
decompositions, then the diagram
\begin{equation}\label{imm}
\begin{CD}
D(\ti u)@>\sigma_{\ti u, \ti w} >>D(\ti w)\\
@V\ti\pi_uVV @VV\ti\pi_wV\\
\ti S_u@>>>\ti S_w\\
\end{CD}
\end{equation}
commutes.

c) When ${\rm char}(k)>0$, the varieties $\ti S_w$ are Frobenius split, compatibly with the closed immersions
$\ti S_u\to \ti S_w$, for all $u\leq w$.

d) For $i>0$, we have $R^i(\ti \pi_w)_*(\O_{D(\ti w)})=(0)$ and $(\ti \pi_w)_*(\O_{D(\ti w)})=\O_{\ti S_w}$.
The variety $\ti S_w$ is Cohen-Macaulay and it has at most rational singularities.
\end{prop}

\begin{Proof}
To prove (a), we note that by Proposition \ref{d1}, the fibers of $\pi_w$ are geometrically connected, and hence $\psi_w$ is a universal homeomorphism. The diagram in (b) arises since (\ref{factorise}) is the Stein factorization of $\pi_w$. Since 
$\psi_w$ and $\psi_u$ are universal homeomorphisms, the morphism $\tilde S_u\to \tilde S_w$ is a universal homeomorphism onto its image $\tilde T_u$ in $\tilde S_w$. The compatible $F$-splitting of $D(\tilde u)$ in $D(\tilde w)$ induces an $F$-splitting of $\tilde T_u$, which is compatible with the $F$-splitting of $\tilde S_u$ induced from the $F$-splitting of $D(\tilde u)$, i.e., the following diagram is commutative,
$$
\begin{CD}
F_*(\mathcal O_{\ti T_u})@>>>   F_*(\mathcal O_{\ti S_u}) \\
@VVV @VVV\\
\mathcal O_{\ti T_u}@>>>\mathcal O_{\ti S_u} 
\end{CD}
$$
in which the horizontal arrows are injections. We can see that this implies that $\mathcal O_{\tilde T_u}=\mathcal O_{\tilde S_u}$, and hence the morphism $\tilde S_u \to \tilde S_w$ is a closed immersion, and the variety $\tilde S_w$ is $F$-split compatibly with the closed subscheme $\tilde S_u$. 

To show (d), it remains to show that
\begin{equation*}
R^i(\ti \pi_w)_*(\O_{D(\ti w)})=(0) \ \  \text{and}\ \  R^i(\ti \pi_w)_*(\omega_{D(\ti w)})=(0) \ \ \text{for}\ \  i>0,
\end{equation*}
and  $(\ti \pi_w)_*(\O_{D(\ti w)})=\O_{\ti S_w}$. Here the assertion concerning $\omega_{D(\ti w)}$ follows from the rest by the Grauert-Riemenschneider theorem
(for Frobenius split varieties in positive characteristic, cf.~[M-K]). The assertions concerning
 $\mathcal O_{D(\ti w)}$ 
are proved by induction on $l(w)$, cf.~[Go2], Lemma 3.13. Let $\ti w=s_i\ti u$ with
$l(u)=l(w)-1$. Let us factor $\ti \pi_w$ as   
\begin{equation*}
D(\ti w)=L^+P_i\times^{L^+ B}D(\ti u) \rightarrow L^+ P_i\times^{L^+ B}\tilde S_u \rightarrow \ti S_w .
\end{equation*}
By induction we may assume that the higher direct images of $\mathcal O_{D(\ti w)}$ under the first morphism vanish.
As in [Fa1], p. 51, we can see
that the second morphism has as geometric fibers either a point or a $\bf P^1$.

By the lemma of Mehta-Srinivas ([Go2], Lemma 3.14), the vanishing of the higher cohomology groups of each fiber implies now, in the presence of a Frobenius splitting, the vanishing of the higher cohomology sheaves. 
This proves (d) when the characteristic of $k$ is positive; the characteristic 
zero case follows using a semicontinuity argument, using the lifting over $W(k)$, c.f. \ref{demazureZ}.
This finishes the sketch of the proof of the proposition.
\endproof
\end{Proof}

\smallskip

Now consider the ind-scheme $S$ over $k$ defined as the inductive limit
of the directed set of the Schubert varieties
$S_w$, $w\in W_a$, taken with the Bruhat order:
$S:=\displaystyle{\lim_{\to}   S_w  }$;  this
is an ind-subscheme of the affine flag variety $\F $.
In fact, since $G$ is simple and simply connected, it coincides with the underlying reduced sub ind-scheme
$(\F)_{\rm red}$ of
$\F $. Using Proposition \ref{d2} we can also define
an ind-scheme $\ti S$ over $k$ as the inductive limit
of the directed set $\ti S_w$, $w\in W_a$, taken with the Bruhat order,
\begin{equation*}
\ti S:=\lim_{\to}  (\ti S_w )\ .
\end{equation*}
We have
\begin{equation*}
\ti S\xrightarrow{\psi} S=(\F)_{\rm red}\hookrightarrow \F
\end{equation*}
and Theorem \ref{normalmain} will follow if we show that, under our conditions,
both maps are isomorphisms.

Now for $\alpha_i$, $i\in I$, a simple affine root, the subgroup scheme $L^+P_i$ contains
the corresponding reflection $s_i$ and also the root subgroups $U_{\alpha_i}$, $U_{-\alpha_i}$. The group scheme $L^+P_i$ acts compatibly on $D(\ti w)$ and $S_w$, provided that $l(s_i w)<l(w)$. Hence  $L^+P_i$ acts also on the normalization $\ti S_w$ when
$l(s_iw)<l(w)$. The proof of the following is straightforward.

\begin{lemma}
For each $w\in W_a$ and each $i\in I$, there is $w'\in W_a$ with $w\leq w'$ and $l(s_iw')<l(w')$.
\endproof
\end{lemma}

The lemma, together with the above discussion and the definition of $\ti S_w$
using the Demazure varieties implies that for each $i\in I$, $L^+P_i$ acts on $\ti S$.
Since $W_a$ is generated by the reflections $s_i\in L^+P_i$, $i\in I$, the above discussion then implies
that the ind-scheme $\ti S$ supports an action of $L^+B$ and of $U_\alpha$ for all affine roots $\alpha\in \Phi_{\aff}$.
The same is of course also true for the ind-scheme $S$; the morphism $\psi: \ti S\to S$ is equivariant
with respect to these actions.

\begin{prop}\label{limitSchu}
Under our assumptions (i.e the group $G$ is absolutely simple, simply connected and splits over a tamely ramified extension of $K$), we have
\begin{equation*}
S=\lim_{\to } S_w=\F\ .
\end{equation*}
\end{prop}

\begin{Proof}
Let us write the ind-scheme $\F$ as a union $\F=\cup_{j}\F_j$
of proper $k$-schemes (Proposition \ref{fixed} (b)). We have
\begin{equation*}
\cup_{w\in W_a}(\F_j\cap S_w)\subset \F_j
\end{equation*}
and given $j$, there is $w\in W_a$ such that $(\F_j)_{\rm red}\subset S_w$.
Let $x$ be a generic point of the support of the ideal sheaf 
of $\O_{\F_j}$ which defines the scheme theoretic intersection
$\F_j\cap S$.  Consider the ideal of definition $I_w$
of $S_w\cap \F_j$ in the local ring $A= \O_{\F_j, x}$ with maximal ideal $M$. Observe that $I_{w'}\subset I_{w}$ if $w'\geq w$
and $A$ is Noetherian. Then by the above, there is a $w\in W_a$ such that
$A /I_w$ is of finite length and is such that if $w'\geq w$ then $I_w=I_{w'}$. We conclude that
there is $n>0$ such that $I_w$,  for all $w\in W_a$,  contains $M^n$ properly.
This gives a point $x_n\in \F_j(A/M^n)$ such that $x_n\not\in S(A/M^n)$. We can write
$x_n=g_n\ {\rm mod}\ L^+B(R)$, where $R=A/M^n$ and $g_n\in LG(R)$ for the Artin
local ring $R=A/M^n$. Proposition \ref{generate} and the fact
that $L^+B(R)$ and $U_{\alpha}(R)$ act on $S$ (see above) shows that
$x_n$ is in $S(R)$ which is a contradiction.\endproof
\end{Proof}

\smallskip

\begin{cor}
If $G$ is semi-simple, simply connected and splits over a tamely ramified extension of $K$,
then $\F=LG/L^+B$ is a reduced ind-scheme. 
\end{cor}

\begin{Proof}
If $G$ is in addition absolutely simple, this 
follows from  Proposition \ref{limitSchu}.
We can reduce to this case  using an argument as in \S \ref{r2}. 
\endproof
\end{Proof}
\smallskip

In view of our reductions in \S 6, this also completes the proof of Theorem \ref{redu}. 

\subsection{}\label{demazureZ}
Suppose now that ${\rm char}(k)=p>0$.
Recall that under our assumptions, we have using \S \ref{tamepar}, a lift of
$LG$, $L^+P$, $\F=LG/L^+B$ over $k$ to $L\und G$, $L^+\und P$, $\und \F=L\und G/L^+\und B$
over the ring of Witt vectors $W=W(k)$. Also, by
\S \ref{9a}, we have the affine root subgroups $U_\alpha\subset LG$ and their lifts 
$\und U_\alpha\subset L\und G$ over $W$ (see Remark \ref{overW}).

Denote by $K_0$ the fraction field of $W$.
For $w\in W_a$, we can define $\und S_w$,
the Schubert variety over $W$, as the reduced subscheme over $W$ that
is the Zariski closure of the ``cell" $L^+\und B\cdot n_w\subset L\und G/L^+\und B$.
For a reduced expression $\ti w=s_{i_1}\cdots s_{i_r}$,
we can again define the Demazure variety over $W$ as the quotient
\begin{equation}
 \und D(\ti w)=(L^+\und P_{i_1}\times\cdots \times L^+\und P_{i_r})/(L^+\und B)^r.
\end{equation}
This is smooth and projective over $W$, is an iterated ${\bf P}^1_W$-fibration,
 and supports
a proper morphism induced by multiplication
\begin{equation}
\und D(\ti w)\xrightarrow{\und \pi_w} \und S_w.
\end{equation}
Now set $\tiuS_w=\Spec((\ti \pi_w)_*(\O_{\und D(\ti w)}))$, i.e., $\tiuS_w$ is the normalization
of $\und S_w$ in $\und D(\ti w)$; then the structure morphism $\tiuS_w\to \Spec(W)$ is flat.
The following can be shown following the arguments of [Fa1] (see [Go2], Prop.~3.15 et seq.~ for some more details):

\begin{prop}\label{liftti}
a) The formation of $\tiuS_w$ commutes with  base change on $W$. In particular,
there is a  natural isomorphism
\begin{equation}
\tiuS_w\otimes_Wk\simeq \ti S_w\ .
\end{equation}

b) If $u\leq w$, then there is a closed immersion $\tiuS_u\to \tiuS_w$
which lifts the natural closed immersion $\und S_u\hookrightarrow \und S_w$
such that a diagram as (\ref{imm}) commutes.\endproof
\end{prop}

Using Proposition \ref{liftti} (b), we can define an ind-scheme $\tiuS$ over $W$ as the inductive limit
of the directed set $\tiuS_w$, $w\in W_a$ taken with the Bruhat order:
\begin{equation*}
\tiuS:=\lim_{\to}  (\tiuS_w )\ .
\end{equation*}
In fact, because of Proposition \ref{liftti} (a), we have
\begin{equation}
\tiuS \otimes_Wk\simeq \ti S \ .
\end{equation}

We can see that the group schemes $L^+\und B$ and $\und U_\alpha$, $\alpha\in \Phi_{\aff}$,
act on $\und S=\displaystyle{\lim_{\to}\und S_w}$ and $\tiuS$.

\bigskip

\subsection{}\label{KacM}
Recall that $K_0$ is the fraction field of $W(k)$. Here we show that, for any $w$, $\ti S_{w}\otimes_WK_0=S_{w}\otimes_WK_0$,
i.e., $S_{w,K_0}:=S_w\otimes_{W}K_0$ is normal. It is enough to show this statement with $K_0$ replaced
by its algebraic closure $\bar K_0$. Observe that the normality of Schubert varieties in affine flag 
varieties associated to (affine) Kac-Moody Lie algebras is known by the work
of Kumar, Mathieu, Littelmann ([Ku1], [Ma1], [Li]). Our strategy is essentially to show that  the (characteristic zero)
Schubert varieties $S_{w, \bar K_0}$ coincide with corresponding Schubert varieties
in the theory of (affine) Kac-Moody Lie algebras.

Recall that  by \S \ref{tamepar} we are dealing with an absolutely simple, simply connected
(quasi-split) group over $\bar K_0((t))$ which is one of the types described in the beginning of that paragraph.
All these types are realized over $\C((t))$  and by the Lefschetz principle it is enough to show
that the Schubert varieties $S_{w, \C}$   are normal  for  groups $G$ over $\C((t))$  as in \S \ref{tamepar}.
For simplicity, in the rest of the paragraph, we will omit the subscript $\C$. Recall our notation 
$\gfr'_{KM}(G)$ (and $\gfr_{KM}(G)$) for the affine Kac-Moody Lie algebra over $\C$ associated to 
the local Dynkin diagram of $G$ (see \ref{KacM1}).

Now consider the embedding $\F_G=LG/L^+B\hookrightarrow \F_H=LH/L^+B_H$ described in \S \ref{tamepar}.
The affine flag variety $\F_H$ for the split group $H=H_0\otimes_\C\C((t))$ is the affine flag variety
for $H_0$ considered e.g in [Fa1]. We can give a projective embedding $\F_H\hookrightarrow {\bf P}(V)$
by suitably embedding $LH$ into $LSL_n$ and composing the resulting closed immersion $\F_H\hookrightarrow \F_{SL_n}$
with the standard projective emdedding of $\F_{SL_n}\hookrightarrow {\bf P}(V)$ given by its lattice description.
(Here $V$ is a countably infinite dimensional complex vector space and ${\bf P}(V)$
stands for the corresponding projective space which is naturally an ind-scheme over $\C$.)
We obtain a projective embedding $\F_G\hookrightarrow {\bf P}(V)$; this induces a line bundle $\L$ on $\F_G$;
which we can pull back to $ S_w$, $\ti S_w$ and $D(\ti w)$ for each $w\in W_a$. 
By Proposition \ref{d2} we see that, when $w'\leq w$, there are pull back (restriction) homomorphisms
\begin{equation}\label{restrictionL}
\Gamma(\ti S_w, \L)\xrightarrow{} \Gamma(\ti S_{w'}, \L);
\end{equation}
 By the Frobenius splitting of the normalizations  in
characteristic $p$ (i.e., over $k$) we can see (cf. [Ma1], Lemma 137)
that  the analogous restriction homomorphisms (\ref{restrictionL}) in characteristic $p$
are surjective.  By semicontinuity the homomorphisms (\ref{restrictionL}) themselves are then also surjective. Now after passing to the $\C$-duals, we can consider the inductive limit of injective maps
\begin{equation}\label{dirlim}
\ti E(\L):=\lim_{\to}(\Gamma(\ti S_{w}, \L)^*)\ .
\end{equation}
The rest of the argument is as in [Fa1] p. 52 (which also follows [Ku1], see also [Ku2]). 
This complex vector space $\ti E(\L)$ admits actions by $L^+P_i$. For each index $i$ there is an $sl_2$-triple $\{e_i, f_i, h_i\}$
in the Lie algebra of $LG$; this is given via the homomorphism $\phi_i: SL_2\to LG$ of \S \ref{9b}. The Lie algebra of $L^+P_i$
contains all the elements $e_j$ and $h_j$ and the element $f_i$ and we can see that these satisfy the additional relations
$[e_j, f_i]=0$ for $i\neq j$, $[h_r, e_s]=a_{rs}e_s$, and $[h_r, f_i]=-a_{ri}f_i$, where 
$\{a_{rs}\}$ is the generalized Cartan matrix. The action of the Lie algebra of $L^+P_i$ on 
$\ti E(\L)$ is locally finite since it respects the direct limit structure of (\ref{dirlim}).
Using [Ku2], Cor.~1.3.10,
we can now see  that the action of $e_i$, $f_i$, $h_i$ on $\ti E(\L)$ combine to a $\gfr'_{KM}(G)$-module structure (cf.  [Ku1], Proposition 2.11).
In fact, $\ti E(\L)$ is an
integrable $\gfr'_{KM}(G)$-module
with a highest weight vector of multiplicity one and hence is irreducible. (This highest
weight vector is given by the evaluation of global sections of $\L$ at the ``origin" $w\in L^+Bw\subset \ti S_w$. Note that
[Ku1], Proposition 2.11 shows a similar statement for $\ti S_w$ replaced by the  Demazure varieties $D(\ti w)$; however, we can see,
using Proposition \ref{d2}, semicontinuity and the projection formula, that $\Gamma(D(\ti w), \L)\simeq \Gamma(\ti S_w, \L)$.
By \ref{compMathieu} our Demazure varieties for $G$  
coincide with the Kac-Moody Demazure varieties in loc. cit.) Similarly, we can consider the inductive limit
$E(\L):=\lim_{\to}(\Gamma(S_{w}, \L)^*)$; a similar argument shows  that this is also
an integrable $\gfr'_{KM}(G)$-module with a highest weight vector of multiplicity one. For each $w\in W_a$
we have a surjective homomorphism $\Gamma(\ti S_{w}, \L)^*\to \Gamma(S_w, \L)^*$
which respects the action of $L^+P_i$, obtained as the dual of the pull-back homomorphism.
Hence,
we obtain a surjective $\gfr'_{KM}(G)$-module homomorphism $\ti E(\L)\to E(\L)$;
by irreducibility this then has to be an isomorphism.
Since $\Gamma(\ti S_{w'}, \L)^*\to \Gamma(\ti S_{w}, \L)^*$ is injective for
$w'\leq w$, we can conclude that $\Gamma(\ti S_{w}, \L)^*\to \Gamma(S_w, \L)^*$ is an isomorphism for
each $w\in W_a$. The same argument applies to $\L^{\otimes n}$, $n>0$: We obtain
$\Gamma( S_{w}, \L^{\otimes n}) \simeq \Gamma(\ti S_w, \L^{\otimes n})$ from which we can easily deduce
$\ti S_w=S_w$, i.e., that $S_w$ is normal
(cf. [Fa1] or [Ku1]). 
%(Actually, from this it also follows that $\F_G$ coincides with
%the  affine flag variety for $\gfr_{KM}(G)$ defined via highest weight modules
%in [Ku1], [Ma1], [Li]; see also the explanation in Remark \ref{CompAtp} below). 

\subsection{}\label{final}
Here we complete the proof of Theorem \ref{normalmain}.
(The argument is due to Faltings [Fa2]). We continue to assume that 
${\rm char}(k)=p>0$.
Consider the morphisms
\begin{equation*}
\und D(\ti w)\xrightarrow{} \tiuS_w\xrightarrow{\und \psi_w} \und S_w
\end{equation*}
over the ring of Witt vectors $W$. Let $e_0: \Spec(W)\to \und S_w$ be  the point that
corresponds to the origin of $\und \F=L^+\und G/L^+\und B$. This lifts to the obvious
neutral point of $D(\ti w)$ and therefore also to a point of $\tiuS _w$. Denote by $A_w$, resp.
$\ti A_w$, $B_w$ the local ring of $\und S_w$, resp. $\tiuS_w$, $\und D(\ti w)$ at the corresponding neutral
points (lifting $\bar e_0:=e_0\otimes _Wk$) of the special fibers. We have
\begin{equation*}
A_w\hookrightarrow \ti A_w\hookrightarrow B_w
\end{equation*}
where
$B_w= W[t_1,\ldots ,t_r]_{(p, t_1,\ldots , t_r)}$. This ring supports a filtration by the ideals $F^n(B_w)=(t_1,\ldots , t_r)^n$; the associated graded pieces are free $W$-modules.
We set $F^n(A_w)=F^n(B_w)\cap A_w$ which is a filtration  of $A_w$ by ideals.
Let us denote by $x\in \und S_w(A_w)$ the tautological point $\Spec(A_w)\to \und S_w$;
it gives points $x_n\in \und S_w(A_w/F^n(A_w))$. For simplicity of notation, we set
$A^n_w:=A_w/F^n(A_w)$.
\medskip

\noindent{\bf Claim:}
{\sl For each $n\geq 1$, there is an element $g_n\in L\und G(A^n_w )$ which belongs to the subgroup generated by
$L^+\und B(A^n_w)$ and $\und U_\alpha(A^n_w)$ for $\alpha \in \Phi_{\rm aff}$ such that $x_n=g_n\cdot e_0$.}
\medskip

\begin{Proof}
We will argue using induction on $n$. For $n=1$ the claim is true since $A^1_w=W$ and $x_1=e_0$.
Suppose now the claim is true for $n$. Let us write $x_n=g_n\cdot e_0$; since $L^+\und B $ and $\und U_\alpha $
are smooth over $W$, we can lift $g_n$ to an element $g'_{n+1}$ which belongs to the subgroup generated by
$L^+\und B(A^{n+1}_w)$ and $\und U_\alpha(A^{n+1}_w)$. Consider $x'_{n+1}=g_{n+1}'\cdot e_0$;
by our construction we have
\begin{equation}
x'_{n+1}\ {\rm mod}\ F^n(A_w)\equiv x_n\ {\rm mod}\ F^n(A_w)\ .
\end{equation}
Since $\und S_w\subset \und \F=L\und G/L^+\und B$,
we can write $x_{n+1}=h_{n+1}\cdot e_0$. Now we can see that
the element $g'_{n+1}\cdot h_{n+1}^{-1}$ of $L\und G(A_w/F^{n+1}(A_w))$ satisfies
$g'_{n+1}\cdot h_{n+1}^{-1}\equiv 1\mod F^n(A_w)/F^{n+1}(A_w)$. Therefore
(cf. (\ref{lie2})), it is given
by an element of
\begin{equation*}
Lie(L\und G)\otimes_W F^n(A_w)/F^{n+1}(A_w)\ .
\end{equation*}
Now our claim follows from the fact (Proposition \ref{prop911}) that the Lie algebra $Lie(L\und G)$ is spanned over $W$
by the Lie algebras of $L^+\und B$ and $\und U_\alpha$ for all $\alpha\in \Phi_{\rm aff}$.
\endproof
\end{Proof}
\smallskip

Since the element $g_n\in L\und G(A^n_w )$ obtained in the Claim above acts on $\tiuS$, we
obtain the following: There is $w'>w$ in $W_a$ such that $x_n\in \und S_w(A^n_w)$
lifts to an element $\tilde x_n\in \tiuS_{w'}(A^n_w)$. Now recall that by \S \ref{KacM}, 
$\tiuS_w\otimes_WK_0=\und S_w\otimes_WK_0$. Hence, we have
\begin{equation*}
\ti x_{n}\otimes_WK_0\in \tiuS_{w'}((A^n_w)\otimes_WK_0)\cap \und S_w((A^n_w)\otimes_WK_0)=\tiuS_{w}((A^n_w)\otimes_WK_0)
\end{equation*}
which implies $x_n\in \tiuS_{w}(A^n_w )$. Furthermore, our arguments show that $\ti x_n$
are compatible when we vary the $n$. Hence, we obtain a splitting $\beta$ of the (injective) ring homomorphism
\begin{equation}
\alpha: \widehat \O_{\und S_w, \bar e_0}\hookrightarrow  \widehat \O_{\tiuS_w, \bar e_0}
\end{equation}
induced by $\psi_w: \tiuS_w\to \und S_w$ on the formal completions (i.e., $\alpha:= \hat\psi^*_w$).
We have $\beta\cdot \alpha={\rm id}$. Now consider
the coherent sheaf
\begin{equation*}
 (\psi_w)_*(\O_{\tiuS_w})/\O_{\und S_w}
\end{equation*}
on $\und S_w\subset \und \F$. This sheaf is $L^+\und B$-equivariant.
Using   \S \ref{KacM} we see that it is annihilated by a power of $p$.
Assume that its support is non-empty. Then it has to contain the unique fixed point $\bar e_0$
of the action of $L^+B$ (in characteristic $p$) and therefore,
\begin{equation*}
\widehat \O_{\tiuS_w, \bar e_0}/\alpha(\widehat \O_{\und S_w, \bar e_0})
\end{equation*}
is non-trivial and $p^m$-torsion for some $m>0$. For $f\in \widehat \O_{\tiuS_w, \bar e_0}$
now write
\begin{equation*}
\alpha\cdot \beta( p^m\cdot f)=\alpha\cdot \beta\cdot \alpha (f')=\alpha (f')=p^m\cdot f\ .
\end{equation*}
Since $\widehat \O_{\tiuS_w, \bar e_0}$ is flat over $W$, we obtain $\alpha\cdot \beta (f)=f$,
so $\alpha\cdot \beta={\rm id}$. This shows $\alpha$ is an isomorphism which is a contradiction.
Therefore, $(\psi_w)_*(\O_{\tiuS_w})=\O_{\und S_w}$ and so $\psi_w$ is an isomorphism.
This concludes the proof when ${\rm char}(k)=p>0$. The case that ${\rm char}(k)=0$
follows by semicontinuity, or more directly by the arguments in \S \ref{KacM}.\endproof
\medskip

\subsection{}\label{CompAtp}
 
Here we discuss the relation with the Schubert varieties in the
theory of Kac-Moody flag varieties. Suppose that $k$ is algebraically closed and $G$ is absolutely simple, simply 
connected and splits over a tamely ramified extension. As in \S \ref{KacM}
we can associate to $G$ (via its extended Dynkin diagram) a  Kac-Moody 
algebra ${\mathfrak g}:={\mathfrak g}_{KM}(G)$ over $k$. 
For $w\in W_a$ with reduced decomposition $\ti w$, there is a Demazure variety 
$D_{\mathfrak g}(\ti w) $ in the Kac-Moody setting (see [Ma1], p. 51) given by a contracted product
similar to the one in \ref{demazure}. By \S \ref{compMathieu}, $D_{\mathfrak g}(\ti w)$ is 
equivariantly isomorphic to the Demazure variety $D(\ti w)$ defined in \S \ref{demazure}. 
Now let $\lambda$ be a regular dominant integral 
weight of ${\mathfrak g}$  
and denote by $L(\lambda)$ the corresponding highest weight module over $k$
with highest weight vector $v_\lambda$. Set $S^{\mathfrak g}_{w,\lambda}$ 
for the Zariski closure of $B\cdot w\lambda$ in the (infinite dimensional) projective space
${\bf P}(L(\lambda))$; this is a projective scheme over $k$. We have a 
birational proper morphism of  $k$-schemes  
$$
p_{w}: D(\ti w)=D_{\mathfrak g}(\ti w)\xrightarrow{\ } S^{{\mathfrak g}}_{w,\lambda}\ .
$$
In this case, Mathieu [Ma1] and Littelmann [Li] prove that 
$S^{\mathfrak g}_{w,\lambda}$ is normal with rational 
singularities and as an abstract scheme is independent 
of $\lambda$. We can see by induction  that the $k$-schemes $S_w$ and 
$S^{\mathfrak g}_{w,\lambda}$ are stratawise isomorphic 
and so homeomorphic as topological spaces.
In fact (cf. [Ma1], proof of Lemme 33), we can find a homeomorphism 
$$
\tau_w:   S_w\xrightarrow{  \sim } S^{\mathfrak g}_{w,\lambda}
$$
such that $\tau_{w}\cdot \pi_{w}=p_w$. By [Ma1] and [Li], resp. Theorem \ref{normalmain}, the structure sheaf of $S^{\mathfrak g}_{w,\lambda}$,
resp. $S_w$, is given as ${p_w}_*(\O_{D(\ti w)})$, resp. ${\pi_w}_*(\O_{D(\ti w)})$.
Since a scheme is determined by the corresponding locally ringed topological space,
we see that $\tau_w$ actually gives an isomorphism of schemes. 
These isomorphisms combine to give an isomorphism of ind-schemes
\begin{equation}\label{twoflags}
LG/L^+B=\lim_{\to} S_w\xrightarrow{  \sim }\lim_{\to}S^{\mathfrak g}_{w,\lambda}\ .
\end{equation}
Recall that the first equality is Proposition \ref{limitSchu}.
On the other hand, $\displaystyle{\lim_{\to}S^{\mathfrak g}_{w,\lambda}=\lim_{\to}\ti S^{\mathfrak g}_{w,\lambda}}$
is by definition the Kac-Moody flag variety of [Ma1], [Li].

\bigskip

\section{ A coherence conjecture}\label{coherence}

\setcounter{equation}{0}

 In this section, we formulate a conjecture which relates amongst each other the
dimensions of the spaces of global sections of natural ``ample"  line bundles on the various
partial affine flag varieties attached to a fixed loop group. In the next section we demonstrate 
an example of the kind of application of this conjecture we have in mind.

We continue to assume that $k$ is algebraically closed and that $G$ is a simply connected semi-simple and 
absolutely simple  group over $K=k((t))$ which splits over the tamely ramified extension $K'$ of $K$.
Since this is the case we are mostly interested in, we also assume ${\rm char}(k)=p>0$. 
We continue to use the notation of the previous section; in particular
 $B$ is an Iwahori subgroup. Recall that we denote by 
$\bf S$ the corresponding set of simple affine roots.

\subsection{}We first explain the construction of certain ample line bundles.
\begin{prop}\label{Pic} Denote by $\F=LG/L^+B$  the full affine flag variety.
There is an isomomorphism 
$$
{\rm Pic} (\F)\xrightarrow {\ \sim\ } {\bf Z}^{\bf S}= \bigoplus_{i\in {\bf S}}\, {\bf Z}\epsilon_i\ .
$$
from the  Picard group of  line bundles on $\F$  to ${\bf Z}^{\bf S}$, defined by  the degrees of the restrictions to the projective lines corresponding to the simple affine roots.
\end{prop}
\begin{proof}  
We can assume that $B$ is an Iwahori subgroup given as in    \S \ref{9a}.
Let us recall that by the projective line ${\bf P}^1_{\alpha_i}\hookrightarrow \F$ that corresponds to the
simple affine root $\alpha_i$ we mean the image of $L^+P_{i}/L^+B$ in $\F$; this coincides with the image 
of ${\bf P}^1\to \F=LG/L^+B$ given by the homomorphism $\phi_{\alpha_i}: SL_2\to LG$ defined in \S \ref{9b}.
The main point now is that our previous results allow us to refer to the Kac-Moody context. Use \S \ref{KacM1} with the choice of positive roots 
and Chevalley basis the same to the choice that gives $B$ as in \S \ref{9a}.
Then, in view of (\ref{twoflags}), our statement is   given by [Ma1], XVIII, Prop. 28 (see also loc. cit. XII, Lemme 76
and its proof). These references
construct a line bundle $\L(\lambda)$ for each weight $\lambda\in (\oplus\ \Z\cdot h_i)^*$. Roughly speaking, 
the line bundle $\L(\lambda)$ can be identified with the homogeneous line bundle on the affine flag variety attached to the character $\lambda$ of the Iwahori subgroup in Kac-Moody theory. The degree of the restriction of $\L(\lambda)$ to the projective line corresponding to the simple affine root
$\alpha_i$ is then $\lambda(h_i)$ and so in our notation $\epsilon_i=h_i^*$.
\end{proof}
\begin{Remarks}{\rm
(i) By restricting a character $\lambda$ of the Iwahori subgroup as above to the center of the Kac-Moody central extension, we obtain a homomorphism
\begin{equation}
c: {\rm Pic} (\F)={\bf Z}^{\bf S}\longrightarrow {\bf Z}\ .
\end{equation}
The image $c(\L)$ of a line bundle $\L$ on $\F$ is called its {\it central charge}. It is given explicitly as follows.

Let us identify $\bf S$ with $\{0,\ldots, l\}$, and let $A=(a_{ij})_{i,j=0,\ldots,l}$ be the corresponding generalized Cartan matrix. By \cite{Kac}, (6.1.1), there exists an element $(a_0^\vee,\ldots,a_{l}^\vee)\in ({\bf Z}_{\geq 0})^{l+1}$ such that 
\begin{equation*}
(a_0^\vee,\ldots,a_{l}^\vee)\cdot A=(0,\ldots,0) \ ,
\end{equation*}
and such that $a_i^\vee=1$ for at least one $i$ with $0\leq i \leq l$. This element is obviously uniquely determined. Then 
\begin{equation}
c(\epsilon_i)=a_i^\vee\ ,\ i=0,\ldots, l \ .
\end{equation} 
Indeed, this follows from the formula $K=\sum\nolimits_i a_i^\vee\alpha_i^\vee$, \cite{Kac}, p.80 for the canonical generator of the center of ${\frak g}_{{\rm KM}}'(G)$, once $\epsilon_i$ is identified with $h_i^*$ as in the proof of Proposition \ref{Pic} above. 

(ii) Let $G$ be split and write $G=G_0\otimes_kK$. Let $P=G_0\otimes_kk[[t]]$ be the corresponding maximal parahoric subgroup. Let $T\subset G_0$ be a maximal torus and $B$ the Iwahori subgroup obtained as the inverse image of a Borel subgroup containing $T$ under the reduction map
\begin{equation*}
L^+P\longrightarrow G_0\ .
\end{equation*}
In \cite{Fa1}, p. 54,  Faltings constructs a central extension 
\begin{equation}\label{centralext}
1\longrightarrow {\bf G}_m\longrightarrow \tilde LG\longrightarrow LG\longrightarrow 1
\end{equation}
which acts on all line bundles on $\F$. The Lie algebra of $\tilde LG$ is the formal version of 
$\frak g_{\rm KM}'(G)$. We may consider $T$ as a subgroup of $L^+P$ ({\it constant loops}). Let $\tilde T$ be the inverse image of $T$ in $\tilde LG$. Then $X^*(\tilde T)$ may be identified with the character group of the Iwahori subgroup in Kac-Moody theory and we obtain a canonical isomorphism
\begin{equation}\label{ident}
X^*(\tilde T)\xrightarrow {\ \sim\ }  {\rm Pic}(\F) \ .
\end{equation}
Let $\{ \alpha_1,\ldots,\alpha_l\}$ be the set of simple roots of $T$ with respect to the fixed Borel subgroup, and let $\alpha_0=1-\theta$ be the remaining simple affine root. Then we may identify $X^*(T)$ with $\bigoplus\nolimits_{i=1}^l {\bf Z} \epsilon_i$, whereas $X^*(\tilde T)$ is identified with $\bigoplus\nolimits_{i=0}^l {\bf Z} \epsilon_i$ via Proposition \ref{Pic} and (\ref{ident}). In terms of these identifications, the canonical inclusion $\iota: X^*(T)\hookrightarrow X^*(\tilde T)$ is given as
\begin{equation}
\iota(\epsilon_i)=\epsilon_i-a_i^\vee \epsilon_0  \ ,\  i=1,\ldots,l \ .
\end{equation}
Indeed, this follows from the equalities
\begin{equation*}
\begin{aligned}
{\rm deg} (\L(\iota(\epsilon_i))\vert {\bf P}_{\alpha_i}^1)&=\delta_{ij} \ ,\ i,j=1,\ldots,l \ ,\\
c(\L(\iota(\epsilon_i))&=0 \ , \ i=1,\ldots,l \ .
\end{aligned}
\end{equation*}
The central extension (\ref{centralext}) has a unique splitting over $L^+P$. This also yields a splitting of the exact sequence
\begin{equation}
0\longrightarrow X^*(T) \xrightarrow{\ \iota\ } X^*(\tilde T)\xrightarrow {\ c\ } {\bf Z}\longrightarrow 0 \ .
\end{equation}
In terms of the identification $X^*(\tilde T)=\bigoplus\nolimits_{i=0}^l {\bf Z}\epsilon_i$,
this splitting is given by sending $1\in {\bf Z}$ to $\epsilon_0$. This follows from \cite{Fa1}, Thm. 7 (the last two bullets)
and the above. 

}
\end{Remarks}

\subsubsection{}\label{restrictLine}
 The result  of Proposition \ref{Pic} in the case of a split group 
$G=G_0\otimes_k K$ is shown in [Fa1], Cor. 12, by using the ``big cell" in $\F$.
Proposition \ref{Pic} can be shown in all cases, except when $G$ is of type $A^{(2)}_{2m}$, by appealing to the 
results of [Fa1]  without using the identification (\ref{twoflags}) of the two types
of affine flag varieties and the results of [Ma1]. In this paragraph,
we take some time to explain this claim; our arguments are also
useful for the rest of the paper.
 
To see the injectivity of the homomorphism in the statement of the proposition, consider a line bundle on $\F$ with zero restriction degrees. This line bundle then  has trivial restriction to any Demazure variety $D(\tilde w)$ and hence to any Schubert variety $S_w$, cf.~[Fa1], proof of Cor.~12. The injectivity assertion therefore follows from Prop.~\ref{limitSchu}. To show the rest of our claim,
we present $G$ as the invariants under $\sigma$ in $H'={\rm Res}_{K'/K}(H_{K'})$, where $H$ is a split group and where $\sigma=\sigma_0\otimes\tau$,
as in section \S7. 

In fact, we will show that under the restriction homomorphism
\begin{equation}
{\rm Pic}(\F_H)\to {\rm Pic}(\F)\hookrightarrow \Z^{\bf S}
\end{equation}
the image of ${\rm Pic}(\F_H)$ is the subgroup $\displaystyle{\oplus_{i\in {\bf S}}\Z\kappa(i)\epsilon_i}$. Here we set $\kappa(i)=2$ or $1$, according as the vector part of the affine root $\alpha_i$ is a multipliable root or not. Recall that $\kappa(i)=1$ for all $i$ unless $G$ is of type $A^{(2)}_{2m}$. 
It is enough to show that the image of the restriction homomorphism contains the elements $\kappa(i)\epsilon_i$.  
Recall from \S 9  that we fixed a Chevalley-Steinberg system relative to the splitting of $G$ over $K'$. We may write any affine root $\alpha$ in the form 
$$
\alpha=\nu(\alpha) + m\ ,\quad m\in \frac{1}{ [K':K]}\cdot \bf Z\ ,
$$
cf.~\S9. Here $a=\nu(\alpha)\in \Phi$ is a root of $(G, S)$. Now $a$ is the restriction of an absolute root $a'$ of $(G\otimes_KK', T')$, and the absolute roots $a'$ with this property form a single orbit $\{a'\}$ under $\sigma$. Recall also from \S 9 
the group homomorphism $\pi_a: G^a\to G$. Set $H^a_{K'}:=G^a\times_KK'$. We have $H^a=\prod_{\{a'\}}{SL_2}_{/K'}$ in case (I)
and $H^a_{K'}\simeq {SL_3}_{/{K'}}$ in case (II). The homomorphism $\pi_a$ induces injective homomorphisms
\begin{equation}
\pi'_a: H^a_{K'}\to H_{K'}\ 
\end{equation}
and we can see that these give corresponding injective maps between affine flag varieties
\begin{equation}
\pi'_a: \F_{H^a}=LH^a/L^+B^{H^a}\to \F_H=LH/L^+B^H\ .
\end{equation}
Note that by definition, ${\bf P}^1_{\alpha}\hookrightarrow \F_G\subset \F_H$ factors through 
$$
{\bf P}^1_{\alpha}\hookrightarrow \F_{G^a}\subset \F_{H^a}\ .
$$
Consider now the simple affine root $\alpha_i=a_i+m_i$ that corresponds to the index $i$.
Note that $\alpha'_i=a'_i+m_i$, for all $a'_i$ in the $\sigma$-orbit that 
correspond to $a_i$, are simple affine roots for $H_{K'}$. Appealing to the case of a split group, we find a line bundle on the affine flag variety $\F_H$ with restriction of degree $1$ to the projective line corresponding to one $\alpha'_i=a'_i+m_i$ in the $\sigma$-orbit, and degree $0$ to the projective lines corresponding to all other simple affine roots of $H_{K'}$. By restriction this induces a line bundle on $\F_G\subset \F_H$.
We claim that this line bundle has image  $\kappa(i)\epsilon_i$ in $\Z^{\bf S}$ under our homomorphism. 

Suppose that 
$\beta=b+m$ is a simple affine root of $G$. First assume that $b$ is not in the ray of the root $a$. Then 
by appealing to the injectivity of the homomorphism of the proposition applied to $H^b$ we can see that the restriction of this line bundle to $\F_{H^b}$ is trivial. Hence, the restriction of this line bundle
to ${\bf P}^1_{\beta}$ is also trivial. 

It remains to deal with the case that $b$ is in the ray of $a_i$. 
Note that the simple affine roots of $H^b$ can be identified with a subset of the set of simple affine roots of $H$.
 
First assume that neither $2a_i$ nor $a_i/2$ are  roots. Then $G^b=G^{a_i}= {\rm Res}_{L/K}SL_2$ 
with $L=K$ or $K'$. The embedding ${\bf P}^1_{\alpha_i}\hookrightarrow \F_{G^{a_i}}\subset \F_G\subset \F_{H}$ factors through
a commutative diagram
\begin{equation} 
\begin{CD}
{\bf P}^1_{\alpha_i} @>>>  \F_{G^{a_i}} \\
@V VV @VV V\\
\prod_{\{\alpha'_i\}}{\bf P}^1_{\alpha'_i} @>>> \prod_{\{a'_i\}} \F_{SL_2}\\
\end{CD}
\end{equation}
Here again $\{a'_i\}$ is the $\sigma$-orbit corresponding to $a_i$, $\{\alpha'_i\}=\{a'_i+m_i\}$, 
and $\F_{H^{a_i}}=\prod_{\{a'_i\}} \F_{SL_2}\subset \F_H$. The left vertical arrow is the diagonal embedding.
We can see that the pull-back of our line bundle on $\F_H$
to $\prod_{\{\alpha'_i\}}{\bf P}^1_{\alpha'_i}$ has degree $1$ on one of the factors and is trivial on all the 
others. Hence our claim follows.

Now assume that either $2a_i$ or $a_i/2$ is a root. Then $G^a=SU_3$, $H^a={SL_3}_{/K'}$ and we can reduce to  the case $G=SU_3$.
In this case, there are two simple affine roots $\beta_1$, $\beta_0$ which correspond to the $\sigma$-orbits
$\{\alpha_1+0,\alpha_2+0\}$ and $\{-(\alpha_1+\alpha_2)+\frac{1}{2}\}$ of affine roots of $SL_3$
(with the standard notation). The  homomorphism $\phi_{\beta_1}: SL_2\to LSU_3$ factors
through the ``constants"  $(SL_3)^{\sigma_0}=SO_3$. (This follows from the definition of $\phi_{\beta_1}$ in \S \ref{9b}
since in this case $m=0$.) 
Therefore the embedding ${\bf P}^1_{\beta_1}\hookrightarrow \F_{SU_3}$
factors through 
a commutative diagram
\begin{equation}\label{cd1}
\begin{CD}
{\bf P}^1_{\beta_1} @>>>  \F_{SU_3} \\
@V VV @VV V\\
SL_3/B @>>>\F_{SL_3}\\
\end{CD}
\end{equation}
and the left vertical morphism identifies ${\bf P}^1_{\beta_1}$ with the space
of isotropic flags. 

The homomorphism $\phi_{\beta_0}: SL_2\to LSU_3\to LSL_3$ agrees (up to an automorphism of the target) 
with the homomorphism $\phi_{\alpha_0}: SL_2\to  LSL_3$ associated to the simple affine root $\alpha_0 =-(\alpha_1+\alpha_2)+\frac{1}{2}$
of ${SL_3}_{/K'}$. This gives a
 commutative diagram
\begin{equation}\label{cd2}
\begin{CD}
{\bf P}^1_{\beta_0} @>>>  \F_{SU_3} \\
@V VV @VV V\\
{\bf P}^1_{\alpha_0} @>>>\F_{SL_3}\\ 
\end{CD}
\end{equation}
where the left vertical map is an isomorphism. 

If $\alpha_i=\beta_1$ then our line bundle restricts under the lower horizontal map in (\ref{cd1}) to one of the two ample generators of 
the Picard group of $SL_3/B$. Either of them has restriction of degree $2$ on 
the ${\bf P}^1$ of isotropic flags and so the degree of its restriction to ${\bf P}^1_{\beta_1}$ is $2$.
Using the commutative diagram (\ref{cd2}) we can see that its restriction to ${\bf P}^1_{\beta_0}$ is $0$.
This shows our claim in this case.  The case of $\alpha_i=\beta_0$ is similar (the restriction of the line bundle
to ${\bf P}^1_{\beta_0}$, resp. ${\bf P}^1_{\beta_1}$, has degree $1$, resp.  degree $0$).

\subsubsection{} It will be convenient for the rest of the paper to introduce the following notation
for the partial affine flag varieties associated to $G$. If $Y$ is a non-empty subset of the set of simple 
affine roots ${\bf S}$, we set $W^Y=W_{{\bf S}-Y}$, $P^Y=P_{{\bf S}-Y}$ and $\F^Y=LG/L^+P^Y=\F_{{\bf S}-Y}$.
Therefore, $\F=\F^{{\bf S}}$. 

Assume that the weight $\lambda = \sum n_i \epsilon_i$ is dominant, i.~e.~ $n_i\geq 0$ for all $i$.
It follows from its construction (see the references in the proof of Prop. \ref{Pic})
that the corresponding line bundle $\mathcal L(\lambda)$ comes by pull-back from an ample line bundle $\mathcal L(\lambda)^Y$ on  $\F^Y$
for $Y=Y(\lambda):=\{ \alpha_i\in {\bf S}\mid n_i>0\}$.
Here ample   is meant in the sense that the restriction of  $\mathcal L(\lambda)^Y$ to any Schubert variety is ample.

For $Y\subset {\bf S}$,  let us now consider the line bundles 
$$\L(Y)= \L(\sum_{i\in Y} \epsilon_i)^Y\ ,\quad
 \L'(Y)=\L(\sum_{i\in Y} \kappa(i)\epsilon_i)^Y
$$ on $\F^Y$. Let $Y=\{ i_1,\dots, i_s \}$ and consider the closed embedding of ind-schemes,
\begin{equation}\label{diagonal}
\F^Y \hookrightarrow \F^{i_1}\times \F^{i_2}\times \dots \times \F^{i_s}\ .
\end{equation}
Then $\mathcal L(Y)$ is the restriction of the line bundle 
\begin{equation}\label{exttens}
\mathcal L(i_1)\boxtimes \mathcal L(i_2)\boxtimes \dots \boxtimes \mathcal L(i_s)
\end{equation}
on $\F^{i_1}\times \F^{i_2}\times \dots \times \F^{i_s}$ and similarly for $\L'(Y)$.
Notice here that the restriction of $\mathcal{L}(Y)$  to a Schubert variety $S_w$ in 
$\F^Y$ is ample, and hence by the Frobenius spliting of $S_w$ we have
${H}^i(S_w, \mathcal{L}(Y))=0$ for $i>0$, and similarly for any positive tensor power of $\mathcal L(Y)$ and for any closed subset of $\mathcal{F}^Y$ which is a finite union of Schubert varieties
(or more generally $B$-orbits). These statements are also true   for $\L'(Y)$.

\begin{Remark}\label{104}
{\rm For $G= SL_n$ resp.~$Sp_{2r}$, the sheaf $\mathcal{L}(Y)=\L'(Y)$ is explicitly 
constructed in [Go2], 3.2.2. As we shall see below we can also obtain
a relatively explicit description for $\L'(Y)$ in the case of the ramified special unitary group of 
dimension $n$ as in \S \ref{unitary}. 
If  $n=2m+1$ is odd, resp. $n=2m$ is even, then we may identify the set of simple affine roots $\bf S$ with $\{0,\ldots, m\}$,
resp. $\{0,\ldots, m-2, m, m'\}$, 
by associating to $i\in \{0,\ldots, m\}$, resp.  $i\in \{0,\ldots, m-2, m, m'\}$, the unique simple reflection that {\sl does not} stabilize the lattice $\lambda_i$.
For $Y\subset {\bf S}$ let  $I=I(Y)$ be the corresponding subset of  $\{0,\ldots, m\}$, resp.  $i\in \{0,\ldots, m-2, m, m'\}$.  Recall the definition of $I^\sharp=\{i_0<\dots<i_k\}$ as in \ref{4b}, and let $\F'_{I^\sharp}$ be the corresponding variety of special
hermitian lattice chains (see (\ref{chains1}) and (\ref{chains2})). By Remark \ref{uniRem} the partial affine flag variety $\F^Y$ for $SU_n$ can be identified with a connected component of $\F'_{I^\sharp}$. 

 Given a Schubert variety 
$S_w$ in $\F'_{I^\sharp}$, we can find $N>0$ such that all corresponding lattices lie between $t^{-N}\lambda_0$ and $t^N\lambda_0$. We get a closed embedding 
\begin{equation}\label{Lexpl}
S_w \hookrightarrow \prod_{i_q\in I^\sharp} {\rm Grass}(t^{-N}\lambda_0/t^N\lambda_0, n_{i_q}),
\end{equation}
where $n_{i_q}={\rm dim}_k \lambda_{i_q}/t^N\lambda_0$.  The restriction of $\mathcal L'(Y)$ to $S_w$ is then the pullback of the exterior tensor product of the very ample generators of the Picard groups  of the various Grassmann varieties occurring in this product. Indeed, this follows from the corresponding expression for these line bundles in the case of $SL_n$ and the fact that we can obtain the line bundle $\L'(Y)$ 
as a restriction via the embedding $\F_{SU_n}\hookrightarrow \F_{SL_n}$ (cf. 
 \ref{restrictLine}. Here we view $SU_n$ as an outer form of $SL_n$.)}
\end{Remark}

\subsection{}
We now  let $T$ be the centralizer torus of $S$ in $G$ and denote by $T_{\rm ad}$ the image of $T$ in the adjoint group.  To $\mu\in X_\ast(T_{\rm ad})$ we attach its image $\lambda$ in the
coinvariants $X_\ast(T_{\rm ad})_I$. By (8.1) we can consider $\lambda$ as an element in the Iwahori-Weyl group $\tilde{W}_{\rm ad}$ of $G_{\rm ad}$.
After the choice of a special vertex $x$ in the chosen alcove we can write a commutative diagram
$$\xymatrix{
W_a\ar@{^{(}->}[d] & \simeq & W(^x\Sigma)\ltimes Q^{\vee}(^x\Sigma)\ar@{^{(}->}[d]\\
\tilde{W}_{\rm ad} & \hookrightarrow & W(^x\Sigma)\ltimes P^{\vee}(^x\Sigma)\quad .
 }
$$
Here $Q^{\vee}(^x\Sigma)$, resp.
$P^{\vee}(^x\Sigma)$, is the group of coroots, resp. of coweights, of the finite root system $^x\Sigma$, cf.~[R], \S 3 or (\ref{semidirect2}), and $\lambda \in P^{\vee}(^x\Sigma).$
In the sequel we write $W_0$ for the finite Weyl group $W(^x\Sigma)$. We  extend 
in the obvious way  the Bruhat 
order from $W_a$ to $\tilde{W}_{\rm ad} = W_a\rtimes\Omega_{\rm ad}$.

We consider  the {\it admissible subset of $\tilde{W}_{\rm ad}$ associated to the coweight} $\mu$, 
\begin{equation}
{\rm Adm} (\mu) = \{ w\in\tilde{W}_{\rm ad}\mid  w\leq w_0(\lambda):=w_0\lambda w_0^{-1}\, \text{ for some}\, w_0\in W_0 \}.
\end{equation}

The set ${\rm Adm} (\mu)$ only depends on the geometric conjugacy class of the one-parameter
subgroup $\mu$, cf.~[R], \S 3. Note that all elements of ${\rm Adm} (\mu)$ have the same image in $\tilde{W}_{\rm ad}/W_a = \Omega_{\rm ad}$.  Denoting by $\tau$ this common image of all elements of 
${\rm Adm} (\mu)$ in $\Omega_{\rm ad}$, we can define the subset ${\rm Adm} (\mu)^{\circ}$ of
$W_a$ by 
$$
{\rm Adm} (\mu)^{\circ}={\rm Adm} (\mu)\cdot \tau^{-1}
$$
so that 
${\rm Adm} (\mu)=\{ w \in \tilde{W}_{\rm ad}\mid w=w' \cdot \tau\ \  \text {with}\ \  w' \in {\rm Adm} (\mu)^{\circ} \}$.

For a non-empty subset $Y$ of the set of simple affine roots, let $Y^{\circ}\subset {\bf S}$ be the subset that corresponds 
to the set of simple reflections of the form $\{ \tau\cdot s_i\cdot \tau^{-1}\ |\ i\in Y\}$ where  $s_i$ is the 
reflection corresponding to the simple affine root parametrized by $i\in Y$.

We may define the subset  ${\rm Adm}^Y (\mu)$ of $\ti W_{\rm ad}$ by 
$W^Y\cdot {\rm Adm} (\mu)\cdot W^{Y}$ and
 the subset  ${\rm Adm}^Y (\mu)^{\circ}$ of
$W_a$ by
$$
{\rm Adm}^Y (\mu)^{\circ}:=W^Y\cdot {\rm Adm} (\mu)^{\circ}\cdot W^{Y^{\circ}}\ .
$$
Since $\tau\cdot W^Y\cdot \tau^{-1}=W^{Y^\circ}$ this is also equal to $
(W^Y\cdot {\rm Adm} (\mu)\cdot W^Y)\cdot \tau^{-1}$.

Let
\begin{equation}
\mathcal{A}(\mu)^{\circ} = \underset{w\in{\rm Adm}(\mu)^{\circ}}{\bigcup}S_w\subset\mathcal{F}, \,\quad
\mathcal{A}^Y(\mu)^{\circ} = \underset{w\in{\rm Adm}^Y(\mu)^{\circ}}{\bigcup}L^+B\cdot n_w\subset\mathcal{F}^{Y^0}
\end{equation}
be the closed reduced subset, union of all Schubert varieties (resp. $L^+B$-orbits) for $w$ in the $\mu$-admissible set (translated into the neutral component). 

We are interested in the dimension
\begin{equation}
h^{(\mu)}_Y(a) = \dim  {H}^0 (\mathcal{A}^Y(\mu)^{\circ}, \mathcal{L}'(Y^\circ)^{\otimes a})\, .
\end{equation}
This is given by a polynomial in $a$. By our previous remarks on the comparison between 
our set-up and the Kac-Moody theory, and using Littelmann's path model [Li], we can write
\begin{equation*}
h^{(\mu)}_Y(a) =\, \vert  \{\hbox{\rm LS-paths $\pi$  of shape}\ a\cdot (\sum_{\alpha_i\in Y^\circ}\kappa(i)\epsilon_i ) \mid i(\pi) \in {{\rm Adm}}^Y(\mu)^{\circ}/W^{Y^0}\}  \vert\,   .
\end{equation*}
Here  $i(\pi)$ denotes  the {\it initial direction} of $\pi$ in $W_a$. We see therefore that
the polynomial $h^{(\mu)}_Y$ has a purely combinatorial description. We want to give a (conjectural) explicit formula for the polynomials $h^{(\mu)}_Y$. 

Assume that $\mu\in X_*(T_{\rm ad})$ is given by a cocharacter $\mu: \Gm\to G_{\rm ad}\otimes_K{K^{\rm sep}}$  over the separable closure 
$ K^{\rm sep}$ of $K$
(here we abuse notation by using  the same symbol for the cocharacter of $G_{\rm ad}$). 
In fact, since $G$ is split over $K'$, $\mu$ is given by $\mu: \Gm\to H_{\rm ad}\otimes_KK'$
(the notation is as in \S \ref{tamepar}). 

Let us first suppose that $\mu$ is miniscule. 
Let $e=[K':K]$ and let $\{\mu=\mu^{(1)}, \ldots, \mu^{(e)}\}$ be 
the orbit of $\mu$ under the automorphism $\sigma_0$ of $H$.
Denote by $P(\mu^{(j)})$ a corresponding maximal parabolic subgroup of $H_{K'}=G_{K'}$. 
Let $X(\mu^{(j)})=H/P(\mu^{(j)})$ be the corresponding homogeneous space over $K'$, and let ${\mathcal L}(\mu^{(j)})$ be the ample generator of the Picard group of 
$X(\mu^{(j)})$. Note that 
\begin{equation}\label{confusion}
X(\mu)=G_{K'}/P(\mu)=(\prod_{j=1}^eX(\mu^{(j)})^{\sigma_0}\subset \prod_{j=1}^eX(\mu^{(j)})
\end{equation}
and that the restriction of $\boxtimes_{j=1}^e\L(\mu^{(j)})$ to this fixed point scheme is $\L(\mu)^{\otimes \, e}$.
 We introduce the polynomial in $a$,
\begin{equation}
h^{(\mu)}(a)={\rm dim}\ {H}^0(X(\mu), \mathcal L(\mu)^{\otimes \, e\cdot a})\ .
\end{equation}

Now assume that $\mu=\mu_1+\dots +\mu_r$ is a sum of minuscules. Then we set 
\begin{equation}
h^{(\mu)}=h^{(\mu_1)}\cdot h^{(\mu_2)}\cdot \dots \cdot h^{(\mu_r)}.
\end{equation}
Note that this is independent of the way $\mu$ is written as a sum of minuscules. Also note that 
$$h^{(\mu)}(a)={\rm dim}\ H^0(X(\mu_1)\times X(\mu_2)\times \dots \times X(\mu_r), (\mathcal L(\mu_1)\boxtimes \mathcal L(\mu_2)\boxtimes \dots \boxtimes \mathcal L(\mu_r))^{\otimes \, e\cdot a})\ .$$
\begin{Remark}
{\rm It is possible to give an explicit expression for the polynomial $h^{(\mu)}$.
For instance, if $G=SL_n$ and $\mu=\mu_r=(1^r, 0^{n-r})$ is a minuscule coweight, then $X(\mu)$  is the Grassmannian $Gr(r, n)$ of $r$-planes in a vector space
of dimension $n$, and $h^{(\mu)}(a)$ equals the degree of the irreducible representation of $PGL_n$ of highest weight $a\cdot (\epsilon_1 + \dots +\epsilon_r)$. }
\end{Remark}

\begin{conjecture}\label{littel}
Assume that $\mu$ is minuscule, or a sum of minuscule coweights. Let $Y\subset {\bf S}$. Then 
\begin{equation*}
h^{(\mu)}_Y(a) = h^{(\mu)}(\vert Y \vert \cdot a)\  ,\ \ \forall \ a= 1, 2,\dots
\end{equation*}
In particular, taking $Y = \{\alpha_i\}$ for  $\alpha_i\in {\bf S}$, the polynomial $h^{(\mu)}_{\{\alpha_i\}}(a)$ is independent 
of $\alpha_i$.
\end{conjecture}
\begin{Remark}
{\rm A somewhat weaker version of this conjecture would assert the existence  of a polynomial $h^{(\mu)}$ for which the above identity holds. It is this coherence statement that is really important for the applications of this conjecture we have in mind. We do not know whether the hypothesis on $\mu$ is necessary for this version of the conjecture to hold true.}
\end{Remark}

\begin{thm}\label{108}
Let $G=SL_n$ or $Sp_{2n}$. Then Conjecture \ref{littel} holds true.
\end{thm}
\begin{proof}
The essential additional ingredient here comes from the results of G\"ortz on 
flatness of local models for Shimura varieties, [Go1], [Go2].
Let us first consider the case where $G=SL_n$ and where  $\mu=(1^r, 0^{n-r})$ is minuscule.  Let us number the set of simple affine roots in the usual way by $s_0, \ldots , s_{n-1}$. Denoting by $K_0$ the fraction
field of $W$, let $\Lambda_i$ be the $W$-lattice in $V=K_0^n$ generated in terms of
the standard basis $e_1,\ldots, e_n$ by $p^{-1}e_1,\ldots, p^{-1}e_i, e_{i+1},\ldots,
e_n$. Then $s_i$ is the unique simple affine reflection that does {\it not} fix the lattice 
$\Lambda_i$. In what follows we will  identify a subset $Y$ of $\bf S$ with a subset of 
$\{0, 1, \ldots , n-1 \}$. 

Consider  the local model associated to $(GL_n, \mu)$ over $\Spec (W)$, and to $I\subset \{ 0,\ldots,
n-1\}$, cf. [Go1]. This is a flat scheme $M_I$ over $\Spec (W)$ with generic fiber the Grassmannian $GL_n/P(\mu)$ and with reduced special fiber. In fact, this special fiber is isomorphic to 
$\mathcal A^I(\mu)^{\circ}$  cf.~[Go1]. For $I=\{ i_1,\ldots, i_s\}$, we have a closed embedding
\begin{equation}\label{diag}M_I\hookrightarrow M_{i_1}\times_W\cdots\times_W M_{i_s}\ .
\end{equation}
(Here we simply write $M_i$ instead of $M_{\{i\}}$.)
Recall that $M_ i$ is a Grassmannian over $W$ which represents the   functor
$$
M_i(R)
=
\{ \mathcal E_i\subset \Lambda_i\otimes_WR\mid \mathcal E_i\ \mbox{locally
on}\ \Spec(R)
\
\mbox{a direct summand of rank}\ r\}\ .
$$
In the generic fiber, this closed embedding is the diagonal embedding of $GL_n/P(\mu)$ into the $s$-fold product of $GL_n/P(\mu)$ with itself.

In the special fiber this closed embedding is part of a commutative diagram, in which we denote by $\F^{I^{\circ}}$ the partial affine flag variety $LSL_n/L^+P^{I^{\circ}}$, 
\begin{equation}
\xymatrix{
\bar M_I\ar@{^{(}->}[d] & \hookrightarrow &\bar M_{i_1}\times\cdots\times \bar M_{i_s}\ar@{^{(}->}[d]\\
\mathcal F^{I^{\circ}} & \hookrightarrow & \mathcal F^{\{i_1\}^{\circ}}\times \cdots \times \mathcal F^{\{i_s\}^{\circ}}\quad .
 }
\end{equation}

{\bf Claim:} {\it The restriction of $\mathcal L'({I^{\circ}})=\mathcal L({I^{\circ}})$ to $ \bar M_I$ lifts to an invertible
sheaf ${\mathfrak L}({I^{\circ}})$ on $M_I$ with generic fiber isomorphic to the $\vert
I\vert$-th power of the ample generator $\mathcal L_0$ of the Picard group of the
Grassmannian.}

Since by (\ref{exttens}) we have $\mathcal L({I^{\circ}})=\mathcal L({\{i_1\}^{\circ}})\boxtimes\cdots \boxtimes \mathcal
L({\{i_s\}^{\circ}})\mid \mathcal F^{I^{\circ}}$, it suffices to consider the case where $I=\{ i\}$
consists of a single element. Then by (\ref{diag}) the claim will follow in general. 
Recall that  $r$ is defined by $\mu=(1^r,0^{n-r})$. The embedding $\bar M_i\hookrightarrow
\mathcal F^{\{i\}^{\circ}}$ is given by associating to the point of $\bar M_i$ represented by the submodule 
$\mathcal E_i$ of $\Lambda_i\otimes_WR=(\lambda_i\otimes_{k[[t]]}k)\otimes_kR$ the
lattice $L_{\mathcal E_i}$ inside ${\lambda_i}_R:=\lambda_i\otimes_{k[[t]]}R[[t]]$ obtained by taking the inverse
image of $\mathcal E_i$ under the surjective map
\begin{equation}
\lambda_i\otimes_{k[[t]]}R[[t]]\longrightarrow \lambda_i\otimes_{k[[t]]}R\ \ .
\end{equation}
Here $\lambda_i\subset k[[t]]^n$ is the $k[[t]]$-lattice generated by
$t^{-1}e_1,\ldots, t^{-1}e_i, e_{i+1},\ldots, e_n$. We obtain
\begin{equation*}
\bar M_i(R)= \{ L\mid L\  \text{lattice with}\ t\cdot {\lambda_i}_R \subset L\subset
{\lambda_i}_R \ \text{with}\ {\rm rk}\ ( L/ t\cdot {\lambda_i}_R) =r\}\ \ .
\end{equation*}
Notice that $P^{\{i\}^\circ}$ is the stabilizer of the $k[[t]]$-lattice  
$\lambda^r_i:=(1,\ldots, 1, t, \ldots ,t)\cdot \lambda_i$ where there are $n-r$ copies of $t$.
The  lattice $\lambda^r_i$ has the same $n$-th exterior power with $L$ and so, locally on $R$, we can write
$L=g\cdot \lambda^r_i$ with $g\in SL_n(R((t)))$. 
The restriction of $\mathcal L(\{i\}^{\circ})$ to this space of lattices has fiber
$\wedge^r(L/t{\lambda_i}_R)$ at the point represented by $L$, comp.~Remark \ref{104}
above. 
It follows that the line bundle ${\mathfrak L}(\{i\}^{\circ})$ on $M_i$ with fiber
$\wedge^r \mathcal E_i$ at the point represented by the submodule $\mathcal E_i$ of $\Lambda_i\otimes R$ is the desired lifting.

From the claim we obtain, by Frobenius splitness and semi-continuity
(for flat schemes), 
\begin{eqnarray*}
 h^{(\mu)}(\vert I\vert \cdot a)
 &=&
 \dim H^0(SL_n/P(\mu), \mathcal L_0^{\otimes\vert I\vert a})\cr
& = & 
\dim H^0(M_I\otimes_WK_0, {\mathfrak L}({I^{\circ}})^{\otimes a})\cr 
& = &
 \dim H^0(M_I\otimes_Wk, \mathfrak L({I^{\circ}})^{\otimes a})\cr
&=&
 \dim H^0(\mathcal A_I(\mu)^{\circ}, \mathcal L({I^{\circ}})^{\otimes a})\cr
&=&
 h_I^{(\mu)}(a)\ ,
\end{eqnarray*}
which proves the conjecture in this case. 

More generally, let $G=SL_n$ and let $\mu=\mu_1+\cdots +\mu_m$ be a sum of minuscule coweights. (Note that, in fact, any coweight of $PGL_n$ is
a sum of minuscule coweights.) In this case let $K'$ be a totally ramified extension of degree $m$ of $K_0={\rm Frac}( W)$, and consider the
local model $M_I$ attached to the group $G=\Res_{K'/K_0}(GL_n)$, the coweight $(\mu_1,\ldots, \mu_m)$, and to $I\subset\{ 0,\ldots, n-1\}$,
comp.\ \cite{P-R2}. Then again $M_I$ is a flat $W$-scheme whose special fiber can be embedded in the partial flag variety $LGL_n/ L^+P^I$. In
fact, the special fiber is reduced and can be identified with $\mathcal A_I(\mu)^{\circ}$. In this case the {\it geometric} generic fiber of $M_I$ is the product of the Grassmannians $SL_n/P(\mu_1),\ldots, SL_n/P(\mu_m)$.
We claim that the same reasoning as before yields the formula
\begin{equation}\label{sumofminusc}
h_I^{(\mu)}(a)= \prod\limits_{j=1}^m \dim H^0 (SL_n/P(\mu_j), \mathcal L_j^{\otimes \vert I\vert
a})\ \ .
\end{equation}
Here $\mathcal L_j$ denotes the ample generator of the Picard group of
$SL_n/P(\mu_j)$.

To see this, we may again assume that $I=\{i\}$ consists of a single element. Then
$M_i$ is contained in the naive local model $M_i^{\rm naive}$. Let $E$ be the reflex
field for $(G,(\mu_1,\ldots, \mu_m))$, cf.~\cite{P-R2}. Then a point of $M_i^{\rm naive}$ with values in
an $\mathcal O_E$-algebra $R$ is given by
\begin{eqnarray*}
M_i^{\rm naive}(R)
&=&
\{ \mathcal E\subset \Lambda_i\otimes_WR\mid \mbox{a}\ \mathcal
O_{K'}\otimes_WR\mbox{-submodule which is locally on Spec($R$})
\\
&&
\mbox{a direct summand as $R$-module with}\ \det(a\vert\mathcal E)
= \prod\limits_j\varphi_j(a)^{r_j}, \forall a\in \mathcal O_{K'}\}\ \ .
\end{eqnarray*}
Here $\Lambda_i\subset K'^n$ denotes the $\mathcal O_{K'}$-lattice generated by
$\pi^{\prime -1}e_1,\ldots, \pi^{\prime-1}e_i, e_{i+1},\ldots, e_n$, where $\pi'$ is uniformizer
of $\O_{K'}$, and
$\varphi_1,\ldots,\varphi_m$ denote the various $K_0$-embeddings of $K'$ into a fixed
algebraic closure $\bar K_0$ of $K_0$ and $r_j$ is defined by $\mu_j=(1^{r_j},0^{n-{r_j}})$.
The special fiber of $M_i^{\rm naive}$ is embedded into $\mathcal F^{\{i\}^{\circ}}$ by
associating to a point represented by the submodule  $\mathcal E$ of
$$\Lambda_i\otimes_WR= \Lambda_i\otimes_{\mathcal O_{K'}}(\mathcal O_{K'} /p\mathcal
O_ {K'})\otimes_kR= (\lambda_i/t^m\lambda_i)\otimes_kR$$
the lattice $L_{\mathcal E}$ equal to the inverse image of $\mathcal E$
under the natural map ${\lambda_i}_R\to {\lambda_i}_R/t^m{\lambda_i}_R$. Here again 
${\lambda_i}_R$ denotes $\lambda_i\otimes_{k[[t]]}R[[t]]$. Hence $L_{\mathcal E}$ is squeezed as
\begin{equation}\label{squeeze}
t^m{\lambda_i}_R\subset L_{\mathcal E}\subset {\lambda_i}_R\ \ .
\end{equation}
The value of $\mathcal L^{\{i\}^{\circ}}$ on $L_{\mathcal E}$ is $\wedge^{\rm max}(L_{\mathcal E}/t^m{\lambda_i}_R)$. Consider the line bundle ${\mathfrak L}^{\rm naive}$ on
$M^{\rm naive}$ with fiber $\wedge^{\rm max}(\mathcal E)$ at the point represented by the submodule $\mathcal
E$. If $\mathcal E$ represents a point of the geometric generic fiber, then
\begin{equation}
\mathcal E=\bigoplus\limits_{j=1}^m\mathcal E_{\varphi_j}\ \ ,
\end{equation}
where $K^{\prime n}\otimes_{K_0}\bar K_0=\bigoplus (\bar K_0^n)_{\varphi_j}$ and
$\mathcal E_{\varphi_j}:=\mathcal E\cap (\bar K_0^n)_{\varphi_j}$ is a subspace of
dimension $r_j$. Now
\begin{equation}
\wedge^{\rm max}\mathcal E=\wedge^{\rm max}(\mathcal E_{\varphi_1})\otimes\cdots\otimes \wedge^{\rm max} (\mathcal
E_{\varphi_m})\ \ ,
\end{equation}
and this is the fiber at $\mathcal E=(\mathcal E_{\varphi_1},\ldots, \mathcal
E_{\varphi_m})$ of the exterior tensor product of the line bundles $\mathcal
L_1,\ldots, \mathcal L_m$ on the factor Grassmannians. It follows that the
restriction ${\mathfrak L}$ of ${\mathfrak L}^{\rm naive}$ to the local model is a lifting of 
the restriction of ${\mathcal L}^{\{i\}^{\circ}}$ to the special fiber of $M_i$, and has geometric generic fiber equal to the exterior tensor product of ${\mathcal
L}_1,\ldots, {\mathcal L}_m$. The formula (\ref{sumofminusc}) above now follows as in the minuscule case.

The case $G=Sp_{2n}$, and when $\mu$ is a multiple of the unique minuscule coweight is analogous. Instead of the local models for $\Res_{K'/K_0}( GL_n)$ one uses the local models for $\Res_{K'/K_0}(GSp_{2n})$, cf.~[Go2] and [P-R2].
\end{proof}
\begin{Remark}
{\rm In the previous proof we used the local models for the group $GL_n$ over $K_0$, resp. for 
${\rm Res}_{K'/K_0} ( GL_n)$. In fact, the same proof works, if  $K_0$ is replaced by the
fraction field of any complete discrete valuation ring with residue field $k$. In particular, one can take this dvr of equal characteristic, e.g. $k[[t]]$.}  
\end{Remark}

\section{Local models for the ramified unitary groups}
\setcounter{equation}{0}

In this section we will show how to embed the special fiber of
the local model attached to a ramified unitary group into
an affine flag variety. The structure results of the
preceding sections have then interesting consequences for
these local models (see also [P-R3]).

We use a  notation that is modeled on that of section \ref{unitary}. Let $F$ be a complete
discretely valued field with ring of integers $\mathcal
O_F$ and perfect residue field $k$ of characteristic $\neq
2$. Let $F'/F$ be a ramified quadratic extension and let
$\pi'\in F'$ be a uniformizer with $\overline{\pi'}=-\pi'$. Let
$V$ be a $F'$-vector space of dimension $n\geq 2$ and let
$$\phi:\  V\times V\longrightarrow F'$$
be a $F'/F$-hermitian form.  We assume that
$\phi$ is split. This means that there
exists a basis $e_1,\ldots, e_n$ of $V$ such that
$$\phi(e_i,e_{n-j+1})= \delta_{ij}\ \ ,\ \ \forall\  i,j=1,\ldots, n\  .$$
As in section \ref{unitary} we have two associated $F$-bilinear forms,
$$(x,y)=\Tr_{F'/F}(\phi(x,y))\ \ ,\ \ \langle x,y\rangle = \Tr_{F'/F} (\pi'^{-1}\cdot\phi (x,y))\ \ .$$
For any $\mathcal O_{F'}$-lattice $\Lambda$ in $V$ we set, as
in section \ref{unitary},
$$\hat\Lambda =\{ v\in V\mid \phi (v,\Lambda)\subset
\mathcal O_{F'}\} =\{ v\in V\mid \langle v,\Lambda\rangle 
\subset \mathcal O_{F'}\} \ .$$
Similarly, as in section \ref{unitary}, we set
$$\hat\Lambda^s =\{ v\in V\mid (v,\Lambda)\subset \mathcal O_F\}\  ,$$
so that $\hat\Lambda^s =\pi'^{-1}\cdot\hat\Lambda$. For
$i=0,\ldots, n-1$, set
$$\Lambda_i={\rm span}_{\mathcal O_{F'}} \{ \pi'^{- 1}
e_1,\ldots, \pi'^{- 1} e_i, e_{i+1},\ldots, e_n\}\ \
.$$ Then we may complete these lattices to a periodic
lattice chain by putting for $j\in\bf  Z$ of the form
$j=kn+ i$ with $0\leq i<n$,
$$\Lambda_j=\pi^{\prime -k}\cdot \Lambda_i\ \ .$$
This lattice chain is selfdual and periodic. More precisely, we have
$\hat\Lambda_j=\Lambda_{-j}$, for all $j$, and $\pi' \cdot \Lambda_j=\Lambda_{j-n}$.

Now we fix non-negative integers $r,s$ with $n=r+s$. We set
$E=F'$ if $r\neq s$ and $E=F$ if $r=s$ (this is the reflex
field of the local model we are about to define). We now
define as follows a functor $ M^{\rm naive}$ on the
category of $\mathcal O_E$-schemes. A point of $
M^{\rm naive}$ with values in an $\mathcal O_E$-scheme $S$
is given by a $\mathcal O_{F'}\otimes_{\mathcal
O_F}\mathcal O_S$ submodule
\begin{equation}\label{Mnaive}
\mathcal E_j\subset \Lambda_j\otimes_{\mathcal
O_F}\mathcal O_S
\end{equation}
for each $j\in\bf  Z$. The following
conditions are imposed.
\begin{itemize}
\item[a)] as an $\mathcal O_S$-module, $\mathcal E_j$ is
locally on $S$ a direct summand of rank $n$.
\item[b)] for each $j<j'$, there is a commutative diagram
$$
\begin{matrix}
\Lambda_j\otimes_{\mathcal O_F}\mathcal O_S
&
\longrightarrow
&
\Lambda_{j'}\otimes_{\mathcal O_F}\mathcal O_S
\\
\cup
&&
\cup
\\
\mathcal E_j & \longrightarrow & \mathcal E_{j'}
\end{matrix}
$$
and for each $j$, the isomorphism $\pi': \Lambda_j\to\Lambda_{j-n}$ induces an isomorphism of  $\mathcal E_j$ with $\mathcal E_{j-n}.$
\item[c)] we have $\mathcal E_{-j}=\mathcal E_j^\bot$ where $\mathcal E_j^\bot$
is the orthogonal complement of $\mathcal E_j$ under the
perfect pairing
$$(\Lambda_{-j}\otimes_{\mathcal O_F}\mathcal O_S)\times
(\Lambda_j\otimes_{\mathcal O_F}\mathcal
O_S)\longrightarrow \mathcal O_S$$ given by $\langle\ ,\
\rangle\otimes_{\mathcal O_F}\mathcal O_S$.
\end{itemize}

 Next note that $\mathcal E_j$ is an $\mathcal
O_{F'}\otimes_{\mathcal O_F}\mathcal O_S$-module, hence
$\mathcal O_{F'}$ and $\mathcal O_E$  act on it. We require
further that
\begin{itemize}
\item[d)] for each $j$, the characteristic polynomial of the action of $\pi'$
equals
$$\det((T\ {\rm id}- \pi')\mid \mathcal E_j)= (T-\pi')^r\cdot
(T+\pi')^s\subset \mathcal O_E[T]\ \ .$$
\end{itemize}
This concludes the definition of the functor $
M^{\rm naive}$, which is obviously representable by a
projective scheme over $\Spec (\mathcal O_E)$. We call
$ M^{\rm naive}$ the {\it naive} local model
associated to the group $GU(V,\phi)$, the signature type
$(r,s)$ and for the complete lattice chain $\Lambda_j$,
$j\in\bf  Z$. Note that the geometric generic fiber of $
M^{\rm naive}$ is isomorphic to the Grassmannian of $r$-planes in $n$-space. 

 We can generalize this definition to
incomplete selfdual periodic lattice chains as follows.  Let $I$ be a non-empty subset of the set of simple roots, which we identify as in Remark \ref{104} with a subset of  $\{ 0, \ldots ,  m \}$ in case $n=2m+1$ is odd, resp. with a subset of $\{ 0, \ldots ,m-2,  m, m'\}$ in case $n=2m \geq 4$ is even. We also recall that we associated to  $I$ a totally ordered subset $I^{\sharp}$ of  $\{ 0, \ldots ,m-1,  m \}$ resp. 
 $\{ 0, \ldots ,m-1,  m, m'\}$. Here $I=I^{\sharp}$ unless $n$ is even and  $\{m, m'\}\subset I$. In the last case, $I^{\sharp}$ is obtained from $I$ by replacing $m'$ by $m-1$ and keeping all other elements. It is obvious that, conversely,  $I^{\sharp}$ determines uniquely  $I$. In order to obtain transition morphisms between naive local models in all cases, we will only consider subsets $I$ which in case $n=2m$ have the property that if $m'\in I$, then $\{ m, m' \}\subset I$. (We refer to Remark \ref{uniRem} C for an explanation of why we do not loose much information by only considering subsets $I$ with this property.) Now  $I$ defines an
incomplete lattice chain by including all $\Lambda_j$ in
the chain with $j$ of the form $j=kn\pm i$ with $i\in I^{\sharp}$.
Correspondingly, we obtain a functor $ M_{I^{\sharp}}^{\rm
naive}$ by only giving the submodules $\mathcal E_j$ for
$j$ in this subset of $\bf  Z$. Therefore the case of complete lattice chains is obtained as $M^{\rm naive}=
M^{\rm naive}_{\{0,\dots ,m\}}.$
We have forgetful
morphisms of projective schemes
\begin{equation}\label{11.1}
 M^{\rm naive} \longrightarrow  M_{I^{\sharp}}^{\rm
naive}\ ,\ \mbox{and}\ \  M_{{I'^{\sharp}}}^{\rm
naive}\longrightarrow  M_{I^{\sharp}}^{\rm naive}\ \mbox{for}\
I'\supset I\ \ .
\end{equation}

We now construct 
natural embeddings of the special fibers $\bar{
M}_{I^{\sharp}}^{\rm naive}= M_{I^{\sharp}}^{\rm
naive}\otimes_{\mathcal O_E}k$ into a
variety of hermitian lattice chains of the type  considered in
section \ref{unitary}. For this we fix identifications
compatible with the actions of $\pi'$, resp.\ $u$,
$$\Lambda_j\otimes_{\mathcal O_F}k=
\lambda_j\otimes_{k[[t]]}k$$ which sends the natural bases
of each side to one another. We therefore also obtain a
$k[[u]] / (u^2)$-module chain isomorphism
\begin{equation}
\Lambda_\bullet\otimes_{\mathcal O_F}k\simeq
\lambda_\bullet\otimes_{k[[t]]}k\ \ ,
\end{equation}
 which is in fact an
isomorphism of polarized periodic module chains [RZ]. Let $R$
be a $k$-algebra. For an $R$-valued point of $
M_{I^{\sharp}}^{\rm naive}$ we have
$$\mathcal E_j\subset \Lambda_j\otimes_{\mathcal O_F}R=
(\lambda_j\otimes_{k[[t]]}k)\otimes_k R , \quad  j=kn\pm i, \, i\in I^{\sharp} \\ .$$ Let
$ L_j := u^{-1}L_{\mathcal E_j}$ where $L_{\mathcal E_j}\subset\lambda_j\otimes_{k[[t]]} R[[t]]$ is the
inverse image of $\mathcal E_j$ under the canonical
projection
\begin{equation}
\lambda_j\otimes_{k[[t]]}R[[t]]\longrightarrow
\lambda_j\otimes_{k[[t]]}R\ \ .
\end{equation}
 If $I^{\sharp}=\{i_0<i_1<\ldots < i_k\}$, then
$$\ldots\subset L_{i_0}\subset  L_{i_1}\subset\ldots
\subset u^{-1} L_{i_0}=  L_{i_0+n}\subset\ldots$$
is a periodic lattice chain in $R((u))^n$ which satisfies
conditions a) and b) in \S \ref{4b}, and which therefore defines an $R$-valued point of $\F_{I^{\sharp}}$.  We obtain in this way a morphism
\begin{equation}
\iota_{I^{\sharp}}: M_{I^{\sharp}}^{\rm naive}\otimes_{\mathcal
O_E}k\hookrightarrow \mathcal F_{I^{\sharp}}\ \ ,
\end{equation}
 which is a closed
immersion (of ind-schemes). 

For $I'\supset I$, we have a commutative diagram
$$
\begin{matrix}
 M_{I'^{\sharp}}^{{\rm naive}}\otimes_{\mathcal O_E}k & \hookrightarrow &
\mathcal F_{ I'^{\sharp}}
\\
\big\downarrow && \big\downarrow
\\
 M_{I^{\sharp}}^{{\rm naive}}\otimes_{\mathcal O_E}k & \hookrightarrow &
\mathcal F_{ I^{\sharp}}
\end{matrix}\ \ .
$$
The horizontal morphisms are equivariant for the actions of the stabilizer groups of \S \ref{unitary}, 
$L^+P_{I'}$ resp.\ $L^+P_I$, in the sense of \cite{P-R2},
section 6.

Examples show [P] that $ M_{ I^{\sharp}}^{\rm naive}$ is not
flat over $\mathcal O_E$ in general. We define $
M_{ I^{\sharp}}^{\rm loc}$ to be the flat closure of the generic fiber
$ M_{ I^{\sharp}}^{\rm naive}\otimes_{\mathcal O_E}E$ in
$ M_{ I^{\sharp}}^{\rm naive}$. Since the generic fibers
$ M_{ I^{\sharp}}^{\rm naive}\otimes_{\mathcal O_E}E$ are all
identical, we obtain for $I'\supset I$ commutative diagrams
of projective morphisms resp.\ closed embeddings,
$$
\begin{matrix}
 M_{ I'^{\sharp}}^{\rm loc} & \hookrightarrow & 
M_{ I'^{\sharp}}^{\rm naive}
\\
\big\downarrow && \big\downarrow
\\
 M_{ I^{\sharp}}^{\rm loc} & \hookrightarrow & 
M_{ I^{\sharp}}^{\rm naive}
\end{matrix}
\  ,
$$
whose generic fibers are all the identity morphism. Note that, since the generic fiber of 
$ M_{ I^{\sharp}}^{\rm loc}$ is connected, so is its special fiber. 

\begin{Remark}{\rm As we have seen in section \ref{unitary} the stabilizer group $P_I$ is not always a parahoric subgroup of $U_n$, but that its intersection $P'_I$ with $SU_n$ always is a parahoric subgroup. Furthermore, the use of the lower index is ``complementary" to its use elsewhere in the paper (e.g., if  $I=\{ 0, \ldots ,  m \}$ resp. 
$I=\{ 0, \ldots ,m-2,  m, m'\}$, then $P_I$ is an Iwahori subgroup of $U_n$).  To avoid a notational conflict, we will  denote for the rest of the paper by $P^I$ the parahoric subgroup of $SU_n$ corresponding to the subset $I$ of the set of simple affine roots, as in \S \ref{coherence} above.  The corresponding parahoric subgroup of $U_n$ will be denoted by $P^I(U_n)$.

}
\end{Remark}

\begin{lemma}\label{lemma11.1}
Let $\mu$ be the geometric cocharacter of $U_n\otimes_K\bar
K\simeq GL_n$ given by $\mu=(1^r, 0^s)$, Let $\mathcal A^I(\mu)=\bigcup_{w\in W^I\backslash {\rm Adm}^I(\mu)/W^I}S_w$, a
closed reduced subset of $LU_n/L^+P^I(U_n)$. \footnote{In \S\ref{coherence},  ${\rm Adm}^I(\mu)$  is only defined in the case of an   adjoint group, but the general definition is the same.}  
Then $\mathcal A^I(\mu)$ is contained in the image
of $ M_{I^{\sharp}}^{\rm loc}\otimes_{\mathcal O_E}k$ under
$\iota_{I^{\sharp}}$.
\end{lemma}

\begin{proof}
It  suffices to prove this for $I=\{ 0,\ldots,
m\}$, resp. $I=\{ 0,\ldots,
m-2, m, m'\}$. In this case, it suffices to lift the points in an open subset of each
maximal stratum $S_w$  to the generic fiber of $ M^{\rm loc}$, where
$w=w_0(\lambda)$ for some $w_0\in W_0$. The   proof of this statement can be found in [P-R3], \S 3.d.
\end{proof}

\begin{thm}\label{theorem11.2} Denote by $\mu_r$ the image in the adjoint group of the cocharacter $(1^r, 0^s)$ of $U_n$. 
Assume the validity of Conjecture \ref{littel} for the pair
$(G,\mu)=(SU_n, \mu_{r})$ and for $I^{\circ}$. Then the special fiber of
$ M_{I^{\sharp}}^{\rm loc}$ is reduced and in fact isomorphic to
$\mathcal A^I(\mu)^{\circ}$. Its irreducible components are normal
and with rational singularities.
\end{thm}
\begin{proof}
As already mentioned above,  $\iota_{I^{\sharp}}$ is
equivariant for the action of $L^+P_I$. Hence the image of
$\iota_{I^{\sharp}}$ is set-theoretically a union of Schubert
varieties. After identifying the special fiber
$\bar{ M}_{I^{\sharp}}^{\rm loc}$ of $ M_{I^{\sharp}}^{\rm loc}$
with its image under $\iota_{I^{\sharp}}$, we have by Lemma
\ref{lemma11.1} a closed embedding
$\mathcal A^I(\mu)\subset \bar{ M}_{I^{\sharp}}^{\rm loc}$. We can identify $\mathcal A^I(\mu)$ with the reduced closed subset $\mathcal A^I(\mu)^{\circ}$ of the 
partial affine flag variety $LSU_n/P^{I^{\circ}}$, and similarly for $\bar{ M}_{I^{\sharp}}^{\rm loc}$ (note that   $\bar{ M}_{I^{\sharp}}^{\rm loc}$ is connected and so by our results in \S \ref{reduced},
$(\bar{ M}_{I^{\sharp}}^{\rm loc})_{\rm red}$ lies in a translate of $LSU_n/P^{I^{\circ}}$ inside $\F_{I^{\sharp}}$, cf. Proposition \ref{GG'prop}).  We consider the natural ample line bundle $\mathcal
L'({I^{\circ}})$ on $LSU_n/P^{I^{\circ}}$  and denote by the same symbol its
restrictions to $\mathcal A^I(\mu)^{\circ}$ and $(\bar{
M}_{I^{\sharp}}^{\rm loc})_{\rm red}$. We will show that the restriction of $ {\mathcal L}'({I^{\circ}})$ to $(\bar{ M}_{I^{\sharp}}^{\rm loc})_{\rm red}$ lifts to a line bundle ${\mathfrak L}'(I^{\circ})$ on ${ M}_{I^{\sharp}}^{\rm loc}$ with geometric generic fiber isomorphic to the $2\vert I\vert$-fold tensor power of the ample generator $\mathcal L(\mu_r)$ of the Picard group of  the Grassmannian of $r$-planes in $n$-space.  Then, to prove that the  inclusions $\mathcal A^I(\mu)^\circ\subset 
(\bar{
M}_{I^{\sharp}}^{\rm loc})_{\rm red}\subset \bar{ M}_{I^{\sharp}}^{\rm loc}$ are
equalities, it suffices to prove that
\begin{equation}\label{claim}
\dim H^0 (\mathcal A^I(\mu)^{\circ}, {\mathcal L}'({I^{\circ}})^{\otimes
a})=\dim H^0(\bar{ M}_{I^{\sharp}}^{\rm loc}, \bar {\mathfrak
L}'({I^{\circ}})^{\otimes a})
\end{equation} for all $a>>0$, or also  $\chi(\mathcal A^I(\mu)^{\circ}, {\mathcal
L}'({I^{\circ}})^{\otimes a})=\chi (\bar{ M}_{I^{\sharp}}^{\rm
loc}, \bar {\mathfrak L}'({I^{\circ}})^{\otimes a})$ (equality of the
Hilbert polynomials). This follows by Conjecture \ref{littel} for $I^{\circ}$, 
since by the flatness of ${ M}_{I^{\sharp}}^{\rm loc}$ over $\O_E$ we have
$$
 \chi(\bar{ M}_{I^{\sharp}}^{\rm loc}, {\mathfrak L}(I^{\circ})^{\otimes
a})= \chi( M_{I^{\sharp}}^{\rm loc}\otimes_{\mathcal O_E}\bar
E, \mathcal L(\mu_r)^{\otimes 2a\vert I\vert})=h^{(\mu)}(\vert I \vert\cdot a)\ .
$$
It remains to construct the line bundle ${\mathfrak L}'={\mathfrak L}'(I^{\circ})$ and for this we may  replace the local model by the naive local model. 
Furthermore, by the Kunneth type argument used in the proof of (\ref{108}) we may assume that  $I^{\sharp}$ consists of a single element $\{i\}$.\footnote{or, that $n=2m$
and  that $I^{\sharp}=\{m-1, m\}$. This case is handled in the same way.} We claim in this case that the line bundle 
 ${\mathfrak L}'$ with value $\wedge^{\rm max}\, \mathcal E$ in the point of $ M_{\{i\}}^{\rm naive}$ represented by the subspace $\E$ of 
 $\Lambda_i\otimes_{\O_F} \mathcal O_S$ has the required property. Let us calculate this line bundle on the special fiber and the geometric generic fiber.
 
 If $\E$ represents a point of the special fiber which is mapped to the lattice $u^{-1} L_{\E}$ under $\iota_{\{i\}}$, then we have an inclusion of lattices as in (\ref{104}),
 $$
 u^N\, \lambda_i \subset u^{-1}L_{\E} \subset u^{-N}\, \lambda_i \, .
$$
By  Remark \ref{104}, the value  of  ${\mathcal L}'(\{i\}^{\circ})$  at $u^{-1}L_{\E}$ is equal to $\wedge^{\rm max}\, (u^{-1} L_\mathcal E/ u^N \, \lambda_i)$. Hence the restriction $\iota_{\{i\}} ^*{\mathcal L}'(\{i\}^{\circ})$ is indeed 
isomorphic to the restriction of ${\mathfrak L}'$ to the special fiber $\bar { M}_{\{i\}}^{\rm naive} $. 

Let us now consider the pullback of ${\mathfrak L}'$ to the geometric generic fiber of ${ M}_{\{i\}}^{\rm naive}$. Now
$$
V\otimes_F F'=V_1\oplus V_2 \, ,
$$
according to the two $F$-automorphisms of $F'$. A point $\E$ of the geometric generic fiber corresponds to a pair of subspaces, $\E_1$ of $V_1$ and
$\E_2$ of $V_2$, of dimension $r$ and $s$ respectively, which are the mutual annihilators of each other under the pairing $V_1\times V_2\to F'$ 
induced by the hermitian form $\phi$ on $V$. We therefore obtain a closed embedding 
$$
{ M}_{\{i\}}^{\rm naive}\otimes_{\mathcal O_E}F'\hookrightarrow Gr(r, n)_{F'}\times Gr(s, n)_{F'}\ .
$$
(Note that this is an example of (\ref{confusion}).)
Now
$$
\wedge^{\rm max}\, {\mathcal E}\, =\wedge^{\rm max}\, {\mathcal E}_1\otimes \wedge^{\rm max}\, {\mathcal E}_2 \, .
$$
Hence the pullback of ${\mathfrak L}'$ to the geometric generic fiber is isomorphic to the exterior tensor product of the two ample generators of the Picard groups of
$Gr(r, n)$, resp. $Gr(s, n)$. By identifying ${ M}_{\{i\}}^{\rm naive}\otimes_{\mathcal O_E}F'$ 
with the Grassmannian $Gr(r, n)_{F'}$, this pull back becomes the second tensor power of the ample generator of 
the Picard group of $Gr(r,n)$.  This proves our claim.

This establishes the first part of the theorem.
The second part follows because the inverse image of
$\mathcal A^{I}(\mu)^{\circ}$ under the projection $LSU_n/L^+B \to
LSU_n/L^+P^{I^{\circ}}$ is a union of Schubert varieties, which, by Theorem \ref{normalmain}, are
normal, simultaneously Frobenius split and with rational
singularities. Hence the same
holds also for the irreducible components of $\mathcal
A^{I}(\mu)^{\circ}$.
\end{proof}

\begin{Remark} {\rm For some index sets $I$, one can prove directly the statements in Theorem \ref{theorem11.2}, and then deduce the conjecture \ref{littel} for $I^{\circ}$. We mention the following fact.}

Suppose that $I=\{m\}(=I^{\sharp})$ if $n=2m$ is even and that $I=\{0\}$ if $n$ is odd. Then
the local model $ M_{I}^{\rm loc}$ has reduced irreducible special fiber which
is normal, Frobenius split, and has only rational
singularities. In fact, $\bar{ M}_{I}^{\rm
loc}=\mathcal A^{I}(\mu)$.

{\rm
Indeed, in  [P-R3], \S 5, we show that, in these cases, the
special fiber is irreducible and contains an open subset which is
reduced and equal to an open subset of the Schubert cell corresponding to
the image  in $ W^{I}\backslash
\tilde W \slash W^{I}$ of the element $\lambda$ in $\tilde W$ arising from $\mu$. The underlying reduced
scheme $(\bar{ M}_{I}^{\rm loc})_{\rm red}$ is
therefore the  Schubert variety ${\mathcal A^{I}}(\mu)$ in $LU_n/L^+P^I(U_n)$. We deduce from Theorem
\ref{normalmain} that $(\bar{ M}_{I}^{\rm
loc})_{\rm red}$ is normal, Frobenius split and has only
rational singularities. Now an application of Hironaka's
lemma (EGA IV.5.12.8) to the flat $\mathcal O_E$-scheme
$ M_{I}^{\rm loc}$ shows that
$(\bar{ M}_{I}^{\rm loc})_{\rm red}
=\bar{ M}_{I}^{\rm loc}=\mathcal A^{I}(\mu)$.

We note that, even though for $I$ as above  the parahoric subgroup $P'_{I}$ is a special maximal parahoric subgroup of $SU_n$, the local model ${ M}_{I}^{\rm loc}$ is rarely smooth. 
This is in contrast to what happens for the groups $G=GSp_{2n}$ and $G=GL_n$, cf.~[Go1], [Go2].  More precisely, 
for such $I$, ${ M}_{I}^{\rm loc}$ is smooth 
only in the following cases:  $r=0$, or $r=n$, or $n=2m$ is even and $r=1$ or $r=n-1$.} 
\end{Remark}


\begin{thebibliography} {XYZU}

\bibitem[An] {An} S.\ Anantharaman:   Sch\'emas en groupes, espaces homog\`enes et espaces alg\'ebriques sur une base de dimension 1. Bull. Soc. Math. France, Mem. {\bf 33} (1973), 5--79, Soc. Math. France, Paris.

\bibitem[B-L] {B-L} A.\ Beauville,  Y.\ Laszlo: Conformal blocks and generalized theta functions. Commun. Math. Phys. {\bf 164} (1994), 385--419.

\bibitem[BLS] {BLS} A.\ Beauville, Y.\ Laszlo, C.\ Sorger: The Picard group of the moduli of $G$-bundles on a curve. Compositio Math. {\bf 112} (1998),  183--216.


\bibitem[B-D] {B-D} A.\ Beilinson, V.\ Drinfeld: Quantization of Hitchin's integrable system and Hecke eigensheaves. Preprint, available under http://www.ma.utexas.edu/~benzvi/

\bibitem[Bo] {Bo}  M.\ Borovoi,  
Abelian Galois cohomology of reductive groups. (English summary) 
Mem. Amer. Math. Soc. 132 (1998), no. 626, 

\bibitem[BLR] {BLR} S.\ Bosch, W.\ L\"utkebohmert, M.\ Raynaud: {\sl N\'eron models.} Ergebnisse der Mathematik und ihrer Grenzgebiete (3), {\bf 21}. Springer-Verlag, Berlin (1990).

\bibitem[B-T] {B-T} F.\ Bruhat, J.\ Tits: Groupes r\'eductifs sur un corps local. Inst. Hautes \'Etudes Sci. Publ. Math. {\bf 41} (1972), 5--251.

\bibitem[B-TII] {B-TII} F.\ Bruhat, J.\ Tits: Groupes r\'eductifs sur un corps local. II. Sch\'emas en groupes. Existence d'une donn\'e radicielle valu\'ee. Inst. Hautes \'Etudes Sci. Publ. Math.  {\bf 60} (1984), 197--376.

\bibitem[C] {C} C.-L. Chai:  N\'eron models for semiabelian varieties: congruence and change of base field.  Asian J. Math. {\bf 4} (2000), no. 4, 715--736. 

\bibitem[D-G] {D-G}  M.\ Demazure, P.\ Gabriel: {\sl Groupes alg\'ebriques. Tome I: G\'eom\'etrie alg\'ebrique, g\'en\'eralit\'es, groupes commutatifs}.  Avec un appendice  {\sl Corps de classes local} par Michiel Hazewinkel. Masson \& Cie, \'Editeur, Paris, North-Holland Publishing Co., Amsterdam (1970).

\bibitem[Ed] {Ed} B.\ Edixhoven: N\'eron models and tame ramification. Compositio Math. {\bf 81} (1992), 291--306.

\bibitem[Fa1] {Fa1} G.\ Faltings: Algebraic loop groups and moduli spaces of bundles. J. Eur. Math. Soc. (JEMS) {\bf 5}  (2003), 41--68.

\bibitem[Fa2] {Fa2} G.\ Faltings: Course at  Bonn University, winter semester 2004. 

\bibitem[G] {G} D.\ Gaitsgory: Construction of central elements in the affine Hecke algebra via nearby cycles.  Invent. Math. {\bf 144} (2001), 253--280.

\bibitem[G-Y] {GY} W.-T.\ Gan,   J.-K.\ Yu:  
Sch\'emas en groupes et immeubles des groupes exceptionnels sur un corps local. I. Le groupe $G\sb 2$.  
Bull. Soc. Math. France {\bf 131} (2003), no. 3, 307--358. 


\bibitem[Ga] {Ga} H.\ Garland:  The arithmetic theory of loop algebras. J. Algebra {\bf 53} (1978), no. 2, 480--551.

\bibitem[Go1] {Go1} U.\ G\"ortz:  On the flatness of models of certain Shimura varieties of PEL-type.  Math. Ann.  {\bf 321}  (2001),  689--727. 

\bibitem[Go2] {Go2} U.\ G\"ortz: On the flatness of local models for the symplectic case. Adv. Math. {\bf 176} (2003), 89--115.

\bibitem[Gr] {Gr} A.\ Grothendieck:  Le groupe de Brauer. III. Exemples et compl\'ements.    Dix Expos\'es sur la Cohomologie des Sch\'emas,  pp. 88--188.  North-Holland, Amsterdam; Masson, Paris (1968)

\bibitem[GrII] {GrII} M.\ Greenberg: Schemata over local rings. II.  Ann. of Math. (2) {\bf 78} (1963), 256--266.

\bibitem[H-R] {H-R} T.\ Haines, M.\ Rapoport: On parahoric subgroups. Appendix to this paper. 

\bibitem[Kac] {Kac} V.\ Kac: {\sl Infinite-dimensional Lie algebras.} Third edition. Cambridge University Press, Cambridge (1990).

\bibitem[K] {K} C.\ Kaiser: Ein getwistetes fundamentales Lemma f\"ur die ${\rm GSp}\sb 4$.  Dissertation, Rheinische Friedrich-Wilhelms-Universit\"at Bonn, Bonn, 1997. Bonner Mathematische Schriften, {\bf 303}. Universit\"at Bonn, Mathematisches Institut, Bonn (1997). 

\bibitem[Ko] {Ko} R.\ Kottwitz: Isocrystals with additional structure. II.
Compositio Math. {\bf 109} (1997), 255--339.

\bibitem[Kon] {Kon}  B.\ Kostant:  Groups over $\Z$.  Algebraic Groups and Discontinuous Subgroups (Proc. Sympos. Pure Math. {\bf 9}, Boulder, Colo., 1965) pp. 90--98. Amer. Math. Soc.,  Providence, R.I. (1966)

\bibitem[Ku1] {Ku1} S.\ Kumar: Demazure character formula in arbitrary Kac-Moody setting.
Invent. Math. {\bf 89} (1987), 395--423.

\bibitem[Ku2] {Ku2} S.\ Kumar: {\sl Kac-Moody groups, their flag varieties and representation theory.} Progress in Mathematics {\bf 204}, Birkh\"auser Boston Inc., Boston, MA (2002).

\bibitem[LS] {LS} Y.\ Laszlo, C.\ Sorger: The line bundles on the moduli of parabolic $G$-bundles over curves and their sections. Ann. Sci. Ecole Norm. Sup. (4) {\bf 30} (1997), 499--525.

\bibitem[Li] {Li} P.\ Littelmann: Contracting modules and standard monomial theory for symmetrizable Kac-Moody algebras. J. Amer. Math. Soc. {\bf 11} (1998), 551--567.

\bibitem[Lu] {Lu} G.\ Lusztig: Study of a $\Z$-form of the coordinate ring of a reductive group. arXiv:0709.1286 
 



\bibitem[Ma1] {Ma1} O.\ Mathieu: Formules de caract\`eres pour les alg\`ebres de Kac-Moody g\'en\'erales.
 Ast\'erisque {\bf 159-160} (1988).
 
 
 \bibitem[Ma2] {Ma2} O.\ Mathieu:  Construction d'un groupe de Kac-Moody et applications.  Compositio Math. {\bf 69} (1989), no. 1, 37--60. 
 
 \bibitem[M-K]{M-K} V.\ Mehta, W.\  van der Kallen: On a Grauert-Riemenschneider
 vanishing theorem for Frobenius split varieties in characteristic $p$, Invent. Math.
 {\bf 108} (1992), 11-13. 
 
 \bibitem[MS] {MS} J.\ Milne,   K.-y.\  Shih:
Conjugates of Shimura varieties.
Hodge cycles, motives, and Shimura varieties, Lect. Notes in Math. {\bf 900}, Springer-Verlag, Berlin (1982), 280- 356.
 
 
\bibitem[P-R1] {P-R1} G.\ Pappas, M.\ Rapoport: Local models in the ramified case I. The EL-case. J. Algebraic Geom.  {\bf 12} (2003), 107--145.

\bibitem[P-R2] {P-R2} G.\ Pappas, M.\ Rapoport: Local models in the ramified case II. Splitting models.
Duke Math. Journal {\bf 127} (2005), 193--250.

\bibitem[P-R3] {P-R3} G.\ Pappas, M.\ Rapoport: Local models in the ramified case III. Unitary groups. Preprint math.AG/0702286  

\bibitem[P-S] {P-S} A.\ Pressley, G.\  Segal: {\it Loop groups}. Oxford Mathematical Monographs. Oxford Science Publications. The Clarendon Press, Oxford University Press, New York (1986) 

\bibitem[P-Y] {P-Y} G. Prasad, J.-K. Yu: On finite group actions on reductive groups and buildings. Invent. Math. {\bf 147} (2002),  545--560. 

\bibitem[R] {R} M.\ Rapoport: A guide to the reduction modulo $p$ of Shimura varieties.  Ast\'erisque {\bf 298} (2005), 271--318.

\bibitem[R-Z] {R-Z} M.\ Rapoport, Th.\ Zink: {\sl Period spaces for $p$--divisible groups}.
Ann.\ of Math. Studies {\bf 141}, Princeton University Press (1996).

\bibitem[Se1] {Se1} J.-P.\ Serre: Groupes de Grothendieck des sch\'emas en groupes r\'eductifs d\'eploy\'es.   Inst. Hautes \'Etudes Sci. Publ. Math. {\bf 34} (1968), 37--52.


\bibitem[Se2] {Se2} J.-P.\ Serre: {\sl Galois cohomology}. Corrected reprint of the 1997 English edition. Springer Monographs in Mathematics. Springer-Verlag, Berlin, 2002.

\bibitem[St] {St} R.\ Steinberg: Generators, relations and coverings of algebraic groups. J. of Algebra {\bf 71} (1981), 527--543.

\bibitem[St2] {St2} R.\ Steinberg:  
Regular elements of semi-simple algebraic groups. 
Inst. Hautes \'Etudes Sci. Publ. Math.  {\bf 25} (1965), p. 49-80 


\bibitem[T1] {T} J.\ Tits: Reductive groups over local fields.  Automorphic forms, representations and $L$-functions. (Proc. Sympos. Pure Math., Oregon State Univ., Corvallis, Ore., 1977), Part 1,  pp. 29--69, Proc. Sympos. Pure Math., XXXIII, Amer. Math. Soc., Providence, R.I. (1979).

\bibitem[T2] {T2} J.\ Tits: Groups and group functors attached to Kac-Moody data. Workshop Bonn 1984 (Bonn, 1984). Lecture Notes in Math. {\bf 1111}, Springer-Verlag, Berlin (1985), 193--223.

\bibitem[T3] {T3} J.\ Tits:  Groupes associ\`es aux alg\'ebres de Kac-Moody.   S\'eminaire Bourbaki, Vol. 1988/89. Ast\'erisque No. {\bf 177}-{\bf 178} (1989), Exp. No. 700, 7--31. 

\bibitem[Yu] {Yu} J.-K.\ Yu: Smooth models associated to concave functions in Bruhat-Tits theory, Preprint (2002).
\end{thebibliography}
\end{document}